% leray

% actadefs
%  simulate the appearance of Acta Mathematica 1934  
%  for the translation of Leray's paper

\nopagenumbers
\magnification=\magstep1
\font\ninerm=cmr9

\def\oddheader#1{\ninerm\hskip2.5cm 
           On the motion of a viscous liquid filling space. \hskip3.3cm #1
           \bigskip\rm}
\def\foddheader#1{\ninerm\hskip2cm 
     Sur le mouvement d'un liquide visqueux emplissant l'espace \hskip2cm #1
           \bigskip\rm}
\def\evenheader#1{\ninerm\noindent #1 \hskip5cm Jean Leray.
           \bigskip\rm}

\def\footrule{\medskip\hrule\smallskip\ninerm}

\def\tripl#1{\int\!\!\!\int\!\!\!\int_{#1}}
\def\intpi{\tripl{\Pi}}
\def\intpipr{\tripl{\Pi'}}
\def\dee#1{\, \delta #1}

\def\partm{\partial^m u_i(x,t)\over
               \partial x_1^h \partial x_2^k \partial x_3^l}

% cover page

.
\vskip50mm
\centerline{\tt A translation of the 1934 paper of Leray}

\vskip20mm
\centerline{\tt August, 2000}

\bigskip
\centerline{\tt Robert E. Terrell}

\vskip20mm
{\ninerm Thanks to Danny Goodman, Shigeru Masuda and Chris Orum
 for helpful comments.}

\vfill
\eject

% f193

\vskip5cm
\centerline{SUR LE MOUVEMENT D'UN LIQUIDE VISQUEUX}
\centerline{EMPLISSANT L'ESPACE$^1$}
\medskip
\ninerm
\centerline{Par}
\medskip
\centerline{JEAN LERAY}
\smallskip
\centerline{\`a Rennes}
\bigskip
\rm
\centerline{{Introduction.}$^2$}
\bigskip

  I. {\it La th\'eorie de la viscosit\'e} conduit \`a admettre que les 
mouvements des liquides visqueux sont r\'egis par les \'equations de
Navier; il est n\'ecessaire de justifier a posteriori cette hypoth\`ese
en \'etablissant {\it le th\`eoreme d'existence suivant:} il exist une
solution des \`equations de Navier qui correspond \`a un \'etat de vitesse
donn\'e arbitrairement \`a l'instant initial. C'est ce qu'a cherch\'e \`a
d\'emontrer M. Oseen$^3$; il n'a r\'eussi \`a \'etablir l'existence d'une 
telle solution que pour une dur\'ee tr\'es br\`eve succ\'edant \`a l'instant
initial. On peut v\'erifier en outre que l'\'energie cin\'etique totale du
liquide reste born\'ee$^4$; mais il ne semble pas possible de d\'eduire de
ce fait que le mouvement lui-m\^eme reste r\'egulier; j'ai m\^eme indiqu\'e
une raison qui me fait croire \`a l'existence de mouvements devenant 
irr\'eguliers au bout d'un temps fini$^5$; je n'ai malheureusement pas
r\'eussi \`a forger un exemple d'une telle singularit\'e.

\footrule

$^1$ Ce m\'emoire a \'et\'e r\'esum\'e dans une note parue aux
Comptes rendus de l'Academie des Sciences, le 20 f\'evrier 1933,
 T. 196, p. 527.

$^2$ Les pages 59--63 de ma These (Journ. de Math. 12, 1933) announcent ce
m\'emoire et en completent l'introduction.

$^3$ Voir Hydrodynamik (Leipzig, 1927), \S 7, p. 66. Acta mathematica T. 34.
Arkiv f\"or matematik, astronomi och fysik. Bd. 6, 1910. Nova acta reg. soc. 
scient. Upsaliensis Ser. IV, Vol. 4, 1917.

$^4$ l. c. 2, p. 59--60.

$^5$ l. c. 2, p. 60--61. Je reviens sur ce sujet an \S 20 du pr\'esent
travail (p. 224).

\smallskip
\hrule
\smallskip
\noindent {\it reset from:}
     25--34198. {\it Acta mathematica.} 63. Imprim\'e le 5 juillet 1934.
\rm

\vfill
\eject
 % e193

\vskip5cm
\centerline{ON THE MOTION OF A VISCOUS LIQUID}
\centerline{FILLING SPACE$^1$}
\medskip
\ninerm
\centerline{by}
\medskip
\centerline{JEAN LERAY}
\smallskip
\centerline{in Rennes}
\rm
\bigskip
\centerline{{Introduction.}$^2$}
\bigskip

  I. {\it The theory of viscosity} leads one to allow that motions of a 
viscous liquid are governed by Navier's equations. It is necessary to
justify this hypothesis a posteriori by establishing {\it the following
existence theorem:} there is a solution of Navier's equations which 
corresponds to a state of velocity given arbitrarily at an initial
instant. That is what Oseen tried to prove$^3$. He only succeeded in
establishing the existence of such a solution for a possibly very short
time after the initial instant. One can also verify that the
total kinetic energy of the liquid remains bounded$^4$ but it does not 
seem possible to deduce from this fact that the motion itself remains
regular$^*$. I have
indicated a reason which makes me believe there are motions which become
irregular in a finite time$^5$. Unfortunately I have not succeeded in
creating an example of such a singularity.

\footrule

$^1$ This paper has been summarized in a note which appeared in
Comptes rendus de l'Academie des Sciences, February 20 1933, vol. 196, p. 527.

$^2$ Pages 59--63 of my Thesis (Journ. de Math. 12, 1933) announce this
paper and complement this introduction.

$^3$ See Hydrodynamik (Leipzig, 1927), \S 7, p. 66. Acta mathematica vol. 34.
Arkiv f\"or matematik, astronomi och fysik. Bd. 6, 1910. Nova acta reg. soc. 
scient. Upsaliensis Ser. IV, Vol. 4, 1917.

$^4$ l. c. 2, p. 59--60.

$^5$ l. c. 2, p. 60--61. I return to this subject in \S 20 of the present
work (p. 224).

\smallskip
\hrule
\smallskip

\noindent {\it reset from:} {\it Acta mathematica.} 63. Printed July 5, 1934.

$^*$ {\tt Translator's note:} ``regular'' is defined on p. 217.

\rm

\vfill
\eject
 % f194

\evenheader{194}

   Il n'est pas paradoxal de supposer en effet que la cause qui r\'egularise le
mouvement --- la dissipation de l'\'energie --- ne suffise pas \`a maintenir
born\'ees et continues les d\'eriv\'ees secondes des conposantes de la vitesse
par rapport aux coordonn\'ees; or la th\'eorie de Navier suppose ces
d\'eriv\'ees secondes born\'ees et continues; M. Oseen lui-m\^eme a d\'ej\`a
insist\'e sur le caract\`ere peu naturel de cette hypoth\'ese; il a 
montr\'e en m\^eme temps comment le fait que le mouvement ob\'eit aux lois
de la m\'ecanique peut s'exprimer \`a l'aide d'\'equations int\'egro-
diff\'erentielles$^1$, o\'u figurent seulement les composantes de la vitesse
et leurs d\'eriv\'ees premi\'eres par rapport aux coordonn\'ees spatiales. Au
cours du pr\'esent travail je consid\`ere justement un syst\'eme de 
relations$^2$ qui \'equivalent aux \'equations int\'egro-diff\'erentielles de
M. Oseen, compl\'et\'es par une in\'egalit\'e exprimant la dissipation de
l'\'energie. Ces relations se d\'eduisant d'ailleurs des \'equations de
Navier \`a l'aide d'int\'egrations par parties qui font dispara\^itre les
d\'eriv\'ees d'ordres les plus \'elev\'es. Et, si je n'ai pu r\'eussir \`a
etablir le th\'eor\`eme d'existence \'enonc\'e plus haut, j'ai pu
n\'eammoins d\'emontrer le suivant$^3$: les relations en question poss\`edent
toujours {\it au moins une solution} qui correspond \`a un \'etat de vitesse
donn\'e initialement et {\it qui est d\'efinie pour une dur\'ee illimit\'e},
dont l'origine est l'instant initial. Peut-\^etre cette solution est-elle
trop peu r\'eguli\`ere pur poss\'eder \`a tout instant des d\'eriv\'ees
seconds born\'ees; alors elle n'est pas, au sens propre du terme, une
solution des \'equations de Navier; je propose de dire qu'elle en constitue
{\it $<<$une solution turbulent$>>$}.$^4$

  Il est d'ailleurs bien remarquable que chaque solution turbulente satisfait
effectivement les \'equations de Navier proprement dites, sauf \`a 
certaines \'epoques d'irr\'egularit\'e; ces \'epoques constituent un ensemble
ferm\'e de mesure nulle; \`a ces \'epoques sont seules v\'erifi\'ees 
certaines conditions de continuit\'e extr\^emement larges.

\footrule

$^1$ Oseen, Hydrodynamik, \S 6, \'equation (1).

$^2$ Voir relations (5.15), p. 240.

$^3$ Voir p. 241.

$^4$ Je me permets de citer le passage suivant de M. Oseen (Hydrodynamik): 
$<<$ A un autre point de vue encore il semble valoir la peine de soumettre \`a
une \'etude attentive les singularit\'es du mouvement d'un liquide visqueux. 
S'il peut surgir des singularit\'es, il nous faut manifestement distinguer
deux esp\`eces de mounvements d'un liquide visqueux. les mouvements
r\'eguliers, c'est-\`a-dire les mouvements sans singularit\'e, et les
mouvements irr\'eguliers, c'est-\`a-dire les mouvements avec singularit\'e.
Or on distingue d'autre part en Hydraulique deux sortes de mouvements, les
mouvements laminaires et les mouvements turbulents. On est des lors tent\'e
de pr\'esumer que les mouvements laminaires fournis par les exp\'eriences
sont identiques aux mouvements r\'eguliers th\'eoriques et que les 
mouvements turbulents exp\'erimentaux s'identifient aux mouvements
irr\'eguliers th\'eoriques. Cette pr\'esumption r\'epond-elle \`a la 
r\'ealit\'e? Seules des recherches ult\'erieures pourront en d\'ecider$>>$.

\rm

\vfill
\eject

 % e194

\evenheader{194}

  In fact it is not paradoxical to suppose that the thing which regularizes the
motion--dissipation of energy--does not suffice to keep the second derivatives
of the velocity components bounded and continuous. Navier's theory assumes
the second derivatives bounded and continuous. Oseen himself had already 
emphasised 
that this was not a natural hypothesis. He showed at the same time how the 
fact 
that the motion obeys the laws of mechanics could be expressed by means of 
integro-differential equations$^1$ which contain only the velocity components 
and
their first spatial derivatives. In the course of the present work I consider
a system of relations$^2$ equivalent to Oseen's integro-differential
equations complemented by an inequality expressing dissipation of energy. 
Moreover,
these relations may be deduced from Navier's equations, using an integration
by parts which causes the higher order derivatives to disappear. 
And, if I have 
not succeeded in establishing the existence theorem stated above, I have 
nevertheless
proved the following$^3$: the relations in question always have {\it at least 
one}
solution corresponding to a given initial velocity and {\it which is defined 
for
an unlimited time} of which the origin is the initial instant. Perhaps that 
solution
is not sufficiently regular to have bounded second partial derivatives at each
instant, so it is not, in a proper sense of the term, a solution to Navier's
equations. I propose to say that it constitutes ``{\it a turbulent solution}''.

  It is moreover quite remarkable that each turbulent solution actually 
satisfies Navier's equations, properly said, except at certain times of 
irregularity. These times constitute a closed$^*$ set of measure zero. At these
times alone must the continuity of the solution be interpreted in a very
generous sense.

\footrule

$^1$ Oseen, Hydrodynamik, \S 6, equation (1).

$^2$ See relations (5.15), p. 240.

$^3$ See p. 241.

$^4$ I allow myself to cite a passage from Oseen (Hydrodynamik): ``From still 
another
point of view it seems worth the trouble to subject the singularities of the 
motion
of a viscous liquid to careful study. If singularities appear, then we must 
distinguish
two types of motion of a viscous liquid, regular motion, which is to say 
motion
without singularity, and irregular motion, which is to say motion with 
singularity.
Now in other parts of Hydraulics one distinguishes two sorts of motion, 
laminar and
turbulent. One is tempted from now on to presume that laminar motions 
furnished
by experiment are identical to theoretical regular motions, and that 
experimental 
turbulent motions are identified with irregular theoretical motion. Does this 
presumption
correspond with reality? Only further research will be able to decide.''

$^*$ [{\tt translator's note:} The set is compact, as is proved on p. 246.]

\rm

\vfill
\eject

 % f195

\foddheader{195}

\noindent  Une {\it solution turbulent} a donc {\it la structure}
suivante: elle se compose {\it d'une succession de solutions r\'egulieres.}

  Si j'avais r\'eussi \`a construire des solutions des \'equations de Navier
qui deviennent irr\'egulieres, j'aurais le droit$^1$ d'affirmer qu'il existe
effectivement des solutions turbulent ne se r\'eduisant pas,tout simplement,
\`a des solutions r\'egulieres. M\^eme si cette proposition \'etait fausse,
la notion de solution turbulente, qui n'aurait d\'es lors plus \`a jouer aucun 
r\^ole dans l'\'etude des liquides visqueux, ne perdrait pas son int\'eret:
il doit bien se pr\'esenter {\it des problems} de Physique math\'ematique
pour lesquels {\it les causes physiques de r\'egularit\'e } ne suffisent
pas \`a justifier {\it les hypoth\'eses de r\'egularit\'e faites lors de la 
mise en \'equation:} \`a ces probl\'ems peuvent alors s'appliquer des
consid\'erations semblables \`a celles que j'expose ici.

  Signalons enfin les deux faits suivants:

  Rien ne permet d'affirmer l'unicit\'e de la solution turbulente qui
correspond \`a un \'etat initial donn\'e. (Voir toutefois Compl\'ements $1^o$,
p. 245; \S 33).  

  La solution qui correspond \`a un \'etat initial suffisamment voisin du 
repos ne devient jamais irr\'eguli\`ere. (Voir les cas de r\'egularit\'e que
signalent les \S 21 et 22, p. 226 et 227).

\medskip

  II.  Le travail pr\'esent concerne les liquides visqueux illimit\'es. Les
conclusions en sont extr\'emement analogues \`a celles d'un autre 
m\'emoire$^2$ que j'ai consacr\'e aux mouvements plans des liquides visqueux
enferm\'es dans des parois fixes convexes; ceci autorise \`a croire que ces
conclusions s'\'etendent au cas g\'en\'eral d'un liquide visqueux \`a deux
ou trois dimensions que limitent des parois quelconques (m\^eme variables).

 L'absense de parois introduit certes quelques complications concernant
l'allure \`a l'infini des fonctions inconnues$^3$, mais simplifie beaucoup
l'expos\'e et met mieux en lumi\`ere les difficult\'es essentielles; le
r\^ole important que joue l'homog\'en\'eit\'e des formules est plus 
\'evident; (les \'equations aux dimensions permettent de pr\'evoir a priori
presque toutes les in\'egalit\'es que nous \'ecrirons).

\footrule

$^1$ En vertu du th\'eorem d'existence du \S 31 (p. 241) et du th\'eoreme
d'unicit\'e du \S 18 (p. 222).

$^2$ Journal de Math\'ematiques, T. 13, 1934.

$^3$ Les conditions \`a l'infini par lesquelles nous caract\'erisons celles
des \'equations de Navier que nous nommons r\'egulieres different 
essentiellement des conditions qu'emploie M. Oseen.
\rm

\vfill
\eject

 % 195

\oddheader{195}

\noindent A {\it turbulent solution} therefore has the following 
{\it structure}:
it is composed of a {\it succession of regular solutions}.

  If I succeed in constructing solutions to Navier's equations which 
become irregular, then I can say that there exist
turbulent solutions which do not simply reduce to regular solutions. Likewise if
this proposition is false, the notion of turbulent solution, which from then on
plays no role in the study of viscous liquids, still does not lose interest: it
serves to present {\it problems} of mathematical Physics for which {\it the 
physical cause of regularity} does not suffice to justify {\it the hypothesis 
of regularity made at the time of writing the equation}; to these problems then
one can apply considerations similar to those which I introduce here.

  Finally let us point out the two following facts:

  Nothing allows one to assert the uniqueness 
of a turbulent solution which 
corresponds to a given initial state. (See however Supplementary information
 $1^o$, p. 245; 
\S 33).

  A solution which corresponds to an initial state sufficiently near rest 
never becomes irregular. (See the case of regularity pointed out in \S 21 and 
\S 22, 
p. 226 and 227.)

\medskip

  II. The present work concerns unlimited viscous liquids. The conclusions are
quite analogous to those of another paper$^2$ that I devoted to plane motion
of a viscous liquid enclosed within fixed convex walls; this leads to the belief
that its conclusions extend to the general case of a viscous liquid in two or 
three
dimensions bounded by arbitrary walls (same variables).

   The absence of walls indeed introduces some complications concerning 
the unknown behavior of functions at infinity$^3$ but greatly 
simplifies the exposition and
brings the essential difficulties more to light. The important role played by
homogeneity of the formulas is more evident (the equations in dimensions allow
us to predict a priori nearly all the inequalities that we write).

\footrule

$^1$ In virtue of the existence theorem of \S 31 (p. 241) and of the 
uniqueness theorem of \S 18 (p. 222).

$^2$ Journal de Math\'ematiques, T. 13, 1934.

$^3$ The conditions at infinity by which we characterize those solutions 
of Navier's equations
which we call regular, differ essentially from those of Oseen.
\rm

\vfill
\eject

% f196

\evenheader{196}

  Rapellons que nous avons d\'ej\`a trait\'e le cas des mouvements plans
illimit\'es$^1$: il est asswz sp\'ecial$^2$; la r\'egularit\'e du mouvement
est alors assur\'ee.

\smallskip

{\it Sommaire du m\'emoire.}

  Le chapitre I  rappelle au Lecteur une s\'erie de propositions d'Analyse,
qui sont importantes, mais qu'on ne peut pas toutes consid\'erer comme 
classiques.

  Le chapitre II \'etablit diverses in\'egalit\'es pr\'eliminaires, ais\'ement
d\'eduites des propri\'etes que poss\`ede la solutions fondamentale 
de M. Oseen.

  Le chapitre III applique ces in\'egalit\'es \`a l'\'etude des solutions 
r\'eguli\`eres des \'equations de Navier.

  Le chapitre IV \'enonce diverses propri\'etes des solutions r\'eguli\`eres,
dont fera usage le chapitr VI.

  Le chapitre V \'etablit qu'\`a tout \'etat initial correspond au moins une
solution turbulente, qui est d\'efinie pendant une dur\'ee illimit\'ee. La
d\'emonstration de ce {\it th\'eorem d'existence} repose sure le {\it principe}
suivant: On n'aborde pas directement le probl\`eme pos\'e, qui est de 
r\'esoudre les equations de Navier; mais on traite d'abord un probl\`eme
voisin dont on peut s'assurer qu'il admet toujours une solution r\'eguli\`ere,
d\'efinie pendant une dur\'ee illimit\'ee; on fait tendre ce probl\`eme
voisin vers le probl\`eme pos\'e et l'on construit la limite 
(ou les limites) de sa solution. Il existe bien une fa\c con elementaire
d'appliquer ce principe: c'est celle qu'utilise mon \'etude des mouvements 
plans des liquides visqueux limit\'es par des parois; mais elle est 
intimement li\'ee \`a cette structure des solutions turbulentes que nous
avons pr\'ec\'edemment signal\'ee; elle ne s'appliquerait pas si cette
structure n'\'etait pas assur\'ee. Nous proc\`ederons ici d'une autre fa\c con,
dont la port\'ee est vraisemblablement plus grande, qui jutifie mieux la
notion de solution turbulente, mais qui fait appel \`a quelques th\'eormes
peu usuels cit\'es au chapitre I.

  Le chapitre VI \'etude la structure des solutions turbulentes.

\footrule

  $^1$ Thesis, Journal de Math\`ematiques 12, 1933; chapitre IV p. 64-82. 
(On peut donner une variante int\'eressante au proc\'ed\'e que nous y
employons en utilisant la notion d'\'etat initial semir\'egulier qu'introduit
le m\'emoire pr\'esent.)

  $^2$  On peut dans ce cs baser l'\'etude du probl\`eme sure la propri\'et\'e
que poss\`ede alors le maximum du tourbillon \`a un instant donn\'e d'\^etre 
une fonction d\'ecroissante du temps.  (Voir: Comptes rendus de l'Acad\'emie 
des Sciences, T. 194; p. 1893; 30 mai 1932). -- M. Wolibner a lui aussi
fait cette remarque.
\rm

\vfill
\eject

 % e196

\evenheader{196}

  Recall that we have already treated the case of unlimited motions in
the plane$^1$. These are special$^2$ because the motion is regular.

\smallskip

{\it Summary of the paper.}

  Chapter I recalls a series of propositions of analysis
which are important, but which are not entirely classical.

  Chapter II establishes several preliminary inequalities easily deduced 
from
properties of Oseen's fundamental solution.

  Chapter III applies the inequalities to the study of regular solutions
of Navier's equation.

  Chapter IV states several properties of regular solutions to be used in 
Chapter VI.

  Chapter V establishes that for each initial state, there is at least one 
turbulent solution defined for unlimited time. The proof of this 
{\it existence theorem} is based
on the following {\it principle}: one doesn't directly approach the problem of 
solving Navier's equations; instead one treats a neighboring problem which can
be proved to have a regular solution of unlimited duration; we let the 
neighboring problem tend toward the original problem to construct the limit 
(or limits) of these solutions. There is an elementary method to apply this
principle: it is the same one which I used in my study of planar motion of a 
viscous liquid within walls, but it
is intimately bound with the structure of turbulent solutions which we have 
previously pointed out. Without this structure the method may not apply. 
Here we proceed in another
fashion whose range is very likely larger, and which justifies the notion of 
turbulent solution, but which requires calling on some of the less ordinary 
theorems of chapter I.

  Chapter VI studies the structure of turbulent solutions.

\footrule

  $^1$ Thesis, Journal de Math\`ematiques 12, 1933; chapter IV p. 64-82. 
(One can give
an interesting variation on the process used there, by using the notion 
of semi-regular
initial state introduced in the present paper.)

  $^2$ In this case one can base the study on the property that the 
maximum swirl is a
decreasing function of time. (See: Comptes rendus de l'Acad\'emie des Sciences, 
T. 194; 
p. 1893; 30 mai 1932). -- Wolibner has also made this remark.

 [{\tt translator's note:\ } 
Wolibner, "Un Theor\`eme sur l'existence du mouvement plan d'un fluide
parfait, homog\`ene, incompressible, pendant un temps infin\`ement long",
Math. Z. 37 (1933), 698-726, (footnote: mark (**) in p.698)
\rm

\vfill
\eject

 % f197

\foddheader{197}

\centerline{\bf I. Pr\'eliminaires}

  {\bf 1.} {\it Notations}

\smallskip
  Nous utiliserons la lettre $\Pi'$ pour d\'esigner un domaine arbitraire
de points de l'espace; $\Pi'$ pourra \^etre l'espace tout entier, que nous
d\'esignerons par $\Pi$; $\varpi$ d\'esignera un domaine born\'e de points
de $\Pi$, dont la fronti\`ere constitue une surface r\'eguli\`ere $\sigma$. 

  Nous repr\'esenterons un point arbitraire de $\Pi$ par $x$, 
ses cooronn\'ees cart\'esiennes par $x_i$ ($i=1$,2,3) sa distance \`a
l'origine $r_0$,
un element de volume qu'il engendre par $\delta x$, un element de surface
qu'il engendre par $\delta x_1$, $\delta x_2$, $\delta x_3$.
Nous d\'esignerons de m\^eme par $y$ un second point arbitraire de $\Pi$;
$r$ repr\'esentera toujours la distance des points nomm\'es $x$ et $y$.

 Nous utileserons la convention $<<$ de l'indice muet$>>$: un term o\`u un
indice figure deux fois repr\'esentera la somme des termes obteunus en 
donnant \`a cet indice successivement les valeurs 1, 2, 3.

 A partir du chapitre II le symbole $A$ nous servira \`a d\'esigner les 
constantes dont nous ne pr\'eciserons pas la valeur num\'erique.

 Nous repr\'esenterons syst\'ematiquement par de grandes lettres les
fonctions que nous supposerons seulement mesurables; par de petites lettres
les fonctions qui sont continues ainsi que leurs d\'eriv\'ees premi\`eres.

  {\bf 2.} {\it Rappelon l'in\'egalit\'e de Schwarz:}

     $$ \left[\intpipr U(x)V(x)\dee{x} \right]^2 \le
                       \intpipr U^2(x) \dee{x} \times \intpipr V^2(x) \dee{x}
     \leqno(1.1)$$

-- On est assur\'e que le premier membre a un sens quand le second est fini. --

Cette in\'egalit\'e est \`a la base de toutes les propri\'et\'es \'enonc\'ees
au cours de ce chapitre.

\medskip
  {\it Premi\`ere application:}

\smallskip
\noindent Si 
           $$ U(x) = V_1(x)+V_2(x) $$
on a:
      $$ \sqrt{\intpipr U^2(x) \dee{x}} \le 
                  \sqrt{\intpipr V_1^2(x) \dee{x}} +
                           \sqrt{\intpipr V_2^2(x) \dee{x}};$$

\vfill
\eject

 % e197

\oddheader{197}

\centerline{\bf I. Preliminaries}

  {\bf 1.} {\it Notation}

\smallskip
  We use the letter $\Pi$' for an arbitrary domain of points in space. $\Pi$' 
may be the entire space, denoted by $\Pi$. $\varpi$ will be a bounded domain 
in $\Pi$ which has
as boundary a regular surface $\sigma$. We represent an 
arbitrary point of $\Pi$
by $x$, which has cartesian coordinates $x_i$ ($i=1$, 2, 3) and distance $r_0$ 
to the origin, and generates  
volume element $\delta x$ and surface element $\delta x_1$, $\delta x_2$, 
$\delta x_3$.  Similarly we use $y$ for a second arbitrary point of $\Pi$. 
$r$ will always represent the distance between points named $x$ and $y$.

  We use the ``silent index'' convention: a term containing the same index 
twice 
represents the sum of terms obtained by successively giving that index the 
values 1, 2, 3.

  Beginning with chapter II the symbol $A$ denotes constants for which we do not
specify the numerical value.

  We systematically use capital letters for measurable functions and lower 
case letters
for functions which are continuous and have continuous first partial 
derivatives.

\medskip
  {\bf 2.} {\it Recall the Schwarz inequality:}

     $$ \left[\intpipr U(x)V(x)\dee{x} \right]^2 \le
                       \intpipr U^2(x) \dee{x} \times \intpipr V^2(x) \dee{x}
     \leqno(1.1)$$

-- The left side is defined whenever the right is finite. --

  This inequality is the foundation of all properties stated in this chapter.

  {\it First application:}

\noindent If
           $$ U(x) = V_1(x)+V_2(x) $$
then
      $$ \sqrt{\intpipr U^2(x) \dee{x}} \le 
                  \sqrt{\intpipr V_1^2(x) \dee{x}} +
                           \sqrt{\intpipr V_2^2(x) \dee{x}};$$

\vfill
\eject

 % f198

\evenheader{198}

\noindent plus g\'en\'eralement si l'on a, $t$ \'etant une constante:

      $$ U(x) = \int_0^t V(x,t')\,dt'$$

alors:

     $$ 
             \sqrt{\intpipr U^2(x) \dee{x}} \leq 
                   \int_0^t dt' \sqrt{\intpipr V^2(x,t') \dee{x}}
     \leqno(1.2)$$

les premiers membres de ces in\'egalit\'es \'etant s\^urement finis quand
les seconds membres le sont.

\medskip
  {\it Seconde application:}
\smallskip

  Soit $n$ constantes $\lambda_p$ et $n$ vecteurs constantes $\vec{\alpha}_p$; 
d\'esignons par $x+\vec{\alpha}_p$ le point obtenu en faisant subir \`a  $x$ 
la translation $\vec{\alpha}_p$; 
nous avons:

$$\intpi \left[\sum_{p=1}^{p=n}\lambda_p U(x+\vec{\alpha}_p)\right]^2 \dee{x}
                     < 
                  \left[\sum_{p=1}^{p=n} |\lambda_p | \right]^2 \times
                    \intpi U^2(x) \dee{x};$$

(cette in\'egalit\'e se d\'emontre ais\'ement en d\'eveloppant les deux carr\'es
qui y figurent et en utilisant l'in\'egalit\'e de Schwarz).  On en d\'eduit
la suivante qui nous sera tr\'es utile: Soit une fonction $H(z)$: nous 
d\'esignerons par $H(y-x)$ la fonction que l'on obtient en substituant aux
coordon\'ees  $z_i$ de $z$ les composantes $y_i-x_i$ du vecteur $\vec{xy}$; 
nous avons:

     $$ 
        \intpi\left[\intpi H(y-x)U(y)\dee{y}\right]^2 \dee{x} 
                 <
         \left[\intpi |H(z)|\dee{z}\right]^2 \times \intpi U^2(y) \dee{y};
     \leqno(1.3)$$ 

et on est assur\'e que le premier membre est fini quand les deux int\'egrals
qui figurent au second membre le sont.

\medskip
  {\bf 3.} {\it Forte convergence en moyenne.}$^1$

\smallskip
  D\'efinition: On dit qu'une infinit\'e de fonctions $U^*(x)$ a pour forte
limite en moyenne sur un domaine $\Pi'$ une fonction $U(x)$ quand:

\footrule

  $^1$ Voir: F. Riesz, Untersuchungen \"uber Systeme integrierbarer Funktionen, 
Math. Ann. T. 69 (1910). Delsarte, M\'emorial des Sciences math\'ematiques, fasicule 57,
Les groupes de transformations lin\'eaires dans l'espace de Hilbert.
\rm

\vfill
\eject

 % e198

\evenheader{198}

\noindent more generally for a constant $t$ if one has

      $$ U(x) = \int_0^t V(x,t')\,dt'$$

then

     $$ 
             \sqrt{\intpipr U^2(x) \dee{x}} \leq 
                   \int_0^t dt' \sqrt{\intpipr V^2(x,t') \dee{x}}
     \leqno(1.2)$$

the left sides of these inequalities being finite when the right sides are.

\medskip
  {\it Second application:}
\smallskip

  Consider $n$ constants $\lambda_p$ and $n$ constant vectors $\vec{\alpha}_p$. 
Write $x+\vec{\alpha}_p$ for the translation of $x$ by the 
vector $\vec{\alpha}_p$. We have

$$\intpi \left[\sum_{p=1}^{p=n}\lambda_p U(x+\vec{\alpha}_p)\right]^2 \dee{x}
                     \leq^*
                  \left[\sum_{p=1}^{p=n} |\lambda_p | \right]^2 \times
                    \intpi U^2(x) \dee{x};$$
(this inequality is easily proved by expanding the two squares and using the 
Schwarz inequality).  From it, one deduces the following very useful one.
Let $H(z)$ be a function.
We denote by $H(y-x)$ the function obtained by substituting for coordinates 
$z_i$ of $z$ the components $y_i-x_i$ of the vector $\vec{xy}$. We have

     $$ 
        \intpi\left[\intpi H(y-x)U(y)\dee{y}\right]^2 \dee{x} 
                \leq^* 
         \left[\intpi |H(z)|\dee{z}\right]^2 \times \intpi U^2(y) \dee{y};
     \leqno(1.3)$$ 

and the left side is finite when the two integrals on the right are finite.

\medskip
  {\bf 3.} {\it Strong convergence in mean.}$^1$

\smallskip
  Definition: One says that an infinity of functions $U^*(x)$ have 
function $U(x)$ 
as strong limit in mean on a domain $\Pi'$ when:

\footrule

  $^1$ See: F. Riesz, Untersuchungen \"uber Systeme integrierbarer Funktionen, 
Math. Ann. vol 69 (1910). Delsarte, M\'emorial des Sciences math\'ematiques, fasicule 57,
Les groupes de transformations lin\'eaires dans l'espace de Hilbert.

$^*$ {\tt Translator's note:} The last two inequality signs appeared
as ``$<$'' in the original. 

\rm

\vfill
\eject

% f199

\foddheader{199}

      $$
             \lim \intpipr [U^*(x)-U(x)]^2 \dee{x} = 0. 
      \leqno(1.4)$$

On a alors quelle que soit la fonction  $A(x)$ de carr\'e sommable sur $\Pi'$:

      $$
            \lim \intpipr U^*(x)A(x) \dee{x} = \intpipr U(x)A(x) \dee{x}.
       \leqno(1.5)$$

De (1.4) et (1.5) r\'esulte

      $$
       \lim \intpipr U^{*2}(x) \dee{x} = \intpipr U^2(x) \dee{x}.
       \leqno(1.6)$$

\medskip
  {\it Faible convergence en moyenne:}

\smallskip
  D\'efinition: Une infinit\'e de fonctions $U^*(x)$ a pour faible limite
en moyenne sur un domaine $\Pi'$ une fonction $U(x)$  quand les deux 
conditions suivantes se trouvent r\'ealis\'ees:

\smallskip
  a) les nombres $\intpipr U^{*2}(x) \dee{x} $ sont born\'es dans leur ensemble;

\smallskip
  b) on a quelle que soit la fonction  $A(x)$ de carr\'e sommable sur $\Pi'$:

      $$ \lim \intpipr U^*(x)A(x) \dee{x} = \intpipr U(x)A(x) \dee{x}.$$

\ninerm

\medskip
  {\it Exemple I.} La suite $\sin x_1$, $\sin 2x_1$, $\sin 3x_1$, $\ldots$ 
converge faiblement vers z\'ero sur tout domaine $\varpi$.

\smallskip
  {\it Exemple II.} Soit une infinit\'e de fonctions $U^*(x)$  admettant une
fonction $U(x)$ comme forte limite en moyenne sure tout domaine  $\varpi$,
elle l'admet comme faible limite en moyenne sur $\Pi$  quand les quantit\'es
$\intpipr U^{*2}(x) \dee{x} $ sont bourn\'ees dans leur ensemble.

\smallskip
  {\it Exemple III.} Soit une infinit\'e de fonctions $U^*(x)$ qui sur un
domaine $\Pi'$ convergent presque partout vers une fonction $U(x)$; cette
fonction est leur faible limite en moyenne quand les quantit\'es  
$\intpipr U^{*2}(x) \dee{x} $ sont bourn\'ees dans leur ensemble.

\rm

\medskip
On a:
       $$
          \lim
         \tripl{{\Pi_1}'} \tripl{{\Pi_2}'} A(x,y)U^*(x)V^*(y) \dee{x}\dee{y} = 
         \tripl{{\Pi_1}'} \tripl{{\Pi_2}'} A(x,y)U(x)V(y) \dee{x}\dee{y}
       \leqno(1.7)$$

\vfill
\eject

 % e199

\oddheader{199}

      $$
             \lim \intpipr [U^*(x)-U(x)]^2 \dee{x} = 0. 
      \leqno(1.4)$$

One then has for any square summable function $A(x)$ on $\Pi'$

      $$
            \lim \intpipr U^*(x)A(x) \dee{x} = \intpipr U(x)A(x) \dee{x}.
       \leqno(1.5)$$

From (1.4) and (1.5)

      $$
       \lim \intpipr U^{*2}(x) \dee{x} = \intpipr U^2(x) \dee{x}.
       \leqno(1.6)$$

\medskip
  {\it Weak convergence in mean:}
\smallskip

  Definition: An infinity of functions $U^*(x)$ has function $U(x)$ as weak
limit in mean on domain $\Pi'$ when the two following conditions hold.
\smallskip
  a) the set of numbers $\intpipr U^{*2}(x) \dee{x} $ is bounded;

  b) for all square summable functions $A(x)$ on $\Pi'$

      $$ \lim \intpipr U^*(x)A(x) \dee{x} = \intpipr U(x)A(x) \dee{x}.$$

\medskip
  {\it Example I.} The sequence $\sin x_1$, $\sin 2x_1$, $\sin 3x_1$, $\ldots$ 
converges weakly to zero on all domains $\varpi$.

\smallskip
  {\it Example II.} If an infinity of functions $U^*(x)$ have strong limit 
$U(x)$ in mean on all domains $\varpi$, then they admit a weak limit in 
mean on $\Pi$ when the set of quantities $\intpipr U^{*2}(x) \dee{x} $ is 
bounded.

\smallskip
  {\it Example III.} Let an infinity of functions $U^*(x)$ on a domain $\Pi'$ 
converge almost everywhere to a function $U(x)$. That function is their weak 
limit in mean when the set of quantities $\intpipr U^{*2}(x) \dee{x} $ 
is bounded.

\medskip
  One has

       $$
          \lim
         \tripl{{\Pi_1}'} \tripl{{\Pi_2}'} A(x,y)U^*(x)V^*(y) \dee{x}\dee{y} = 
         \tripl{{\Pi_1}'} \tripl{{\Pi_2}'} A(x,y)U(x)V(y) \dee{x}\dee{y}
       \leqno(1.7)$$

\vfill
\eject

 % f200

\evenheader{200}

\noindent quand on suppose que les $U^*(x)$ convergent faiblement en moyenne
vers  $U(x)$ sur ${\Pi_1}'$, les $V^*(x)$ vers $V(x)$ sur ${\Pi_2}'$ 
et que l'int\'egral

       $$ \tripl{{\Pi_1}'}\tripl{{\Pi_2}'}A^2(x)\dee{x}\dee{y}$$
est fini. 

On a:

       $$ 
             \lim \intpipr A(x)U^*(x)V^*(x) \dee{x} =
             \intpipr A(x)U(x)V(x) \dee{x}
       \leqno(1.8)$$ 
quand on suppose, sur $\Pi'$, $A(x)$ born\'e , $U(x)$ forte limite des $U^*(x)$
et $V(x)$ faible limite des $V^*(x)$.

  Il est d'autre part \'evident que l'on a, si les fonctions $U^*(x)$ 
convergent  faiblement en moyenne vers $U(x)$ sur un domaine ${\Pi_1}'$:

$$ \lim \{\intpipr [U^*(x)-U(x)]^2\dee{x}
              -\intpipr U^{*2}(x)\dee{x}+\intpipr U^2(x)\dee{x}\}=0;$$
d'o\`u r\'esultent l'in\'egalit\'e:

       $$
             \liminf \intpipr U^{*2}(x)\dee{x} \geq \intpipr U^2(x)\dee{x},
       \leqno(1.9)$$
et le {\it crit\`ere de forte convergence:}

  Les fonctions $U^*(x)$ convergent fortement en moyenne sur le domaine $\Pi'$
vers la fonction $U(x)$ quand elles convergent faiblement en moyenne vers cette
fonction sur ce domaine et qu'en outre:

       $$
              \limsup \intpipr U^{*2}(x)\dee{x} \leq \intpipr U^2(x)\dee{x}.
       \leqno(1.10)$$

   De m\^eme: Les composantes $U_i^*(x)$ d'un vecteur convergent fortement en
moyenne sur le domaine $\Pi'$ vers celles d'un vecteur $U_i(x)$ quand elles
convergent faiblement en moyenne vers ces composantes sur ce domaine
et qu'en outre$^1$:

\footrule

$^1$ Rappelons que le symbole $U_i(x)U_i(x)$ repr\'esent l'expression
       $\sum_{i=1}^{i=3}U_i(x)U_i(x)$.
\rm

\vfill
\eject

 % e200

\evenheader{200}

\noindent when the $U^*(x)$ converge weakly in mean to $U(x)$ on ${\Pi_1}'$ and
$V^*(x)$ to $V(x)$ on ${\Pi_2}'$ and the integral

       $$ \tripl{{\Pi_1}'}\tripl{{\Pi_2}'}A^2(x)\dee{x}\dee{y}$$
is finite. One has

       $$ 
             \lim \intpipr A(x)U^*(x)V^*(x) \dee{x} =
             \intpipr A(x)U(x)V(x) \dee{x}
       \leqno(1.8)$$ 
when on $\Pi'$, $A(x)$ is bounded, $U(x)$ is the strong limit of the $U^*(x)$ 
and
$V(x)$ is the weak limit of the $V^*(x)$.

  It is also evident that, if the functions $U^*(x)$ converge weakly in 
mean to $U(x)$ on a domain ${\Pi_1}'$

$$ \lim \{\intpipr [U^*(x)-U(x)]^2\dee{x}
              -\intpipr U^{*2}(x)\dee{x}+\intpipr U^2(x)\dee{x}\}=0$$
from which one gets the inequality

       $$
             \liminf \intpipr U^{*2}(x)\dee{x} \geq \intpipr U^2(x)\dee{x}
       \leqno(1.9)$$
and the {\it criteria for strong convergence:}

  The functions $U^*(x)$ converge strongly in mean on domain $\Pi'$ to 
$U(x)$ when 
they converge weakly in mean to $U(x)$ on the domain and in addition

       $$
              \limsup \intpipr U^{*2}(x)\dee{x} \leq \intpipr U^2(x)\dee{x}.
       \leqno(1.10)$$

   Equivalently, the components $U_i^*(x)$ of the vector converge strongly in 
mean on
domain $\Pi'$ to those of a vector $U_i(x)$ when they converge weakly in mean 
to the
components on the domain and in addition$^1$

\footrule

$^1$ Recall that the symbol $U_i(x)U_i(x)$ represents the expression
       $\sum_{i=1}^{i=3}U_i(x)U_i(x)$.
\rm

\vfill
\eject

 % f201

\foddheader{201}

     $$ \limsup \intpipr U_i^*(x) U_i^*(x) \dee{x}
                \leq
                \intpipr U_i(x) U_i(x) \dee{x}.
     \leqno(1.10')$$

 Ce crit\`ere de faible convergence appliqu\'e \`a l'Exemple III fournit
la propri\'et\'e suivante:

\smallskip
   {\it Lemme 1.} Soit une infinit\'e de fonctions $U^*(x)$ [ou de vecteurs
$U_i^*(x)$] qui converge presque partout sur un domaine $\Pi'$ vers une
fonction $U(x)$
[ou un  vecteur $U_i^*(x)$]; elles [ils] convergent fortement en moyenne
vers cette limite quand l'in\'egalit\'e (1.10) [ou (1.10')] est v\'erifi\'ee.

\medskip
   {\it Th\'eoreme de F. Riesz:} Une infinit\'e de fonctions $U^*(x)$
poss\`ede une faible limite en moyenne sure un domaine $\Pi'$ si les deux
conditiones suivant sont v\'erifi\'ees:

\smallskip
    a) les nombres $\intpipr U^{*2}(x)\dee{x}$ sont born\'ees dans leur 
         ensemble;

\smallskip
    b) pour chaque fonction $A(x)$ de carr\'e sommable sur $\Pi'$ 
    les quantit\'es
       $\tripl{c} U^*(x)A(x) \dee{x}$ ont une seule valeur limite.

\medskip
On peut substituer \`a la condition b) la suivante:

\smallskip
    b') Pour chaque cube $c$ dont les ar\`etes sont parall\`eles aux axes de
coordonn\'es et dont les sommets ont des coordonn\'ees rationnelles les
quantit\'es
       $\intpipr U^*(x) \dee{x}$ ont une seule valeur limite.

\medskip
\ninerm

    La d\'emonstration de ce th\'eoreme fait usage des travaux de
M. Legesgue sure les fonctions sommables.

\rm

\medskip
    {\bf 4.}{\it Proc\'ed\'e diagonal de Cantor.}

\smallskip
   Soit une infinit\'e d\'enombrable de quantit\'es d\'ependant chacune de
l'indice entier
$n$:$a_n$, $b_n$, $\ldots (n=1,2,3 \ldots)$. Supposons les  $a_n$
born\'es dans leur ensemble, les $b_n$ born\'es dans leur ensemble, etc. Le 
proc\'ed\'e de Cantor permet de trouver une suite d'entiers 
$m_1$, $m_2$, $\dots$ telle que chacun des suites 
$a_{m_1}$, $a_{m_2}$, $\ldots$; $b_{m_1}$, $b_{m_2}$, $\ldots$;
$\ldots$ converge vers une limite.

\medskip
\ninerm

   Rappelons br\^evement quel est ce proc\'ed\'e: on construit 
une premi\`ere suite d'entiers
$n_1^1$, $n_2^1$, $n_3^1$ $\ldots$ telle que les quantit\'es  $a_{n_1^1}$,
$a_{n_2^1}$, $a_{n_3^1}$, $\ldots$ convergent vers une limite; 
on constitue avec des \'el\'ements de cette premi\`ere suite une seconde suite
$n_1^2$, $n_2^2$, $n_3^2$ $\ldots$, telle que les quantit\'es $b_{n_1^2}$, 
$b_{n_2^2}$, $b_{n_3^2}$, $\ldots$ 

\rm

\vfill
\eject
 % e201

\oddheader{201}

     $$ \limsup \intpipr U_i^*(x) U_i^*(x) \dee{x}
                \leq
                \intpipr U_i(x) U_i(x) \dee{x}.
     \leqno(1.10')$$

  The weak convergence criteria applied in Example III gives the following.

   {\it Lemma 1.} If an infinity of functions $U^*(x)$ [or vectors
 $U_i^*(x)$] converge almost everywhere on domain $\Pi'$ to a function $U(x)$
[or a vector $U_i^*(x)$] and satisfy inequality (1.10) [or (1.10')],
then they converge strongly in mean.

\medskip
   {\it Theorem of F. Riesz:} An infinity of functions $U^*(x)$ have a weak
limit in mean on domain $\Pi'$ if the two following conditions are satisfied:

\smallskip
    a) the set of numbers $\intpipr U^{*2}(x)\dee{x}$ is bounded;

\smallskip
    b) for each square summable function $A(x)$ on $\Pi'$ the quantities
       $\intpipr U^*(x)A(x) \dee{x}$ have a single limiting value.

\medskip
  Condition b) may be replaced by the following:

\smallskip
    b') for each cube $c$ with sides parallel to the coordinate axes and
   rational vertices,  the quantities
       $\tripl{c} U^*(x) \dee{x}$ have a single limiting value.

\medskip
     \ninerm

     The proof of this theorem makes use of the work of Lebesgue on
summable functions.

     \rm

\medskip
    {\bf 4.}{\it Cantor's diagonal method.}

\smallskip
   Consider a countable infinity of quantities each dependent on integer
indices $n$: \ $a_n$, $b_n$, $\ldots (n=1,2,3 \ldots)$. Suppose the $a_n$
are bounded, the $b_n$ are bounded, etc. Cantor's diagonal method allows us
to find a sequence of integers $m_1$, $m_2$, $\dots$, such that each of the
sequences $a_{m_1}$, $a_{m_2}$, $\ldots$; $b_{m_1}$, $b_{m_2}$, $\ldots$;
$\ldots$ converge to a limit.

\medskip
\ninerm

   Recall this process briefly: one constructs a first sequence of integers
$n_1^1$, $n_2^1$, $n_3^1$ $\ldots$ such that  
the quantities $a_{n_1^1}$, $a_{n_2^1}$, $a_{n_3^1}$, $\ldots$ converge to
a limit; one then constructs using elements of the first sequence a second
sequence $n_1^2$, $n_2^2$, $n_3^2$ $\ldots$, such that the quantities 
$b_{n_1^2}$, $b_{n_2^2}$, $b_{n_3^2}$, $\ldots$ 

\rm

\vfill
\eject

% f202

\evenheader{202}

\ninerm

\noindent converge vers une limite; etc. On choisit alors $m_p$ 
\'egal \`a $n_p^p$, qui est le $p^{ieme}$ terme 
de la diagonale du tableau infini des $n_i^j$.

\rm

\medskip
   Application: Du th\'eoreme cit\'e au parapraphe pr\'ec\'edent r\'esulte
le suivant:

\smallskip
   {\it Th\'eoreme fondamental de M. F. Riesz:} Soit une infinit\'e de
fonctions $U^*(x)$ d\'efinie sur un domaine $\Pi'$ 
et telles que les quantit\'es $\intpipr U^{*2}(x) \dee{x}$
soient born\'ees dans leur ensemble; on peut toujours en extraire une suite
illimit\'e de fonctions poss\'edant une faible limite en moyenne.

\smallskip
\ninerm

   En effet la condition a)  est satisfait et le Proc\'ed\'e diagonal de 
Cantor permet de construire une suite de fonctions $U^*(x)$ qui v\'erifient
la condition b').

\rm

\medskip
{\bf 5.}{\it Divers modes de continuit\'e d'une fonction par rapport
       \`a un param\`etre.}

\smallskip
   Soit une fonction $U(x,t)$ d\'ependant d'un param\`etre $t$.  Nous dirons
qu'elle est 
{\it uniform\'ement continue} en $t$ 
quand les trois conditions suivantes seront r\'ealis\'ees:

\smallskip
   a) elle est continue par rapport \`a $x_1$, $x_2$, $x_3$, $t$;

\smallskip
   b) pour chaque valeur particuli\`ere $t_0$ de $t$ le maximum de $U(x,t_0)$
      est fini; 

\smallskip
   c) \'etant donn\'e un nombre positif $\epsilon$, 
          on peut trouver un nombre positif $\eta$
      tel que l'in\'egalit\'e $|t-t_0|<\eta$ entra\^ine:

                  $$|U(x,t)-U(x,t_0)| < \epsilon.$$
Le maximum de $|U(x,t)|$ sur  $\Pi$ est alors une fontion continue de $t$.

\medskip
  Nous dirons que $U(x,t)$ est {\it fortement continue} en $t$ quand,
pour chaque valeur particuli\`ere $t_0$ de $t$, $\intpi U^2(x,t_0) \dee{x}$
est fini et qu'on peut, \'etant donn\'e $\epsilon$, trouver  $\eta$ tel que
l'in\'egalit\'e $|t-t_0|<\eta$ entra\^ine:

              $$\intpi [U(x,t)-U(x,t_0)]^2 \dee{x} < \epsilon.$$
L'int\'egral $\intpi U^2(x,t_0) \dee{x}$ est donc une fonction continue de
$t$. Inversement le lemme I nous apprend qu'une fonction $U(x,t)$ continue
par rapport aux variables $x_1$, $x_2$, $x_3$, $t$ est fortement continue
en $t$ quand l'int\'egral pr\'ec\'edent est une fonction continue de  $t$.

\vfill
\eject

 % e202

\evenheader{202}

\ninerm

\noindent converge to a limit; etc. One then chooses $m_p$ equal to
$n_p^p$, which is the $p$-th term of the diagonal of the infinite table
of $n_i^j$.

\rm

\medskip
   Application: The following results from the theorem of the preceeding
paragraph.

\smallskip
   {\it Fundamental theorem of F. Riesz} Let an infinity of functions $U^*(x)$
on a domain $\Pi'$ be such that the quantities $\intpipr U^{*2}(x) \dee{x}$
are bounded. Then one can always extract from them a sequence of functions
which have a weak limit in mean.

\smallskip
\ninerm

   In fact condition a) is satisfied and Cantor's diagonal process allows
construction of a sequence of functions $U^*(x)$ which satisfy condition b').

\rm

\medskip
{\bf 5.}{\it Various modes of continuity of a function with respect to
      a parameter.}

\smallskip
   Let a function $U(x,t)$ depend on a parameter $t$. We say it is 
{\it uniformly continuous} in $t$ when the following three conditions
hold:

\smallskip
   a) it is continuous with respect to $x_1$, $x_2$, $x_3$, $t$;

\smallskip
   b) for each particular value $t_0$ of $t$ the maximum of $U(x,t_0)$
       is finite;

\smallskip
   c) given a positive number $\epsilon$, one can find a positive $\eta$
      such that the inequality $|t-t_0|<\eta$ implies

                  $$|U(x,t)-U(x,t_0)| < \epsilon.$$
The maximum of $|U(x,t)|$ on $\Pi$ is then a continuous function of $t$.

\medskip
  We say that $U(x,t)$ is {\it strongly continuous} in $t$ when, for each
particular value $t_0$ of $t$, $\intpi U^2(x,t_0) \dee{x}$ is finite and
for each $\epsilon$ there is $\eta$ such that the inequality $|t-t_0|<\eta$ 
implies

              $$\intpi [U(x,t)-U(x,t_0)]^2 \dee{x} < \epsilon.$$
The integral $\intpi U^2(x,t_0) \dee{x}$ is therefore a continuous function
of $t$. Conversely we learn from lemma 1 that a function $U(x,t)$ continuous
with respect to variables $x_1$, $x_2$, $x_3$, $t$ is strongly continuous
in $t$ when the preceeding integral is a continous function of $t$.

\vfill
\eject

 % f203

\foddheader{203}

  {\bf 6.} {\it Relations entre une fonction et ses d\'eriv\'ees}

\smallskip

   Consid\'erons deux fonctions $u(x)$ et $a(x)$ poss\'edent des d\'eriv\'ees
premi\`eres continues qui soient, comme ces fonctions elles-m\^emes, de
carr\'es sommables sur $\Pi$. $s$ \'etant la surface d'une sph\`ere $S$
dont le centre est l'origine et dont le rayon $r_0$ augmente ind\'efiniment,
posons:

    $$\varphi(r_0) = \int\!\!\int_s u(x)a(x) \dee{x_i};$$
nous avons:

    $$\varphi(r_0) = \tripl{S}\left[u(y){\partial a(y) \over \partial y_i}
             + {\partial u(y) \over \partial y_i}a(y)\right]\dee{y}.$$
La second expression de $\varphi(r_0)$ prouve que 
cette quantit\'e tend vers une
limite $\varphi(\infty)$ quand $r_0$ augmente ind\'efiniment. La premi\`ere
expression de $\varphi(r_0)$ nous donne:

    $$|\varphi(r_0)| \leq \int\!\!\int_s |u(x)a(x)| {x_i\dee{x_i} \over r_0}$$
d'o\`u:

 $$\int_0^\infty |\varphi(r_0)|\,dr_0
            \leq 
         \intpi |u(x)a(x)| \dee{x}.$$
Par suite $\varphi(\infty)=0$; en d'autres termes:

  $$\intpi\left[u(y){\partial a(y) \over \partial y_i}
             + {\partial u(y) \over \partial y_i}a(y)\right]\dee{y}
             = 0;
  \leqno(1.11)$$
il en r\'esulte que plus g\'en\'eralement:

  $$\tripl{\Pi-\varpi}\left[u(y){\partial a(y) \over \partial y_i}
             + {\partial u(y) \over \partial y_i}a(y)\right]\dee{y}
             = -\int\!\!\int_\sigma u(y)a(y) \dee{y_i}.
  \leqno(1.12)$$

   Choisissons comme domaine $\varpi$ une sph\`ere de rayon infiniment petit
dont nous nommerons le centre $x$; faisons$^1$ dans (1.12) 
$a(y) = {1 \over 4\pi}{\partial({1 \over r}) \over \partial y_i}$;
ajoutons les

\footrule
 
  $^1$ r repr\'esente la distance des points $x$ et $y$.

\vfill
\eject

 % e203

\oddheader{203}

  {\bf 6.} {\it Relations between a function and its derivatives}

\smallskip

   Consider two functions 
    $u(x)$ and $a(x)$ 
with continuous first derivatives, with the functions and the first derivatives
square summable on 
$\Pi$. $s$ is the surface of a sphere $S$
with center at the origin and for which the radius $r_0$ may become arbitrarily
large. Let

    $$\varphi(r_0) = \int\!\!\int_s u(x)a(x) \dee{x_i};$$
We have$^*$

    $$\varphi(r_0) = \tripl{S}\left[u(y){\partial a(y) \over \partial y_i}
             + {\partial u(y) \over \partial y_i}a(y)\right]\dee{y}.$$
The second expression shows that $\varphi(r_0)$ tends to 
a limit
$\varphi(\infty)$ when $r_0$ increases indefinitely. The first
expression for $\varphi(r_0)$ gives us

    $$|\varphi(r_0)| \leq \int\!\!\int_s |u(x)a(x)| {x_i\dee{x_i} \over r_0}$$
from which

 $$\int_0^\infty |\varphi(r_0)|\,dr_0
            \leq 
         \intpi |u(x)a(x)| \dee{x}.$$
As a result $\varphi(\infty)=0$. In other words

  $$\intpi\left[u(y){\partial a(y) \over \partial y_i}
             + {\partial u(y) \over \partial y_i}a(y)\right]\dee{y}
             = 0
  \leqno(1.11)$$
and from this we have more generally

  $$\tripl{\Pi-\varpi}\left[u(y){\partial a(y) \over \partial y_i}
             + {\partial u(y) \over \partial y_i}a(y)\right]\dee{y}
             = -\int\!\!\int_\sigma u(y)a(y) \dee{y_i}.
  \leqno(1.12)$$

   Choose as domain $\varpi$ a sphere of infinitely small radius
and center $x$ and take$^1$ in (1.12) 
$a(y) = {1 \over 4\pi}{\partial({1 \over r}) \over \partial y_i}$;
add the

\footrule
 
  $^1$ r is the distance between the points $x$ and $y$.

 $^*${\tt translator's note:} It seems that $\dee{x_1}$ means
    $dx_2 dx_3$ etc.

\rm
\vfill
\eject

 % f204

\evenheader{204}

\noindent relations qui correspondent aux valeurs 1, 2, 3 de $i$; nous 
obtenons l'identit\'e importante:

  $$u(x)={1 \over 4\pi}
         \tripl{} 
         {\partial({1 \over r}) \over \partial y_i}
         {\partial u \over \partial y_i}\dee{y}.
  \leqno(1.13)$$

  Si maintenant nous faison dans (1.11) $a(y)={y_i-x_i \over r^2} u(y)$
et si nous ajoutons les relations qui correspondent aux valeurs 1, 2, 3 
de $i$, il vient:

  $$2\tripl{\Pi}
      {y_i-x_i \over r^2}
      {\partial u \over \partial y_i}u(y)\dee{y} 
           = 
     -\tripl{\Pi} 
       {1 \over r^2}
       u^2(y)\dee{y};$$
en appliquant l'in\'egalit\'e de Schwarz au premier membre de cette identit\'e
nous obtenons une in\'egalit\'e qui nous sera utile:

  $$\tripl{\Pi} 
       {1 \over r^2}
       u^2(y)\dee{y}
            \leq
     4\tripl{\Pi} 
       {\partial u \over \partial y_i}
       {\partial u \over \partial y_i}\dee{y}.
  \leqno(1.14)$$

\medskip

   {\bf 7.} {\it Quasi-d\'eriv\'ees.}

\smallskip

   Soit une infinit\'e de fonctions $u^*(x)$ poss\'edant des d\'eriv\'ees
premi\`eres continues qui soient, comme ces fonctions elles-m\^emes, de
carr\'es sommable sur $\Pi$. Supposon que les d\'eriv\'ees
${\partial u^* \over \partial x_1}$, ${\partial u^* \over \partial x_2}$,
${\partial u^* \over \partial x_3}$ convergent faiblement en moyenne sur
$\Pi$ vers des fonctions $U_{,1}$, $U_{,2}$, $U_{,3}$. Soit $U(x)$ la
fonction mesurable d\'efinie presque partout par la relation:

  $$U(x)={1 \over 4\pi}
         \tripl{\Pi} 
         {\partial({1 \over r}) \over \partial y_i}
          U_{,i}(y)\dee{y}.$$
Nous avons:

  $$\tripl{\varpi}[u^*(x)-U(x)]^2 \dee{x}
             =
   -\tripl{\Pi} \tripl{\Pi} 
       K_{ij}(y,y')
       \left[{\partial u^* \over \partial y_i}-U_{,i}(y)\right]
       \left[{\partial u^* \over \partial {y_j}'}-U_{,j}(y')\right]
      \dee{y}
      \dee{y'}
  \leqno(1.15)$$
en posant$^1$:

\footrule

  $^1$ $r'$ repr\'esente la distance des points $x$ et $y'$.

\rm
\vfill
\eject

 % e204

\evenheader{204}

\noindent relations for values 1, 2, 3 of $i$
to obtain the important indentity

  $$u(x)={1 \over 4\pi}
         \tripl{} 
         {\partial({1 \over r}) \over \partial y_i}
         {\partial u \over \partial y_i}\dee{y}.
  \leqno(1.13)$$

  We now take $a(y)={y_i-x_i \over r^2} u(y)$ in (1.11)
and add these relations for values 1, 2, 3 
of $i$, giving

  $$2\tripl{\Pi}
      {y_i-x_i \over r^2}
      {\partial u \over \partial y_i}u(y)\dee{y} 
           = 
     -\tripl{\Pi} 
       {1 \over r^2}
       u^2(y)\dee{y}.$$
By applying the Schwarz inequality to the left side we get the useful
 inequality

  $$\tripl{\Pi} 
       {1 \over r^2}
       u^2(y)\dee{y}
            \leq
     4\tripl{\Pi} 
       {\partial u \over \partial y_i}
       {\partial u \over \partial y_i}\dee{y}.
  \leqno(1.14)$$

\medskip

   {\bf 7.} {\it Quasi-derivatives.}

\smallskip

  Let  $u^*(x)$ be an infinity of square summable functions with continuous
square summ-able first derivatives on $\Pi$.
Suppose that the derivatives 
${\partial u^* \over \partial x_1}$, ${\partial u^* \over \partial x_2}$,
${\partial u^* \over \partial x_3}$ converge weakly in mean on 
$\Pi$ to functions $U_{,1}$, $U_{,2}$, $U_{,3}$. 
Let $U(x)$ be the measurable function defined almost everywhere by

  $$U(x)={1 \over 4\pi}
         \tripl{\Pi} 
         {\partial({1 \over r}) \over \partial y_i}
          U_{,i}(y)\dee{y}.$$
We have

  $$\tripl{\varpi}[u^*(x)-U(x)]^2 \dee{x}
             =
   -\tripl{\Pi} \tripl{\Pi} 
       K_{ij}(y,y')
       \left[{\partial u^* \over \partial y_i}-U_{,i}(y)\right]
       \left[{\partial u^* \over \partial {y_j}'}-U_{,j}(y')\right]
      \dee{y}
      \dee{y'}
  \leqno(1.15)$$
where$^1$

\footrule

  $^1$ $r'$ is the distance between the points $x$ and $y'$.

\rm
\vfill
\eject

% f205

\foddheader{205}

  $$K_{ij}(y,y')=
     {1 \over 16\pi^2}
     \tripl{\varpi}
         {\partial({1 \over r}) \over \partial y_i}
         {\partial({1 \over r'}) \over \partial {y_j}'}
      \dee{x};$$
cette expression de $K$ permet d'\'etablir ais\'ement que 
l'int\'egrale

  $$\tripl{\Pi} \tripl{\Pi} 
       K_{ij}(y,y') K_{ij}(y,y') 
          \dee{y}\dee{y'}$$
est finie;
  le seconde membre de (1.15) a donc bien un sens; et il tend vers
z\'ero d'apr\`es la relation (1.7). Donc les fonctions $u^*(x)$ ont
sur tout domaine $\varpi$ une forte limite en moyenne: la fonction
$U(x)$. Et, si les int\'egrals $\intpi U^{*2}(x)\dee{x}$ sont
born\'ees dans leur ensemble, $U(x)$ est sur $\Pi$ faible limite
en moyenne des fonctions $u^*(x)$ (Cf. \S 3. Exemple II, p. 199);
on d\'eduit alors de (I.11) l'\'egalit\'e:

  $$\tripl{\Pi}
      \left[
       U(y){\partial a \over \partial y_i}
         +
       U_{,i}(y)a(y)
       \right]
       \dee{y}=0.
  \leqno(1.16)$$

\medskip

   Posons \`a ce propos la d\'efinition suivante:

\smallskip

  {\it D\'efinition des quasi-d\'eriv\'ees:} Soient deux fonctions
de carr\'es sommables sur $\Pi$, $U(y)$ et $U_{,i}(y)$; nous dirons
que $U_{,i}(y)$ est la quasi-d\'eriv\'ee de $U(x)$ par rapport \`a
$y_i$ quand la relation (1.16) sera v\'erifi\'ee; rappelons que 
dans cette relation (1.16) $a(y)$ repr\'esente une quelconque des
fonctions admettant des d\'eriv\'ees premi\`eres continues qui sont, 
comme ces fonctions elles-m\^emes, de carr\'es sommables sur $\Pi$.

  R\'esumons les r\'esultats acquis au cours de ce paragraphe:

\smallskip

  {\it Lemme 2.}  Soit une infinit\'e de fonctions $u^*(x)$ continues 
ainsi que leurs d\'eriv\'ees premi\`eres. Supposons les int\'egrales
$\intpi u^{*2}(x)\dee{x}$ born\'ee dans leur ensemble; supposons que 
chacune des d\'eriv\'ees ${\partial u^*(x) \over \partial x_i}$ ait
sur $\Pi$ une faible limite en moyenne $U_{,i}(x)$. Alors les fonctions
$u^*(x)$ convergent en moyenne vers une fonction $U(x)$, dont les 
fonctions $U_{,i}(x)$ sont des quasi-d\'eriv\'ees; cette convergence
est forte sure tout domaine $\varpi$; elle est faible$^1$ sur $\Pi$.

\footrule

  $^1$ Ou forte.

\rm
\vfill
\eject

 % e205

\oddheader{205}

  $$K_{ij}(y,y')=
     {1 \over 16\pi^2}
     \tripl{\varpi}
         {\partial({1 \over r}) \over \partial y_i}
         {\partial({1 \over r'}) \over \partial {y_j}'}
      \dee{x}.$$
This expression for $K$ allows an easy proof that the integral

  $$\tripl{\Pi} \tripl{\Pi} 
       K_{ij}(y,y') K_{ij}(y,y') 
          \dee{y}\dee{y'}$$
is finite,
  so the right side of (1.15) is defined. It tends to zero by
(1.7). Therefore the $u^*(x)$ have $U(x)$ as strong limit in mean on
all domains
 $\varpi$. And, if the integrals $\intpi U^{*2}(x)\dee{x}$ 
are bounded, $U(x)$ is the weak limit in mean of the
 $u^*(x)$ on $\Pi$ (Cf. \S 3. Example II, p. 199).
One then gets from (1.11) the equality

  $$\tripl{\Pi}
      \left[
       U(y){\partial a \over \partial y_i}
         +
       U_{,i}(y)a(y)
       \right]
       \dee{y}=0.
  \leqno(1.16)$$

\medskip

   We make the following definition:

\smallskip

  {\it Definition of quasi-derivatives:} Consider two square
summable functions $U(y)$ and $U_{,i}(y)$ on $\Pi$. We say
that $U_{,i}(y)$ is the quasi-derivative of $U(x)$ with respect
to
$y_i$ when (1.16) holds. Recall that in (1.16) $a(y)$ 
is any square summable function with continuous square summable
first derivatives on $\Pi$.

  Let us summarize the results of preceeding paragraph.

\smallskip

  {\it Lemma 2.}  Suppose we have an infinity of continuous 
functions $u^*(x)$ with continuous first derivatives. Suppose
the integrals
$\intpi u^{*2}(x)\dee{x}$ are bounded and that each of the 
derivatives
${\partial u^*(x) \over \partial x_i}$ has a weak limit in mean
 $U_{,i}(x)$ on $\Pi$. Then the 
$u^*(x)$ converge in mean to a function $U(x)$ for which the
 $U_{,i}(x)$ are the quasi-derivatives. This convergence is
strong on all domains $\varpi$. It is weak$^1$ on $\Pi$.

\footrule

  $^1$ Or strong.

\rm
\vfill
\eject

 % f206

\evenheader{206}

   De m\^eme que nous avons d\'efini les quasi-d\'eriv\'ees, nous allons
d\'efinir comme suit la {\it quasi-divergence $\Theta(x)$} d'un vecteur
$U_i(x)$ dont les composants sont de carr\'es sommables sur $\Pi$:
c'est, quand elle existe, une fonction de carr\'e sommable
v\'erifiant la relation:

  $$\intpi
       \left[
             U_i(y){\partial a \over \partial y_i}
             + \Theta(y)a(y)
       \right]
    \dee{y} = 0.
  \leqno(1.17)$$

\medskip

  {\bf 8.} {\it Approximation d'une fonction mesurable par une
            suite de fonctions r\'eguli\`eres.} Soit une quantit\'e
positive arbitraire $\epsilon$. Choisissons$^1$ une fonction $\lambda(s)$
continue, positive, d\'efinie pour $0\le s$, identique \`a z\'ero pour
$1\le s$, poss\'edant des d\'eriv\'ees de tous les ordres et telle que:

  $$ 4\pi\int_0^1 \lambda(\sigma^2)\sigma^2 \, d\sigma 
        =
     1.$$
$U(x)$ \'etant une fonction sommable sur tout domaine $\varpi$, nous
poserons

  $$\overline{U(x)}
          =
            {1 \over \epsilon^3}
            \intpi
            \lambda \left({r^2 \over \epsilon^2}\right)
            U(y)\dee{y}
  \leqno(1.18)$$
\centerline{($r=$ distance des points $x$ et $y$)}
\smallskip
\noindent Cette fonction $\overline{U(x)}$ poss\`ede des d\'eriv\'ees do tous 
les ordres:

  $${\partial^{l+m+n}\overline{U(x)} \over
     \partial x_1^l \partial x_2^m \partial x_3^n}
         =
       {1 \over \epsilon^3} 
            \intpi 
            {\partial^{l+m+n}\lambda \left({r^2 \over \epsilon^2}\right) 
                 \over
             \partial x_1^l \partial x_2^m \partial x_3^n}
           U(y)\dee{y}. 
  \leqno(1.19)$$
Supposons $U(x)$ born\'ee sur $\Pi$; nous avons manifestement:

  $$\textstyle{minimum\,\, de\,\, } U(x)
        \leq
        \overline{U(x)}
        \leq
    \textstyle{maximum\,\, de\,\, } U(x).
  \leqno(1.20)$$
Suppossons $U(x)$ de carr\'e sommable sur $\Pi$; l'in\'egalit\'e
(1.3) appliqu\'ee \`a (1.18) nous donne:

  $$\intpi \overline{U(x)^2} \dee{x} < \intpi U^2(x) \dee{x};
  \leqno(1.21)$$

\footrule

  $^1$ Pour fixer les id\'ees nous prendrons $\lambda(s)=Ae^{1 \over s-1}$,
$A$ \'etant une constante convenable, pour $0 < s < 1$.

\rm
\vfill
\eject

 % e206

\evenheader{206}

   Following our definition of quasi-derivatives we are going to define
the {\it quasi-divergence} $\Theta(x)$ of a vector
$U_i(x)$ with square summable components on $\Pi$.
When it exists, it is a square summable function with

  $$\intpi
       \left[
             U_i(y){\partial a \over \partial y_i}
             + \Theta(y)a(y)
       \right]
    \dee{y} = 0.
  \leqno(1.17)$$

\medskip

  {\bf 8.} {\it Approximation of a measurable function by a sequence
    of regular functions.} Let $\epsilon > 0$.
We choose a positive continuous function  $\lambda(s)$ defined for
$0\le s$,  identically zero for
$1\le s$ and having derivatives of all orders such that

  $$ 4\pi\int_0^1 \lambda(\sigma^2)\sigma^2 \, d\sigma 
        =
     1.$$
If $U(x)$ is summable on all domains $\varpi$, let

  $$\overline{U(x)}
          =
            {1 \over \epsilon^3}
            \intpi
            \lambda \left({r^2 \over \epsilon^2}\right)
            U(y)\dee{y}
  \leqno(1.18)$$
\centerline{($r=$ distance between $x$ and $y$)}
\smallskip
\noindent $\overline{U(x)}$ has derivatives of all orders

  $${\partial^{l+m+n}\overline{U(x)} \over
     \partial x_1^l \partial x_2^m \partial x_3^n}
         =
       {1 \over \epsilon^3} 
            \intpi 
            {\partial^{l+m+n}\lambda \left({r^2 \over \epsilon^2}\right) 
                 \over
             \partial x_1^l \partial x_2^m \partial x_3^n}
           U(y)\dee{y}. 
  \leqno(1.19)$$
If  $U(x)$ is bounded on $\Pi$ then we clearly have

  $$\min U(x)
        \leq
        \overline{U(x)}
        \leq
    \max U(x).
  \leqno(1.20)$$
If $U(x)$ is square summable on $\Pi$ the inequality
(1.3) applied to (1.18) gives$^*$ 

  $$\intpi \overline{U(x)}^2 \dee{x} \leq \intpi U^2(x) \dee{x}.
  \leqno(1.21)$$

\footrule

  $^1$ To fix ideas we take $\lambda(s)=Ae^{1 \over s-1}$,
$A$ any suitable constant, $0 < s < 1$.

  $^*$ [{\tt translator's note:} The bar was extended too far in the
  original.]

\rm
\vfill
\eject

 % f207

\foddheader{207}

\noindent appliqu\'ee \`a (1.19) elle prouve que les d\'eriv\'ees partielles
de $\overline{U(x)}$ sont de carr\'es sommables sur $\Pi$.

  Notons enfin que nous avons, si $U(x)$ et $V(x)$ sont de 
carr\'es sommables sur $\Pi$:

  $$\intpi \overline{U(x)} V(x) \dee{x}
        =
       \intpi U(x)\overline{V(x)} \dee{x}.
  \leqno(1.22)$$

   Si $V(x)$ est continue, $\overline{V(x)}$ tend uniform\'ement vers
$V(x)$ sur tout domaine $\varpi$ quand $\epsilon$ tend vers z\'ero;
on a alors d'apr\`es (1.22):

  $$\lim\intpi\overline{U(x)} V(x) \dee{x}
           =
       \intpi U(x) V(x) \dee{x};$$
on en d\'eduit que sur $\Pi$ $\overline{U(x)}$ convergent faiblement en
moyenne vers $U(x)$ quand $\epsilon$ tend vers z\'ero; 
l'in\'egalit\'e (1.21) et le crit\`ere de forte convergence \'enonc\'e
p. 200 autorisent m\^eme une conclusion plus pr\'ecise:

\medskip

   {\it Lemme 3.} Soit une fonction $U(x)$ de carr\'es sommables sur $\Pi$;
$\overline{U(x)}$ converge sur $\Pi$ fortement en moyenne vers $U(x)$ 
quand $\epsilon$ tend vers z\'ero.

\smallskip

    On \'etablit de m\^eme la proposition suivante:

\medskip

   {\it G\'en\'eralization du lemme 3.} Soit une suite de fonctions
$U_\epsilon (x)$ qui sur $\Pi$ convergent fortement (ou faiblement)
vers une limite $U(x)$ quand $\epsilon$ tend vers z\'ero; les fonctions
$\overline{U_\epsilon (x)}$ convergent fortement (ou faiblement) ver
cette m\^eme limite.

\medskip

  {\bf 9.} {\it Quelques lemmes concernant les quasi-d\'eriv\'ees.}

\smallskip

   Soit une fonction $U(x)$ de carr\'e sommable sur $\Pi$; supposons

  $$ \intpi U(x) a(x) \dee{x} = 0 $$
quelle que soit la fonction $a(x)$ de carr\'e sommable sur $\Pi$ dont les
d\'eriv\'ees de tous les ordres existent et sont de
carr\'es sommables sur $\Pi$; nous avons alors:

  $$\intpi U(x) \overline{U(x)} \dee{x} = 0 $$

\vfill
\eject

 % e207

\oddheader{207}

\noindent The same applied to (1.19) 
 proves  that the partial derivatives of $\overline{U(x)}$ 
are square summable on $\Pi$.

  Finally note that we have, if $U(x)$ and $V(x)$  
are square summable on $\Pi$

  $$\intpi \overline{U(x)} V(x) \dee{x}
        =
       \intpi U(x)\overline{V(x)} \dee{x}.
  \leqno(1.22)$$

   If $V(x)$ is continuous $\overline{V(x)}$ 
tends uniformly to $V(x)$ 
on all domains $\varpi$ when $\epsilon$ tends to zero. One therefore has 
from (1.22)

  $$\lim\intpi\overline{U(x)} V(x) \dee{x}
           =
       \intpi U(x) V(x) \dee{x}.$$
From this one deduces that $\overline{U(x)}$ converges weakly
in mean to $U(x)$ on $\Pi$ when $\epsilon$ approaches zero. 
Inequality (1.21) and the criteria for strong convergence on
p. 200 similarly give a more precise conclusion:

\medskip

   {\it Lemma 3.} Let $U(x)$ be square summable on $\Pi$.
$\overline{U(x)}$ converges 
strongly in mean to $U(x)$ 
on $\Pi$ when $\epsilon$ tends to zero.

\smallskip

    Similarly one establishes the following proposition.

\medskip

   {\it Generalization of lemma 3.} Suppose a sequence of functions
$U_\epsilon (x)$ converge strongly (or weakly) in mean on $\Pi$ to a limit
$U(x)$ as $\epsilon$ tends to zero. Then the functions 
$\overline{U_\epsilon (x)}$ converge strongly (or weakly) 
to the same limit.

\medskip

  {\bf 9.} {\it Some lemmas on quasi-derivatives.}

\smallskip

   Let $U(x)$ be square summable on $\Pi$. Suppose that for all
square summable functions
$a(x)$ having square summable  derivatives of all orders

  $$ \intpi U(x) a(x) \dee{x} = 0 $$
then

  $$\intpi U(x) \overline{U(x)} \dee{x} = 0 $$

\vfill
\eject

% f208

\evenheader{208}

\noindent d'o\`u, en faisant tendre $\epsilon$ vers z\'ero:

  $$ \intpi U^2(x) \dee{x} = 0.$$
La fonction $U(x)$ est donc nulle presque partout.

  Ce fait permet d'\'etablir les propositions suivantes: La 
quasi-d\'eriv\'ee d'une fonction par rapport \`a la variable $x_i$ est unique 
quand elle existe. (Nous consid\'erons comme identique deus fonctions \'egales
presque partout.)

  La quasi-divergence d'un vecteur est unique quand elle existe.

\medskip

  {\it Lemme 4.} Soit une fonctions $U(x)$ admettant une
quasi-d\'eriv\'ee $U_{,i}(x)$; je dis que 
${\partial\overline{U(x)}
    \over
   \partial x_i} = \overline{U_{,i}(x)}.$

\medskip
\ninerm
    Il suffit de prouver que:

  $$\intpi
        {\partial \overline{U(x)} \over \partial x_i}a(x)
    \dee{x}
       =
     \intpi
        U_{,i}(x)a(x)
    \dee{x}.$$
Or on d\'eduit ais\'ement de (1.18) que:

  $${\partial \overline{a(x)} \over \partial x_i}
         =
      \overline{\partial a(x) \over \partial x_i};$$
cette formule et les formules (1.11), (1.16), (1.22) justifient les
transformations:

  $$ \intpi
        {\partial \overline{U(x)} \over \partial x_i}a(x)
       \dee{x}
       =
       -\intpi
          \overline{U(x)} {\partial a(x) \over \partial x_i}
        \dee{x}
       = 
       -\intpi
          U(x) \overline{\left({\partial a(x) \over \partial x_i}\right)}
        \dee{x}
       = $$
  $$
       -\intpi
          U(x) {\partial \overline{a(x)} \over \partial x_i}
        \dee{x}
      =
        \intpi
           U_{,i}(x)\overline{a(x)}
        \dee{x}
       = 
        \intpi
           \overline{U_{,i}(x)} a(x)
        \dee{x}. \eqno{C.Q.F.D.}$$

\rm

\medskip

   {\it Lemme 5.} Soient deux fonctions de carr\'es sommables sur $\Pi$,
$U(x)$ et $V(x)$, qui poss\`edent les quasi-d\'eriv\'ees $U_{,i}(x)$ et
$V_{,i}(x)$; je dis que:

  $$\intpi [U(x)V_{,i}(x)+U_{,i}(x)V(x)]\dee{x}
      = 0.
  \leqno(1.23)$$

\medskip

\ninerm
    
  Cette formule s'obtient en appliquant le lemme 3 \`a la formule:

  $$\intpi [U(x)\overline{V_{,i}(x)}+U_{,i}(x)\overline{V(x)}]\dee{x} 
      = 0.$$
qui elle-m\^eme r\'esulte de la relation (1.16) et du lemme 4.

\rm
\vfill
\eject

 % e208

\evenheader{208}

\noindent from which one gets as $\epsilon$ tends to zero

  $$ \intpi U^2(x) \dee{x} = 0.$$
The function $U(x)$ is therefore zero almost everywhere.

  That fact allows us to establish the following propositions. 1) When the
quasi-derivative of a function with respect to $x_i$ exists, it is unique.
(We consider two functions identical if they are equal almost everywhere.)
2) The quasi-divergence of a vector is unique if it exists.

\medskip

  {\it Lemma 4.} Let $U(x)$ have a quasi-derivative $U_{,i}(x)$. Then I claim that 
${\partial\overline{U(x)}
    \over
   \partial x_i} = \overline{U_{,i}(x)}.$

\medskip
\ninerm
    It suffices to prove that

  $$\intpi
        {\partial \overline{U(x)} \over \partial x_i}a(x)
    \dee{x}
       =
     \intpi
        U_{,i}(x)a(x)
    \dee{x}.$$
Because one easily deduces from (1.18) that

  $${\partial \overline{a(x)} \over \partial x_i}
         =
      \overline{\partial a(x) \over \partial x_i}$$
and this formula with (1.11), (1.16), and (1.22) justify the
transformations

  $$ \intpi
        {\partial \overline{U(x)} \over \partial x_i}a(x)
       \dee{x}
       =
       -\intpi
          \overline{U(x)} {\partial a(x) \over \partial x_i}
        \dee{x}
       = 
       -\intpi
          U(x) \overline{\left({\partial a(x) \over \partial x_i}\right)}
        \dee{x}
       = $$
  $$
       -\intpi
          U(x) {\partial \overline{a(x)} \over \partial x_i}
        \dee{x}
      =
        \intpi
           U_{,i}(x)\overline{a(x)}
        \dee{x}
       = 
        \intpi
           \overline{U_{,i}(x)} a(x)
        \dee{x}. \eqno{Q.E.D.}$$

\rm

\medskip

   {\it Lemma 5.} Suppose that two square summable functions 
$U(x)$ and $V(x)$ have quasi-derivatives $U_i(x)$ and
$V_i(x)$ on $\Pi$. I claim that

  $$\intpi [U(x)V_{,i}(x)+U_{,i}(x)V(x)]\dee{x}
      = 0.
  \leqno(1.23)$$

\medskip

\ninerm
    
  This is obtained by applying lemma 3 to the formula

  $$\intpi [U(x)\overline{V_{,i}(x)}+U_{,i}(x)\overline{V(x)}]\dee{x} 
      = 0.$$
which follows from (1.16) and lemma 4.

\rm
\vfill
\eject

 % f209

\foddheader{209}

  {\it Lemme 6.} Soit un vecteur $U_i(x)$ admettant une 
quasi-divergence $\Theta(x)$; on a: divergence 
$\overline{U_i(x)} = \overline{\Theta(x)}$.

\smallskip
\ninerm
   (La d\'emonstration de ce lemme est tr\`es analogue \`a celle
du lemme 4).
\rm
\medskip

  {\it Lemme 7.} Soit un vecteur $U_i(x)$ de quasi-divergence nulle.
Supposons $$\intpi U_i(x) a_i(x)\dee{x}=0$$ quel que soit le 
vecteur $a_i(x)$, de divergence nulle, dont les composantes ainsi que leurs
d\'eriv\'ees de tous les ordres sont de carr\'es sommables sur $\Pi$. Je 
dis que $U_i(x)=0$.

\medskip
\ninerm
  En effet le lemme 4 nous autorise \`a choisir $a_i(x)=\overline{U_i(x)}$;
or quand $\epsilon$ tend vers z\'ero la relation
$\intpi U_i(x)\overline{U_i(x)}\dee{x}=0$ se r\'eduit \`a la suivante:

  $$\intpi U_i(x)U_i(x)\dee{x}=0.$$

\rm
\smallskip
  {\it Corollaire.} Une infinit\'e de vecteurs $U_i^*(x)$, de quasi-divergence
nulle, poss\`ede sur $\Pi$ une faible limite en moyenne unique si les
deus conditions suivantes sont v\'erifi\'ees:

\smallskip

   a) Les nombres $\intpi U_i^*(x)U_i^*(x)\dee{x}$ sont born\'es dans
leurs ensemble;

\smallskip

   b) Pour chaque vecteur $a_i(x)$ de divergence nulle, dont les composantes,
ainsi que leurs d\'eriv\'ees de tous les oedres, sont de carr\'es 
sommables sur $\Pi$, les quantit\'es 
$\intpi U_i^*(x)a_i(x)\dee{x}$ ont une seule valeur limite.

\medskip
\ninerm
    Sinon le Th\'eorem fondamental de M. F. Riesz (p. 202) permettrait
d'extraire de la suite $U_i^*(x)$ deus suites partielles poss\'edant
deux faibles limites distinctes, dont l'existence contredireit le lemme 7.

\rm

\bigskip
\centerline{\bf II. Mouvements infiniment lents.}

\medskip

   {\bf 10.} On d\'esigne par $<<$ \'equations de Navier lin\'earis\'ees $>>$
les \'equations suivantes:

  $$\leqno(2.1)$$
  $$\nu\Delta u_i(x,t)
          -{\partial u_i(x,t)\over \partial t}
          -{1\over \rho}{\partial p(x,t)\over \partial x_i}
          =
          -X_i(x,t)
          \left[\Delta=
            {\partial^2\over \partial x_k \partial x_k}
          \right]$$
  $$ {\partial u_j(x,t)\over \partial x_j}=0. $$

\vfill
\eject

 % e209

\oddheader{209}

  {\it Lemma 6.} If a vector $U_i(x)$ has quasi-divergence
$\Theta(x)$  then the divergence of
$\overline{U_i(x)} = \overline{\Theta(x)}$.

\smallskip
\ninerm
   (The proof is very much analogous to that of lemma 4.)
\rm
\medskip

  {\it Lemma 7.} Suppose a vector  $U_i(x)$ has quasi-divergence 0,
and that $$\intpi U_i(x) a_i(x)\dee{x}=0$$  for all square summable vectors
$a_i(x)$ which have 0 divergence and square summable derivatives of all
orders on
$\Pi$. Then I claim that   $U_i(x)=0$.

\medskip
\ninerm
  In fact lemma 4 allows us to choose $a_i(x)=\overline{U_i(x)}$ because when
$\epsilon$ tends to 0 the relation
$\intpi U_i(x)\overline{U_i(x)}\dee{x}=0$ 
reduces to

  $$\intpi U_i(x)U_i(x)\dee{x}=0.$$

\rm
\smallskip
  {\it Corollary.} An infinity of vectors
 $U_i^*(x)$, of quasi-divergence 0 has, on $\Pi$, a unique weak limit in mean
if the two following conditions hold:

\smallskip

   a) the numbers $\intpi U_i^*(x)U_i^*(x)\dee{x}$ 
are bounded

\smallskip

   b) for all square summable vectors
$a_i(x)$ which have 0 divergence and square summable derivatives of all
orders on
$\Pi$,
the quantities $\intpi U_i^*(x)a_i(x)\dee{x}$ have a single limiting value.

\medskip
\ninerm
   If not, then the fundamental theorem of F. Riesz (p. 202) allows
extraction from the sequence 
 $U_i^*(x)$ two subsequences  having distinct limits. This 
contradicts lemma 7. 

\rm

\bigskip
\centerline{\bf II. Infinitely slow motion.}

\medskip

   {\bf 10.} The ``linearised Navier equations'' are
the following

  $$\nu\Delta u_i(x,t)
          -{\partial u_i(x,t)\over \partial t}
          -{1\over \rho}{\partial p(x,t)\over \partial x_i}
          =
          -X_i(x,t)
          \left[\Delta=
            {\partial^2\over \partial x_k \partial x_k}
          \right]
  \leqno(2.1)$$
  $$ {\partial u_j(x,t)\over \partial x_j}=0. $$

\vfill
\eject
 % f210

\evenheader{210}

\noindent $\nu$ et $\rho$ sont les constantes donn\'ees, $X_i(x,t)$
est un vecteur donn\'e qui repr\'esente les forces ext\'erieures;
$p(x,t)$ repr\'esente la pression, $u_i(x,t)$ la vitesse des
mol\'ecules du liquide.

\medskip

  {\it Le probl\'eme} que pose la th\'eorie des liquides visqueux est
le suivant:

  Construire pour $t>0$ la solution de (2.1) qui correspond \`a des
valeurs initiales donn\'ees, $u_i(x,0)$.

  Nous allons rappeler la solution de ce probl\`eme et quelques-unes
des propri\'et\'es qu'elle poss\`ede. Nous posserons:

  $$W(t) = \intpi u_i(x,t)u_i(x,t)\dee{x} $$

  $$J_m^2(t) 
       =
       \intpi
       {\partial^m u_i(x,t)\over 
         \partial x_k \partial x_l \ldots}
        {\partial^m u_i(x,t)\over    
         \partial x_k \partial x_l \ldots}
       \dee{x}.$$ 

  $V(t)=$ Maximum de $\sqrt{u_i(x,t)u_i(x,t)}$ \`a l'instant $t$.

  $D_m(t)=$ Maximum des fonctions 
              $\left|
                 {\partial^m u_i(x,t)\over    
                   \partial x_1^h \partial x_2^k \partial x_3^l}
               \right|$
        \`a l'instant $t$ ($h+k+l=m$).

\smallskip

   Nous ferons les hypoth\`eses suivantes relativement aux donn\'ees:
Les fonctions $u_i(x,t)$ et leurs d\'eriv\'ees premi\`eres sont continues;
${\partial u_j(x,0)\over \partial x_j}=0$; les quantit\'es $W(0)$ 
et $V(0)$ sont finies; $|X_i(x,t)-X_i(y,t)|<r^{1\over 2}C(x,y,t)$,
$C(x,y,t)$ \'etant une fonction continue;
$\intpi X_i(x,t)X_i(x,t)\dee{x}$ est une fonction continue de $t$, ou
est inf\'erieure \`a une fonction continue de $t$.

   Les lettres $A$ et $A_m$ nous serviront d\'esormais \`a d\'esigner 
les constantes et les fonctions de l'indice $m$ dont nous ne pr\'eciserons
pas les valeurs num\'eriques.

\medskip

  {\bf 11. } {\it Premier cas particulier:} $X_i(x,t)=0$.

\smallskip

   La th\'eorie de la Chaleur fournit dans ce cas la solution suivante
du syst\`eme (2.1):

  $$u'_i(x,t) = 
         {1\over (2\sqrt{\pi})^3}
         \intpi
          {
           e^{-{r^2\over 4\nu t}} 
             \over
           (\nu t)^{3\over 2}
          }
          u_i(y,0)\dee{y}; \quad p'(x,t)=0.
  \leqno(2.2)$$

  Les int\'egrales $u'_i(x,t)$ sont uniform\'ement continues en $t$
(cf. \S 5, p. 202) pour $0<t$, et l'on a:

\vfill
\eject

 % e210

\evenheader{210}

\noindent $\nu$ and $\rho$ are given constants, $X_i(x,t)$
is a vector which represents external forces,
$p(x,t)$ is the pressure, and $u_i(x,t)$ the speed of the molecules
of the liquid.

\medskip

  {\it The problem} posed by the theory of viscous liquids is the following:
Construct for $t>0$ the solution of (2.1) which has given initial values
$u_i(x,0)$.

  We recall the solution of this problem and some of its properties. 
Write

  $$W(t) = \intpi u_i(x,t)u_i(x,t)\dee{x} $$

  $$J_m^2(t) 
       =
       \intpi
       {\partial^m u_i(x,t)\over 
         \partial x_k \partial x_l \ldots}
        {\partial^m u_i(x,t)\over    
         \partial x_k \partial x_l \ldots}
       \dee{x}.$$ 

  $V(t)=$ Maximum of $\sqrt{u_i(x,t)u_i(x,t)}$ at time $t$.

  $D_m(t)=$ Maximum of the function 
              $\left|
                 {\partial^m u_i(x,t)\over    
                   \partial x_1^h \partial x_2^k \partial x_3^l}
               \right|$
        at time  $t$ ($h+k+l=m$).

\smallskip

   We make the following assumptions:
The functions $u_i(x,t)$ and their first derivatives are continuous,
${\partial u_j(x,0)\over \partial x_j}=0$, the quantities $W(0)$ 
and $V(0)$ are finite, $|X_i(x,t)-X_i(y,t)|<r^{1\over 2}C(x,y,t)$,
where $C(x,y,t)$ is a continuous function, and
$\intpi X_i(x,t)X_i(x,t)\dee{x}$ is a continuous function of $t$, or
is less than a continuous function of $t$.

   From now on the letters $A$ and $A_m$ denote constants and functions
with index $m$ for which we do not specify the numerical value.

\medskip

  {\bf 11. } {\it First case:} $X_i(x,t)=0$.

\smallskip

   The theory of heat gives the following solution$^*$ to system (2.1):

  $$u'_i(x,t) = 
         {1\over (2\sqrt{\pi})^3}
         \intpi
          {
           e^{-{r^2\over 4\nu t}} 
             \over
           (\nu t)^{3\over 2}
          }
          u_i(y,0)\dee{y}; \quad p'(x,t)=0.
  \leqno(2.2)$$

  The integrals  $u'_i(x,t)$ are uniformly continuous in $t$
(cf. \S 5, p. 202) for $0<t$, and from this one has

\footrule

[{\tt translator's note:} where $r=|y-x|$.

\vfill
\eject

% f211

\foddheader{211}

  $$V(t) < V(0).
  \leqno(2.3)$$

  Quand $J_1(0)$ est fini, l'appliction de l'in\'egalit\'e (1.14) et de
l'in\'egalit\'e de Schwarz (1.1) \`a (2.2) permet d'obtenir une 
seconde borne de $V(t)$:

  $$V^2(t) 
        <
       4J_1^2(0)
       {1\over (4\pi)^3}
        \intpi
          {
           e^{-{r^2\over 2\nu t}} 
             \over
           (\nu t)^3
          }
          r^2 \dee{y},$$
c'est-\`a-dire:

  $$ V(t) < {AJ_1(0) \over
                    (\nu t)^{1\over 4}}.
  \leqno(2.4)$$

  L'in\'egalit\'e (1.3) appliqu\'ee \`a (2.2) prouve que l'on a:

  $$ W(t) < W(0)
  \leqno(2.5)$$
les int\'egrales $u_i'(x,t)$ sont fortement continues en $t$ 
(cf. \S 5, p. 202)  m\^eme pour $t=0$. Appliqu\'ee \`a la relation:

  $${\partial u_i'(x,t)
       \over
     \partial x_k}
        =
       {1\over (2\sqrt{\pi})^3}
         \intpi
          {
          \partial
          \over { \partial x_k }
          }
          \left[ { 
           e^{-{r^2\over 4\nu t}} 
             \over
           (\nu t)^{3\over 2}
          } \right]
          u_i(y,0)\dee{y}$$
cette in\'egalit\'e (1.3) prouve que:

  $$ J_1(t) < J_1(0);
  \leqno(2.6)$$
les d\'eriv\'ees pr\'emi\`eres ${\partial u_i' \over \partial x_k}$
sont fortement continues en $t$, m\^eme pour $t=0$ si $J_1(0)$ est fini.

  Pour les raisons analogues les d\'eriv\'ees de tous les ordres des
int\'egrales $u_i'(x,t)$ sont uniform\'ement et fortement continues
en $t$ pour $t>0$; et plus pr\'ecis\'ement:

  $$ D_m(t) < {A_m\sqrt{W(0)}\over (\nu t)^{2m+3\over 4}},
  \leqno(2.7)$$
  $$ J_m(t) < {A_m\sqrt{W(0)}\over (\nu t)^{m\over 2}}.
  \leqno(2.8)$$

\medskip

  {\bf 12. } {\it Second cas particulier; } $u_i'(x,0)=0$.

\smallskip

  La solution fondamental de M. Oseen$^1$, $T_{ij}(x,t)$, fournit la
solution suivante du syst\`eme (2.1):

\footrule

  $^1$ Voir: Oseen: Hydrodynamik \S 5; Acta mathematice T. 34.

\rm

\vfill
\eject

 % e211

\oddheader{211}

  $$V(t) < V(0).
  \leqno(2.3)$$

  If $J_1(0)$ is finite, inequality (1.14) and the Schwarz inequality (1.1)
applied to  (2.2) give  a second bound on
$V(t)$:

  $$V^2(t) 
        <
       4J_1^2(0)
       {1\over (4\pi)^3}
        \intpi
          {
           e^{-{{\bf r}^2\over 2\nu t}} 
             \over
           (\nu t)^3
          }
          r^2 \dee{y},$$
which is to say

  $$ V(t) < {AJ_1(0) \over
                    (\nu t)^{1\over 4}}.
  \leqno(2.4)$$

  Inequality (1.3) applied to (2.2) proves:

  $$ W(t) < W(0);
  \leqno(2.5)$$
the integrals $u_i'(x,t)$ are strongly continuous in $t$ 
(cf. \S 5, p. 202)  including  $t=0$. Inequality (1.3) applied to

  $${\partial u_i'(x,t)
       \over
     \partial x_k}
        =
       {1\over (2\sqrt{\pi})^3}
         \intpi
          {
          \partial
          \over { \partial x_k }
          }
          \left[ { 
           e^{-{r^2\over 4\nu t}} 
             \over
           (\nu t)^{3\over 2}
          } \right]
          u_i(y,0)\dee{y}$$
proves that

  $$ J_1(t) < J_1(0);
  \leqno(2.6)$$
the first derivatives ${\partial u_i' \over \partial x_k}$
are strongly continuous in $t$, including  $t=0$ if $J_1(0)$ is finite.

  For analogous reasons the derivatives of all orders of
$u_i'(x,t)$ are uniformly and strongly continuous in
$t$ for $t>0$ and more precisely

  $$ D_m(t) < {A_m\sqrt{W(0)}\over (\nu t)^{2m+3\over 4}},
  \leqno(2.7)$$
  $$ J_m(t) < {A_m\sqrt{W(0)}\over (\nu t)^{m\over 2}}.
  \leqno(2.8)$$

\medskip

  {\bf 12. } {\it Second particular case } $u_i'(x,0)=0$.

\smallskip

  Oseen's fundamental solution$^1$, $T_{ij}(x,t)$, furnishes the following
solution to system (2.1):

\footrule

  $^1$ See: Oseen: Hydrodynamik \S 5; Acta mathematica vol. 34. 

[{\tt translator's note:} which gives  on p. 41

\noindent 
 $t_{jk}=\delta_{jk}{1\over 2\nu}{E(r,t^{(0)}-t)\over t^{(0)}-t}
            +
         {\partial^2\Phi\over \partial x_j \partial x_k}$, 
$\Phi = {1\over r}\int_0^r E(\alpha, t^{(0)}-t)\,d\alpha$,
$E(r,t^{(0)}-t)={
           e^{-{{\bf r}^2\over 4\nu (t^{(0)}-t)}} 
             \over
           \sqrt{t^{(0)}-t}
          }$.]
\rm
\vfill
\eject

 % f212

\evenheader{212}

  $$\leqno(2.9)$$
  $$ u_i ''(x,t) = \int_0^t \,dt' \,\intpi T_{ij}(x-y,t-t')X_j(y,t')\dee{y}$$
  $$ p''(x,t) = -{\rho\over 4\pi}{\partial\over\partial x_j}
                  \intpi {1\over r}X_j(y,t)\dee{y}$$
Nous avons:

  $$\leqno(2.10)$$
  $$|T_{ij}(x-y,t-t')| < {A\over [r^2+\nu(t-t')]^{3\over 2}};$$
  $$\left| {\partial^m T_{ij}(x-y,t-t')\over
         \partial x_1^h \partial x_2^k \partial x_3^l}\right|
          <
         {A_m \over [r^2+\nu(t-t')]^{m+3\over 2}}; (t'<t).$$

  Nous remarquerons en premier lieu que les int\'egrals (1.2) et (1.3)
appliqu\'ees en m\^eme temps que (2.10) \`a la formule:

  $${\partial u_i ''(x,t)\over
        \partial x_k} = 
       \int_0^t \,dt' \,\intpi
           {\partial T_{ij}(x-y,t-t')\over
            \partial x_k}X_j(y,t')\dee{y}
  \leqno(2.11)$$
prouvent que les d\'eriv\'ees premi\`eres 
${\partial u_i ''\over \partial x_k}$ sont fortement continues en $t$ pour
$t\geq 0$, et que:

  $$J_1(t) < A \int_0^t \,{dt'\over\sqrt{\nu (t-t')}}
              \sqrt{\intpi X_i(x,t')X_i(x,t')\dee{x}}.
  \leqno(2.12)$$

  Ceci fait, adjoignons aux hypoth\'eses d\'ej\`a \'enonc\'ees la suivante:
le maximum de $\sqrt{X_i(x,t')X_i(x,t')}$ \`a l'instant $t$ est une fonction
continue de $t$, ou est inf\'erieur \`a une fonction continue de $t$;
il n'y a aucune difficult\'e \`a d\'eduire de (2.9) que
$u_i ''(x,t)$ et ${\partial u_i ''\over \partial x_k}$ sont alors 
uniform\'ement continue en $t$ pour $t\geq 0$, et \`a pr\'eciser par exemple
que

  $$D_1(t) < A \int_0^t \,{dt'\over\sqrt{\nu (t-t')}}
              \max\sqrt{X_i(x,t')X_i(x,t')};
  \leqno(2.13)$$
cette in\'egalit\'e (2.13) peut \^etre compl\'et\'ee comme suit: nous avons

\vfill
\eject

 % e212

\evenheader{212}

  $$ u_i ''(x,t) = \int_0^t \,dt' \,\intpi T_{ij}(x-y,t-t')X_j(y,t')\dee{y}
  \leqno(2.9)$$
  $$ p''(x,t) = -{\rho\over 4\pi}{\partial\over\partial x_j}
                  \intpi {1\over r}X_j(y,t)\dee{y}$$
We have

  $$|T_{ij}(x-y,t-t')| < {A\over [r^2+\nu(t-t')]^{3\over 2}}
  \leqno(2.10)$$
  $$\left| {\partial^m T_{ij}(x-y,t-t')\over
         \partial x_1^h \partial x_2^k \partial x_3^l}\right|
          <
         {A_m \over [r^2+\nu(t-t')]^{m+3\over 2}}; \ \  (t'<t).$$

  We remark in the first place that integrals (1.2) and (1.3)
applied with (2.10) to the formula

  $${\partial u_i ''(x,t)\over
        \partial x_k} = 
       \int_0^t \,dt' \,\intpi
           {\partial T_{ij}(x-y,t-t')\over
            \partial x_k}X_j(y,t')\dee{y}
  \leqno(2.11)$$
prove that the first derivatives
${\partial u_i ''\over \partial x_k}$ are strongly continuous
in  $t$ for 
$t\geq 0$, and that

  $$J_1(t) < A \int_0^t \,{dt'\over\sqrt{\nu (t-t')}}
              \sqrt{\intpi X_i(x,t')X_i(x,t')\dee{x}}.
  \leqno(2.12)$$

  This done, we add to previously stated hypotheses the assumption
that the maximum of  
$\sqrt{X_i(x,t')X_i(x,t')}$ at time $t$ is a continuous function of
$t$, or is less than a continuous function of  $t$.
Then there is no difficulty in deducing from (2.9) that
$u_i ''(x,t)$ and ${\partial u_i ''\over \partial x_k}$ are uniformly
continuous in 
$t$ for $t\geq 0$, and more precisely for example

  $$D_1(t) < A \int_0^t \,{dt'\over\sqrt{\nu (t-t')}}
              \max\sqrt{X_i(x,t')X_i(x,t')}.
  \leqno(2.13)$$
Inequality (2.13) may be complemented as follows. We have

\vfill
\eject

 % f213

\foddheader{213}

  $$ {\partial u_i''(x,t)
         \over
      \partial x_k} -
     {\partial u_i''(y,t)
         \over
      \partial y_k} = \int_0^t \,dt'\tripl{\varpi}
            {\partial T_{ij}(x-z,t-t')
               \over
             \partial x_k}X_j(z,t')\dee{z}$$
  $$ -  \int_0^t \,dt'\tripl{\varpi}
            {\partial T_{ij}(y-z,t-t')
               \over
             \partial y_k}X_j(z,t')\dee{z}$$

  $$ + \int_0^t \,dt'\tripl{\Pi-\varpi}
            \left[{\partial T_{ij}(x-z,t-t')
               \over
             \partial x_k}
            -{\partial T_{ij}(y-z,t-t')
               \over
             \partial y_k}\right]X_j(z,t')\dee{z},$$
$\varpi$ \'etant le domaine des points $z$ situ\'es \`a une distance 
de $x$ ou de $y$ inf\'erieure \`a $2r$; appliquons la formule des 
accroissements finis au crochet:

  $$\left[{\partial T_{ij}(x-z,t-t')
               \over
             \partial x_k} -
            {\partial T_{ij}(y-z,t-t')
               \over
             \partial y_k}\right]$$
et majorons les trois int\'egrales pr\'ec\'edentes en rempla\c cant les
diverses fonctions qui y figurent par des majorantes de leurs valeurs
absolues; nous v\'erifions  ais\'ement que:

  $$\left|{\partial u_i''(x,t)
         \over
      \partial x_k} -
     {\partial u_i''(y,t)
         \over
      \partial y_k}\right| <
  \leqno(2.14)$$
  $$ < Ar^{1\over 2}\int_0^t {dt'\over [\nu (t-t')]^{3\over 4}}
              \max\sqrt{X_i(x,t')X_i(x,t')}.$$
--- Nous dirons {\it qu'une fonction $U(x,t)$ satisfait une condition H}
quand elle v\'erifie une in\'egalit\'e analogue \`a la pr\'ec\'edente:

  $$|U(x,t)-U(y,t)| < r^{1\over 2}C(t),
  \leqno(2.15)$$
o\'u $C(t)$ est inf\'erieur \`a une fonction continue de $t$. Nous nommerons
coefficient de la condition H celle des fonctions $C(t)$ dont les valeurs
sont les plus faibles possibles. ---

  Supposons maintenant que les fonctions $X_i(x,t)$ satisfassent une telle
condition H, de coefficient $C(t)$; les d\'eriv\'ees secondes 
${\partial^2 u_i''(x,t)\over \partial x_k \partial x_l}$, qui sont
donn\'ees par les formules:

\vfill
\eject

 % e213

\oddheader{213}

  $$ {\partial u_i''(x,t)
         \over
      \partial x_k} -
     {\partial u_i''(y,t)
         \over
      \partial y_k} = \int_0^t \,dt'\tripl{\varpi}
            {\partial T_{ij}(x-z,t-t')
               \over
             \partial x_k}X_j(z,t')\dee{z}$$
  $$ -  \int_0^t \,dt'\tripl{\varpi}
            {\partial T_{ij}(y-z,t-t')
               \over
             \partial y_k}X_j(z,t')\dee{z}$$

  $$ + \int_0^t \,dt'\tripl{\Pi-\varpi}
           \left[ {\partial T_{ij}(x-z,t-t')
               \over
             \partial x_k}
          -  {\partial T_{ij}(y-z,t-t')
               \over
             \partial y_k}\right]X_j(z,t')\dee{z},$$
$\varpi$ being the domain of points at distance less than $2r$ to $x$ or $y$.
We apply the formula of finite differences to the bracket

  $$\left[{\partial T_{ij}(x-z,t-t')
               \over
             \partial x_k} -
            {\partial T_{ij}(y-z,t-t')
               \over
             \partial y_k}\right]$$
and majorize the preceeding three integrals by replacing the various
functions there by the majorants of their absolute values. We easily verify

  $$\left|{\partial u_i''(x,t)
         \over
      \partial x_k} -
     {\partial u_i''(y,t)
         \over
      \partial y_k}\right| <
  \leqno(2.14)$$
  $$ < Ar^{1\over 2}\int_0^t {dt'\over [\nu (t-t')]^{3\over 4}}
              \max\sqrt{X_i(x,t')X_i(x,t')}.$$
--- We say that {\it a function $U(x,t)$ satisfies condition H}
if an inequality analogous to the preceeding holds:

  $$|U(x,t)-U(y,t)| < r^{1\over 2}C(t),
  \leqno(2.15)$$
where $C(t)$ is smaller than a continuous function of
$t$. We call the weakest possible $C(t)$, the condition H coefficient. ---

  Now suppose that the functions $X_i(x,t)$ satisfy 
condition H with  coefficient $C(t)$. Then the second derivatives
${\partial^2 u_i''(x,t)\over \partial x_k \partial x_l}$,
given by the formulas

\vfill
\eject

% f214

\evenheader{214}

  $$ {\partial^2 u_i''(x,t)
         \over
      \partial x_k \partial x_l} = 
      \int_0^t \,dt'\intpi
            {\partial^2 T_{ij}(x-y,t-t')
               \over
             \partial x_k \partial x_l}
            [X_j(y,t')-X_j(x,t')]\dee{y}$$
sont alors des fonctions uniform\'ement continues en $t$, et l'on \`a:

  $$ D_2(t) < A \int_0^t {C(t')\, dt' \over [\nu (t-t')]^{3\over 4}}.
  \leqno(2.16)$$

  Plus g\'en\'eralement:

  Supposons que les d\'eriv\'ees d'ordre $m$ des fonctions $X_i(x,t)$ 
par
rapport \`a $x_1$, $x_2$, $x_3$ existent, soient continues et soient
inf\'erieures en valeur absolue \`a une fonctions continue $\varphi_m(t)$. 
Alors les d\'eriv\'ees d'ordre $m+1$, par rapport \`a $x_1$, $x_2$, $x_3$,
des fonctions $u_i''(x,t)$ existent, sont uniform\'ement continues en
$t$; on a:

  $$ D_{m+1}(t) < A \int_0^t {\varphi_m(t')\, dt'
                \over \sqrt{\nu (t-t')}}
  \leqno(2.17)$$
enfin ces d\'eriv\'ees d'ordre $m+1$ satisfont condition H de coefficient:

  $$ C_{m+1}(t) < A \int_0^t {\varphi_m(t')\, dt' \over [\nu (t-t')]^{3\over 4}}.
  \leqno(2.18)$$
Si de plus:

  $$\intpi\left[{\partial^m X_i(x,t)\over 
            \partial x_1^h \partial x_2^k \partial x_3^l}\right]^2\dee{x}
           < \psi_m^2(t),$$
$\psi_m(t)$ \'etant une fonction continue (positive), alors les
d\'eriv\'ees d'ordre $m+1$ par rapport \`a $x_1$, $x_2$, $x_3$ des fonctions
$u_i(x,t)$ sont fortement continues en $t$ et v\'erifient l'in\'egalit\'e:

  $$ J_{m+1}(t) < A_m \int_0^t {\psi_m(t')\, dt' 
                \over \sqrt{\nu (t-t')}}
  \leqno(2.19)$$
Supposons maintenant que les d\'eriv\'ees d'ordre $m$ des fonctions $X_i(x,t)$
par rapport \`a $x_1$, $x_2$, $x_3$ existent, soient inf\'erieures en valeur
absolue \`a une fonction continue de $t$ et v\'erifient une condition H de
coefficient $\theta_m(t)$. Alors les d\'eriv\'ees d'ordre $m+2$ des fonctions
$u_i(x,t)$ par rapport \`a $x_1$, $x_2$, $x_3$ existent, sont uniform\'ement
continues en $t$ et v\'erifient l'in\'egalit\'e:

\vfill
\eject

 % e214

\evenheader{214}

  $$ {\partial^2 u_i''(x,t)
         \over
      \partial x_k \partial x_l} = 
      \int_0^t \,dt'\intpi
            {\partial^2 T_{ij}(x-y,t-t')
               \over
             \partial x_k \partial x_l}
            [X_j(y,t')-X_j(x,t')]\dee{y},$$
are then uniformly continuous in $t$ and

  $$ D_2(t) < A \int_0^t {C(t')\, dt' \over [\nu (t-t')]^{3\over 4}}.
  \leqno(2.16)$$

  More generally:

  Suppose the $m$-th order derivatives of the $X_i(x,t)$ 
with respect to  $x_1$, $x_2$, $x_3$ 
exist, are continuous, and are smaller than some
continuous functions
$\varphi_m(t)$. Then the derivatives of order $m+1$ of the 
$u_i''(x,t)$ with respect to $x_1$, $x_2$, $x_3$ 
exist and are uniformly continuous in 
$t$. We have

  $$ D_{m+1}(t) < A \int_0^t {\varphi_m(t')\, dt'
                \over \sqrt{\nu (t-t')}}
  \leqno(2.17)$$
and finally the derivatives of order $m+1$ 
satisfy condition H with coefficient

  $$ C_{m+1}(t) <
        A \int_0^t {\varphi_m(t')\, dt' \over [\nu (t-t')]^{3\over 4}}.
  \leqno(2.18)$$
If further

  $$\intpi\left[{\partial^m X_i(x,t)\over 
            \partial x_1^h \partial x_2^k \partial x_3^l}\right]^2\dee{x}
           < \psi_m^2(t),$$
where $\psi_m(t)$ is a (positive) continuous function, then the
derivatives of order $m+1$ with respect to $x_1$, $x_2$, $x_3$ of the
$u_i(x,t)$ are strongly continuous in $t$ and satisfy the inequality

  $$ J_{m+1}(t) < A_m \int_0^t {\psi_m(t')\, dt' 
                \over \sqrt{\nu (t-t')}}
  \leqno(2.19)$$
Now suppose that the derivatives of order $m$ of the functions $X_i(x,t)$
with respect to $x_1$, $x_2$, $x_3$ exist, are smaller in absolute value
than a continuous function of $t$, and satisfy condition H with coefficient
$\theta_m(t)$. Then the derivatives of order $m+2$ of the 
$u_i(x,t)$ 
with respect to  $x_1$, $x_2$, $x_3$ 
exist, are uniformly continuous, and satisfy the inequality

\vfill
\eject

 % f215

\foddheader{215}

  $$  D_{m+2}(t) < A \int_0^t {\theta_m(t')\, dt' 
                \over [\nu (t-t')]^{3\over 4}}.
  \leqno(2.20)$$

\medskip

  {\bf 13. } {\it Cas G\'en\'eral.}

\smallskip

  Pour obtenir une solutions $u_i(x,t)$ de (2.1) correspondant 
\`a des valeurs initiales donn\'ees $u_i(x,0)$, il suffit d'ajouter
les deux  solutions particuli\`eres pr\'ec\'edentes, c'est-\`a-dire
de prendre:

  $$u_i(x,t) = u_i'(x,t)+u_i''(x,t);\,\, p(x,t) = p''(x,t).$$
Nous nous proposons de compl\'eter les renseignements que fournissent
les deux paragraphes pr\'ec\'edents en \'etablissant que $u_i(x,t)$ 
est fortement continue en $t$ et en majorant $W(t)$.

  Cette fort continuit\'e est \'evidente dans le cas o\`u $X_i(x,t)$ est 
nul hors d'un domaine $\varpi$; quand $x$ s'\'eloigne ind\'efiniment
$u_i''(x,t)$, ${\partial u_i''(x,t)\over x_k}$ et $p(x,t)$ tendent
alors vers z\'ero respectivement comme 
$(x_i x_i)^{-{3\over 2}}$, $(x_i x_i)^{-2}$ et $(x_i x_i)^{-1}$; et il
suffit d'int\'egrer les deux membres de l'\'egalit\'e:

  $$\nu u_i \Delta u_i - 
        {1\over 2}{\partial\over\partial t}(u_i u_i) -
        {1\over\rho}u_i{\partial p \over \partial x_i}
        = -u_iX_i$$
pour obtenir $<<$ {\it la r\'elation de dissipation de l'\'energie:} $>>$

  $$\nu\int_0^t J_1^2(t')\, dt'+{1\over 2}W(t)-{1\over 2}W(0)
      = \int_0^t\, dt'\intpi u_i(x,t')X_i(x,t')\dee{x}
  \leqno(2.21)$$
d'o\`u r\'esulte l'in\'egalit\'e:

  $${1\over 2}W(t) \leq {1\over 2}W(0)
                    +\int_0^t\, dt'\sqrt{W(t')}
                     \sqrt{\intpi X_i(x,t')X_i(x,t')\dee{x}}.$$
$W(t)$ est donc inf\'erieur ou \'egal \`a la solution $\lambda(t)$
de l'\'equation

  $${1\over 2}\lambda(t) = {1\over 2}W(0)
                    +\int_0^t\, dt'\sqrt{\lambda(t')}
                     \sqrt{\intpi X_i(x,t')X_i(x,t')\dee{x}}$$
c'est-\`a-dire:

  $$\sqrt{W(t)} \leq \int_0^t
                  \sqrt{\intpi X_i(x,t')X_i(x,t')\dee{x}}
                 \,dt' + \sqrt{W(0)}.
  \leqno(2.22)$$

\vfill
\eject

 % e215

\oddheader{215}

  $$  D_{m+2}(t) < A \int_0^t {\theta_m(t')\, dt' 
                \over [\nu (t-t')]^{3\over 4}}.
  \leqno(2.20)$$

\medskip

  {\bf 13. } {\it  General case.}

\smallskip

  To obtain solutions $u_i(x,t)$ of (2.1) corresponding to given initial
values $u_i(x,0)$, it suffices to superpose the two preceeding
particular solutions, taking

  $$u_i(x,t) = u_i'(x,t)+u_i''(x,t);\,\, p(x,t) = p''(x,t).$$
We propose to complete the information of the two preceeding paragraphs
by establishing that
$u_i(x,t)$  is strongly continuous in $t$ and is majorised by 
$W(t)$.

  This strong continuity is evident in the case where 
$X_i(x,t)$ is zero outside of a domain
$\varpi$. When $x$ moves indefinitely far away,  
$u_i''(x,t)$, ${\partial u_i''(x,t)\over x_k}$ and $p(x,t)$ approach zero 
as $(x_i x_i)^{-{3\over 2}}$, $(x_i x_i)^{-2}$ and  $(x_i x_i)^{-1}$ 
respectively, and it suffices to integrate  

  $$\nu u_i \Delta u_i - 
        {1\over 2}{\partial\over\partial t}(u_i u_i) -
        {1\over\rho}u_i{\partial p \over \partial x_i}
        = -u_iX_i$$
to obtain  {\it the relation of dissipation of energy} 

  $$\nu\int_0^t J_1^2(t')\, dt'+{1\over 2}W(t)-{1\over 2}W(0)
      = \int_0^t\, dt'\intpi u_i(x,t')X_i(x,t')\dee{x}
  \leqno(2.21)$$
from which we get the inequality

  $${1\over 2}W(t) \leq {1\over 2}W(0)
                    +\int_0^t\, dt'\sqrt{W(t')}
                     \sqrt{\intpi X_i(x,t')X_i(x,t')\dee{x}}.$$
$W(t)$ is therefore less than or equal to the solution $\lambda(t)$
of the equation

  $${1\over 2}\lambda(t) = {1\over 2}W(0)
                    +\int_0^t\, dt'\sqrt{\lambda(t')}
                     \sqrt{\intpi X_i(x,t')X_i(x,t')\dee{x}}$$
which is to say

  $$\sqrt{W(t)} \leq \int_0^t
                  \sqrt{\intpi X_i(x,t')X_i(x,t')\dee{x}}
                 \,dt' + \sqrt{W(0)}.
  \leqno(2.22)$$

\vfill
\eject

 % f216

\evenheader{216}

  Quand $X_i(x,t)$ n'est pas nul hors d'un domaine $\varpi$, on peut
approcher les fonctions $X_i(x,t)$ par une suite de fonctions $X_i^*(x,t)$
nulles hors de domaines $\varpi^*$, et par ce proc\'ed\'e \'etablir
ques les relations (2.21) et (2.22) sont encore valables. La relation
(2.21) prouve que $W(t)$ est continue; les fonctions $u_i(x,t)$
sont donc fortement continues en $t$ pour $t\geq 0$.

\medskip

  {\bf 14. } $u_i(x,t) = u_i'(x,t)+u_i''(x,t)$ est la seule solution
du probl\`eme pos\'e au parapraphe 10 pour laquelle $W(t)$ est 
inf\'erieure \`a une fonction continue de $t$;
cette proposition r\'esulte de la suivante

\medskip

  {\it Th\'eor\`eme d'unicit\'e:} Le syst\`eme

  $$ \nu\Delta u_i(x,t)
          -{\partial u_i(x,t)\over \partial t}
          -{1\over \rho}{\partial p(x,t)\over \partial x_i}
          = 0;\,\,
           {\partial u_j(x,t)\over \partial x_j}=0
  \leqno(2.23)$$
admet une seule solution, d\'efinie et continue pour $t\geq 0$, nulle
pour $t=0$, telle que $W(t)$ soit inf\'erieure \`a une fonction continue
de $t$; c'est $u_i(x,t)=0$.

\medskip

\ninerm

  En effet les fonctions

  $$v_i(x,t)=\int_0^t\overline{u_i(x,t')}\, dt', \quad
     q(x,t)=\int_0^t\overline{p(x,t')}\, dt'$$
constituent des solutions du m\^eme syst\`eme (2.23); les d\'eriv\'ees

  $${\partial^m v_i(x,t)\over
       \partial x_1^h \partial x_2^k \partial x_3^l}
       \quad {\rm et} \quad
      {\partial^{m+1} v_i(x,t)\over
       \partial t \partial x_1^h \partial x_2^k \partial x_3^l}$$
existent et sont continues; on a \'evidement $\Delta q=0$ et par suite:

  $$\nu\Delta\Delta v_i-{\partial\over\partial t}(\Delta v_i)=0;$$
la Th\'eorie de la Chaleur permet d'en d\'eduire $\Delta v_i=0$.
D'autre part les in\'egalit\'ees (1.2) et (1.21) prouvent que
l'int\'egral $\intpi v_i(x,t) v_i(x,t)\dee{x}$ est finie.
Donc $v_i(x,t)=0$. Et par suite $u_i(x,t)=0$.

\rm

\medskip

  Enon\c cons un corollaire qu'utilisera le paragraph suivant:

\smallskip

  {\it Lemme 8.} Supposons que nous ayons pour $\Theta\leq t < T$ le
syst\`eme de relations:

  $$\nu\Delta u_i(x,t)
          -{\partial u_i(x,t)\over \partial t}
          -{1\over \rho}{\partial p(x,t)\over \partial x_i}
          =
          -{\partial X_{ik}(x,t)\over \partial x_k};
  \,\, {\partial u_j(x,t)\over \partial x_j}=0. $$
Supposons les d\'eriv\'ees 
${\partial^2 X_{ik}(x,t)\over \partial x_j \partial x_l}$
continues et les int\'egrals

\vfill
\eject

 % e216

\evenheader{216}

  When $X_i(x,t)$ is not zero outside a domain $\varpi$, 
one can approach the functions
$X_i(x,t)$ by a sequence of functions $X_i^*(x,t)$
zero outside domains $\varpi^*$, and establish by the preceeding
that relations (2.21) and (2.22) still hold. Then
(2.21) shows that $W(t)$ is continuous. The $u_i(x,t)$
are therefore strongly continuous in $t$ for $t\geq 0$.

\medskip

  {\bf 14. } $u_i(x,t) = u_i'(x,t)+u_i''(x,t)$ is the only solution
to the problem posed in paragraph 10, for which $W(t)$ is less than a
continuous function of $t$. This proposition results from the following

\medskip

  {\it Uniqueness theorem} The system

  $$ \nu\Delta u_i(x,t)
          -{\partial u_i(x,t)\over \partial t}
          -{1\over \rho}{\partial p(x,t)\over \partial x_i}
          = 0;\,\,
           {\partial u_j(x,t)\over \partial x_j}=0
  \leqno(2.23)$$
has just one solution defined and continuous for $t\geq 0$, zero for 
$t=0$, such that $W(t)$ is less than a continuous function of
$t$. This solution is $u_i(x,t)=0$.

\medskip

\ninerm

  In fact the functions

  $$v_i(x,t)=\int_0^t\overline{u_i(x,t')}\, dt', \quad
     q(x,t)=\int_0^t\overline{p(x,t')}\, dt'$$
are solutions to the same system (2.23). The derivatives

  $${\partial^m v_i(x,t)\over
       \partial x_1^h \partial x_2^k \partial x_3^l}
       \quad {\rm and} \quad
      {\partial^{m+1} v_i(x,t)\over
       \partial t \partial x_1^h \partial x_2^k \partial x_3^l}$$
exist and are continuous. One evidently has $\Delta q=0$ and 
it follows that

  $$\nu\Delta\Delta v_i-{\partial\over\partial t}(\Delta v_i)=0.$$
The theory of heat allows us to deduce that $\Delta v_i=0$.
Further, inequalities (1.2) and (1.21) show that
the integral
$\intpi v_i(x,t) v_i(x,t)\dee{x}$ is finite. Therefore
$v_i(x,t)=0$. And then $u_i(x,t)=0$.

\rm

\medskip

 We state a corollary to be used in the following paragraph.

\smallskip

  {\it Lemma 8.} Suppose we have for $\Theta\leq t < T$ 
the system of relations

  $$\nu\Delta u_i(x,t)
          -{\partial u_i(x,t)\over \partial t}
          -{1\over \rho}{\partial p(x,t)\over \partial x_i}
          =
          -{\partial X_{ik}(x,t)\over \partial x_k};
  \,\, {\partial u_j(x,t)\over \partial x_j}=0. $$
Suppose the derivatives
${\partial^2 X_{ik}(x,t)\over \partial x_j \partial x_l}$
are continuous and the integrals

\vfill
\eject

% f217

\foddheader{217}

  $$\intpi X_{ik}(x,t)X_{ik}(x,t)\dee{x},
       \,\,\intpi u_i(x,t)u_i(x,t)\dee{x}$$
inf\'erieures \`a des fonctions de $t$ continues pour $\Theta\leq t < T$.
Nous avons alors:

  $$u_i(x,t)={1\over (2\sqrt{\pi})^3}
            \intpi 
               {e^{-{{\bf r}^2 \over 4\nu t}}\over
                   (\nu t)^{3\over 2}}u_i(y,t_0)
             \dee{y}+$$
  $$ {\partial\over\partial x_k}
                 \int_{t_0}^t\,dt'
                   \intpi 
                     T_{ij}(x-y,t-t')X_{jk}(y,t)
                   \dee{y};$$
  $$p(x,t) = -{\rho\over 4\pi}
               {\partial\over\partial x_k}
                     \intpi
                        {1\over r}X_{ik}(y,t)
                     \dee{y};
          \,\,(\Theta\leq t_0 < t < T).$$

\bigskip

\centerline{\bf III. Mouvements r\'eguliers.}

\medskip

  {\bf 15.} {\it D\'efinitions:} Les mouvements des liquides visqueux sont
r\'egis par les \'equations de Navier:

  $$\nu\Delta u_i(x,t)
          -{\partial u_i(x,t)\over \partial t}
          -{1\over \rho}{\partial p(x,t)\over \partial x_i}
          =
          u_k(x,t) {\partial u_i(x,t)\over \partial x_k};
        \,\, {\partial u_k(x,t)\over \partial x_k}=0,
  \leqno(3.1)$$
o\`u $r$ et $\rho$ sont des constantes, $p$ la pression, $u_i$ les
composantes de la vitesse. Nous poserons:

  $$W(t) = \intpi u_i(x,t)u_i(x,t)\dee{x},$$
  $$V(t) = \max\sqrt{u_i(x,t)u_i(x,t)}.$$

  Nous dirons qu'{\it une solution $u_i(x,t)$ de ce syst\`eme est r\'eguliere}
dans un intervalle de temps$^1$ $\Theta < t < T$ si dans cet intervalle
de temps les fonctions $u_i$, la fonction $p$ correspondante et les
d\'eriv\'ees ${\partial u_i\over \partial x_k}$, 
${\partial^2 u_i\over \partial x_k \partial x_l}$,
${\partial u_i\over \partial t}$, ${\partial p\over \partial x_i}$ sont
 continues
par rapport \`a $x_1$, $x_2$, $x_3$, $t$ et si en outre les fonctions
$W(t)$ et $V(t)$ sont inf\'erieures \`a des fonctions de $t$ continues 
pour $\Theta < t < T $.

  Nous utiliserons les conventions suivantes:

  La fonction $D_m(t)$ sera d\'efinie pour chaque valeur de $t$ au 
voisinage de

\footrule

  $^1$ Le cas o\`u $T=+\infty$ n'est pas exclu.

\rm
\vfill
\eject

 % e217

\oddheader{217}

  $$\intpi X_{ik}(x,t)X_{ik}(x,t)\dee{x},
       \,\,\intpi u_i(x,t)u_i(x,t)\dee{x}$$
less than some continuous functions of $t$ for $\Theta\leq t < T$.
We have then

  $$u_i(x,t)={1\over (2\sqrt{\pi})^3}
            \intpi 
               {e^{-{{\bf r}^2 \over 4\nu t}}\over
                   (\nu t)^{3\over 2}}u_i(y,t_0)
             \dee{y}+$$
  $$ {\partial\over\partial x_k}
                 \int_{t_0}^t\,dt'
                   \intpi 
                     T_{ij}(x-y,t-t')X_{jk}(y,t)
                   \dee{y};$$
  $$p(x,t) = -{\rho\over 4\pi}
               {\partial\over\partial x_k}
                     \intpi
                        {1\over r}X_{ik}(y,t)
                     \dee{y};
          \,\,(\Theta\leq t_0 < t < T).$$

\bigskip

\centerline{\bf III. Regular motions.}

\medskip

  {\bf 15.} {\it Definitions:} Motions of viscous liquids
are governed by Navier's equations

  $$\nu\Delta u_i(x,t)
          -{\partial u_i(x,t)\over \partial t}
          -{1\over \rho}{\partial p(x,t)\over \partial x_i}
          =
          u_k(x,t) {\partial u_i(x,t)\over \partial x_k};
        \,\, {\partial u_k(x,t)\over \partial x_k}=0,
  \leqno(3.1)$$
where $\nu$ and $\rho$ are constants, $p$ is the pressure, $u_i$ the
components of the velocity. We set

  $$W(t) = \intpi u_i(x,t)u_i(x,t)\dee{x},$$
  $$V(t) = \max\sqrt{u_i(x,t)u_i(x,t)}.$$

  We say that {\it a solution $u_i(x,t)$ of this system is regular}
in an interval of time$^1$ $\Theta < t < T $ if in this interval
the functions  $u_i$, the corresponding $p$ and the derivatives
${\partial u_i\over \partial x_k}$, 
${\partial^2 u_i\over \partial x_k \partial x_l}$,
${\partial u_i\over \partial t}$, ${\partial p\over\partial x_i}$ 
are continuous with respect to
$x_1$, $x_2$, $x_3$, $t$
and if in addition the functions
$W(t)$ and $V(t)$
are less than some continuous functions
of $t$ for $\Theta < t < T $.

  We use the following conventions.

  The function $D_m(t)$ will be defined for each value of $t$ in
a neighborhood in 

\footrule

  $^1$ The case where $T=+\infty$ is not excluded.

\rm
\vfill
\eject

 % f218

\evenheader{218}

\noindent laquelle les d\'eriv\'ees $\partm$ existent et sont uniform\'ement
continues en $t$; elle sera \'egale \`a la borne sup\'erieure de leurs
valeurs absolues.

  La fonction $C_0(t)$ [ou $C_m(t)$] ser d\'efinie pour toutes les valeurs
de $t$ au voisinage desquelles les fonctions $u_i(x,t)$
[ ou les d\'eriv\'ees $\partm$ ] v\'erifient 
une m\^eme condition H; elle
en sera le coefficient.

  Enfin la fonction $J_m(t)$ sera d\'efinie pour chaque valeur de $t$ au
voisinage de Laquelle les d\'eriv\'ees $\partm$ existent et sont fortement
continues en $t$; nous poserons:

  $$J_m^2(t)=\intpi{\partial^m u_i(x,t)\over
               \partial x_k \partial x_l \ldots}
                \,\,{\partial^m u_i(x,t)\over
               \partial x_k \partial x_l \ldots} \dee{x}.$$

  Le lemme 8 (p. 216) s'applique aux solutions r\'eguli\`eres du 
syst\`eme (3.1) et nous donne les relations:

  $$u_i(x,t)={1\over (2\sqrt{\pi})^3}
            \intpi 
               {e^{-{{\bf r}^2 \over 4\nu t}}\over
                   (\nu t)^{3\over 2}}u_i(y,t_0)
             \dee{y}+
  \leqno(3.2)$$
  $$ {\partial\over\partial x_k}
                 \int_{t_0}^t\,dt'
                   \intpi 
                     T_{ij}(x-y,t-t')u_j(y,t')u_k(y,t')
                   \dee{y};$$
  $$p(x,t) = {\rho\over 4\pi}
               {\partial^2\over \partial x_k \partial x_j}
                     \intpi
                        {1\over r}u_k(y,t)u_j(y,t)
                     \dee{y};
          \,\,(\Theta < t_0 < t < T).
  \leqno(3.3)$$

  Les paragraphs 11 et 12 permettent de d\'eduire de la relation (3.2)
les faits suivants: les fonctions $u_i(x,t)$ sont uniform\'ement et
fortement continues en $t$ pour $\Theta <  t < T$; la fonction $C_0(t)$
est d\'efinie pour $\Theta <  t < T$ et l'on a [cf. (2.7) et (2.18)]:

  $$C_0(t) < {A\sqrt{W(t_0)} \over \nu(t-t_0)}
                  +
              A\int_{t_0}^t {V^2(t')\,dt'\over
                                      [\nu(t-t')]^{3\over 4}}.$$
Ce r\'esultat port\'e dans (3.2) prouve que $D_1(t)$ existe pour
$\Theta <  t < T$ et fournit l'in\'egalit\'e [cf. (2.7) et (2.16)]:

  $$D_1(t) < {A\sqrt{W(t_0)} \over [\nu(t-t_0)]^{5\over 4}}
                  +
              A\int_{t_0}^t {V(t')C_0(t')\,dt'\over
                                      [\nu(t-t')]^{3\over 4}}.$$

\vfill
\eject

 % e218

\evenheader{218}

\noindent which the derivatives
$\partm$
exist and are uniformly continuous in 
$t$; it will be the upper bound of their absolute values.

  The function $C_0(t)$ [or $C_m(t)$] will be defined for all values of $t$ 
in a neighborhood in which the functions $u_i(x,t)$
[ or the derivatives $\partm$ ] satisfy the same condition H; it will be the
coefficient.

  Finally the function $J_m(t)$ will be defined for each value of $t$ in a 
neighborhood
in which the derivatives $\partm$ exist and are strongly continuous in
$t$. We set

  $$J_m^2(t)=\intpi{\partial^m u_i(x,t)\over
               \partial x_k \partial x_l \ldots}
                \,\,{\partial^m u_i(x,t)\over
               \partial x_k \partial x_l \ldots} \dee{x}.$$

  Lemma 8 (p. 216) applies to regular solutions to system (3.1) and gives us
the relations

  $$u_i(x,t)={1\over (2\sqrt{\pi})^3}
            \intpi 
               {e^{-{{\bf r}^2 \over 4\nu t}}\over
                   (\nu t)^{3\over 2}}u_i(y,t_0)
             \dee{y}+
  \leqno(3.2)$$
  $$ {\partial\over\partial x_k}
                 \int_{t_0}^t\,dt'
                   \intpi 
                     T_{ij}(x-y,t-t')u_j(y,t')u_k(y,t')
                   \dee{y};$$
  $$p(x,t) = {\rho\over 4\pi}
               {\partial^2\over \partial x_k \partial x_j}
                     \intpi
                        {1\over r}u_k(y,t)u_j(y,t)
                     \dee{y};
          \,\,(\Theta < t_0 < t < T).
  \leqno(3.3)$$

  Paragraphs 11 and 12 allow us to conclude 
from (3.2) that
the functions $u_i(x,t)$ are uniformly and strongly continuous in $t$ 
for $\Theta <  t < T$,  $C_0(t)$
is defined for $\Theta <  t < T$ and we have [cf. (2.7) and (2.18)]

  $$C_0(t) < {A\sqrt{W(t_0)} \over \nu(t-t_0)}
                  +
              A\int_{t_0}^t {V^2(t')\,dt'\over
                                      [\nu(t-t')]^{3\over 4}}.$$
This result with (3.2) shows that $D_1(t)$ exists for 
$\Theta <  t < T$ and gives the inequality [cf. (2.7) and (2.16)]

  $$D_1(t) < {A\sqrt{W(t_0)} \over [\nu(t-t_0)]^{5\over 4}}
                  +
              A\int_{t_0}^t {V(t')C_0(t')\,dt'\over
                                      [\nu(t-t')]^{3\over 4}}.$$

\rm

\vfill
\eject

 % f219

\foddheader{219}

   Poursuivons par r\'ecurrence:

  L'existence de $D_1(t),\dots,D_{m+1}(t)$ assure celle de $C_{m+1}(t)$
et l'on a [cf. (2.7) et (2.18)]:

  $$C_{m+1}(t) <
      {A_m\sqrt{W(t_0)}
                 \over
        [\nu(t-t_0)]^{m+3\over 2}} +
      A_m\int_{t_0}^t
            {V(t')D_{m+1}(t')+\sum_{\alpha+\beta=m+1}D_\alpha(t')D_\beta(t')
             \over
            [\nu(t-t')]^{3\over 4}}
         \,dt'.$$

  L'existence de  $D_1(t),\dots,D_{m+1}(t),C_0(t),\dots,C_{m+1}(t)$ assure
celle de $D_{m+2}(t)$ et l'on peut majorer cette derni\`ere fonction \`a 
l'aide des pr\'ec\'edentes [cf. (2.7) et (2.20)].

  Les fonctions $D_m(t)$ et $C_m(t)$ sont donc d\'efinies pour
$\Theta <  t < T$, quelque grand que soit $m$.

  D'autre part les paragraphes 11 et 12 permettent de d\'eduire de (3.2)
l'existence de $J_1(t)$ pour toutes ces valeurs de $t$; et nous avons 
[cf. (2.8 et (2.19)]:

  $$J_1(t) <      
       {A\sqrt{W(0)}
                 \over
        [\nu(t-t_0)]^{1\over 2}} +
      A\int_{t_0}^t
            {W(t')D_1(t')
             \over
            \sqrt{\nu(t-t')}}
         \,dt'.$$
Plus g\'en\'eralement l'existence de $D_1(t),\dots,D_m(t)$, 
$J_1(t),\dots,J_{m-1}(t)$ assure celle de $J_m(t)$ [cf. (2.8) et (2.19)].

  Il nous est maintenant ais\'e d'\'etablir par l'interm\'ediaire de (3.3)
que la fonction $p(x,t)$ et ses d\'eriv\'ees 
${\partial^m p(x,t) \over \partial x_k \partial x_j \ldots }$ sont
uniform\'ement et fortement continues en $t$ pour $\Theta <  t < T$.
D'apr\`es les \'equations de Navier il en est de m\^eme pour les fonctions
${\partial u_i \over \partial t}$,
${\partial^{m+1} u \over \partial t \partial x_k \partial x_j \ldots }$.

  Plus g\'en\'eralement les \'equations (3.1) et (3.3) permettent de
ramener l'\'etude de d\'eriv\'ees qui sont d'ordre $n+1$ par rapport
\`a $t$ \`a l'\'etude des d\'eriv\'ees qui sont d'ordre $n$ par
rapport \`a $t$. Et l'on aboutit finalement au {\it th\'eoreme} suivante:

\medskip

  {\it Si les fonctions $u_i(x,t)$ constituent une solution des \'equations
de Navier r\'eguli\`ere pour $\Theta <  t < T$, alors toutes leurs 
d\'eriv\'ees partielles existent; ces d\'eriv\'ees partielles et les 
fonctions $u_i(x,t)$ elles-m\^emes sont uniform\'ement et fortement
continues en $t$ pour $\Theta <  t < T$.

\vfill
\eject

 % e219

\oddheader{219}

   We proceed by recurrence:

  The existence of $D_1(t),\dots,D_{m+1}(t)$ guarantees that of $C_{m+1}(t)$
and one has [cf. (2.7) and (2.18)]

  $$C_{m+1}(t) <
      {A_m\sqrt{W(t_0)}
                 \over
        [\nu(t-t_0)]^{m+3\over 2}} +
      A_m\int_{t_0}^t
            {V(t')D_{m+1}(t')+\sum_{\alpha+\beta=m+1}D_\alpha(t')D_\beta(t')
             \over
            [\nu(t-t')]^{3\over 4}}
         \,dt'.$$

  The existence of $D_1(t),\dots,D_{m+1}(t),C_0(t),\dots,C_{m+1}(t)$ 
guarantees that of $D_{m+2}(t)$ and one can majorize this last function
using the preceeding
[cf. (2.7) and (2.20)].

  The functions $D_m(t)$ and $C_m(t)$ are therefore defined
for
$\Theta <  t < T$, however large  $m$ may be.

  Further, paragraphs 11 and 12 allow us to deduce from (3.2)
the existence of $J_1(t)$ 
for all values of $t$ and we have
[cf. (2.8) and  (2.19)]

  $$J_1(t) <      
       {A\sqrt{W(0)}
                 \over
        [\nu(t-t_0)]^{1\over 2}} +
      A\int_{t_0}^t
            {W(t')D_1(t')
             \over
            \sqrt{\nu(t-t')}}
         \,dt'.$$
More generally the existence of $D_1(t),\dots,D_m(t)$, 
$J_1(t),\dots,J_{m-1}(t)$ guarantees that of $J_m(t)$ [cf. (2.8) and (2.19)].

  It is now easy for us to establish by the intermediary (3.3) that 
$p(x,t)$  and its derivatives
${\partial^m p(x,t) \over \partial x_k \partial x_j \ldots }$ 
are uniformly and strongly continuous in
$t$ for $\Theta <  t < T$.
By Navier's equations it is the same for the functions
${\partial u_i \over \partial t}$,
${\partial^{m+1} u \over \partial t \partial x_k \partial x_j \ldots }$.

  More generally, equations (3.1) and (3.3) allow us to
reduce the study of the order $n+1$ derivatives with respect to $t$
to the study of the order $n$ derivatives with respect to $t$.
So finally we achieve the following {\it theorem}.

\medskip

  {\it If the functions $u_i(x,t)$ are a regular solution
of Navier's equations for $\Theta <  t < T$, then all their
partial derivatives exist, and the derivatives as well as the
$u_i(x,t)$ are uniformly and strongly continuous in 
$t$ for $\Theta <  t < T$.

\vfill
\eject

% f220

\evenheader{220}

  {\bf 16.} Le paragraphe pr\'ec\'edent nous apprend plus: il nous 
apprend \`a majorer les fonctions $u_i(x,t)$ et leurs d\'eriv\'ees
partielles de tous les ordres au moyen des seules fonctions $W(t)$
et $V(t)$. Il en r\'esulte:

\medskip

  {\it Lemme 9.} Soit une infinit\'e de solutions des \'equations
de Navier, $u_i^*(x,t)$, r\'eguli\`eres dans un m\^eme intervalle de
temps $(\Theta, T)$. Supposons les diverses fonctions
$V^*(t)$ et $W^*(t)$ inf\'erieures \`a une m\^eme fonction de $t$,
continue dans $(\Theta, T)$. De cette infinit\'e de solutions on peut
alors extraire une suite partielle telle que les fonctions $u_i^*(x,t)$
de cette suite et chaqune de leurs d\'eriv\'ees convergent
respectivement vers certaines fonctions $u_i(x,t)$ et vers leurs
d\'eriv\'ees. Chacune de ces convergences est uniforme sur tout 
domaine $\varpi$ pour $\Theta+\eta < t < T-\eta \,\, (\eta > 0)$.
Les fonctions $u_i(x,t)$ constituent une solution des \'equations
de Navier r\'eguli\`ere dans $(\Theta, T)$.

\medskip

\ninerm

   En effet le Proc\'eed\'e diagonal de Cantor (\S 4, p. 201) permet
d'extraire une suite de fonctions $u_i^*(x,t)$ telle que ces fonctions
$u_i^*(x,t)$ et leurs d\'eriv\'ees convergent pour tous les 
syst\`emes rationnels de valeurs donn\'ees \`a $x_1$, $x_2$, $x_3$, $t$.
Cette suite partielle poss\`ede les propri\'et\'es qu'\'enonce le lemme.

\rm

\medskip

  {\bf 17. } La quantit\'e $W(t)$ et la quantit\'e $J_1(t)$ --- que
d\'esormais nous d\'esignerons pour simplifier par $J(t)$ --- sont li\'ees
par une relation importante; elle s'obtient en rempla\c cant dans (2.21)
$X_i$ par $u_k {\partial u_i\over\partial x_k}$ et en remarquant que:

  $$\intpi u_i(x,t')u_k(x,t'){\partial u_i(x,t')\over\partial x_k}\dee{x}
      =
    {1\over  2}
     \intpi u_k(x,t'){\partial u_i(x,t')u_i(x,t') \over\partial x_k}\dee{x}
      = 0;$$
c'est $<<$ {\it la relation de la dissipation de l'\'energie} $>>$:

  $$\nu\int_{t_0}^t J^2(t')\,dt' + {1\over  2}W(t)={1\over  2}W(t_0).
  \leqno(3.4)$$
Cette relation et les deus paragraphes ci-dessus prouvent que les fonctions
$W(t)$, $V(t)$, et $J(t)$ jouent un r\^ole essentiel. Aussi retiendrons-nous
de toutes les in\'egalit\'es qu'on peut d\'eduire du chapitre II uniquement
quelques-unes o\`u figurent ces trois fonctions, sans plus nous occuper des 
quantit\'es $C_m(t)$, $D_m(t)$, $\dots$

  Avant d'\'ecrire ces quelques in\'egalit\'es fondamentales posons 
{\it une d\'efinition: }

  {\it Une solution $u_i(x,t)$ des \'equations de Navier sera dite 
r\'eguli\`ere pour $\Theta\leq t < T$ } quand elle sera r\'eguli\`ere
pour $\Theta < t < T$ et qu'en outre les circonstances suivantes

\vfill
\eject

 % e220

\evenheader{220}

  {\bf 16.} The preceeding paragraph teaches us more: we learn that
it is possible to bound the functions 
$u_i(x,t)$ and their partial derivatives of all orders by means of
just
$W(t)$
and $V(t)$. The result is:

\medskip

  {\it Lemma 9.} Let $u_i^*(x,t)$ be an infinity of solutions to 
Navier's equations, all regular in the same interval $(\Theta, T)$.
Suppose the various $V^*(t)$ and $W^*(t)$ all less than one function 
of $t$, continuous in $(\Theta, T)$. Then one can extract a subsequence
such that the $u_i^*(x,t)$ and each of their derivatives converge
respectively to certain functions $u_i(x,t)$ and their derivatives.
Each of the convergences is uniform on all domains $\varpi$ for
 $\Theta+\eta < t < T-\eta \,\, (\eta > 0)$. The functions
 $u_i(x,t)$ are a regular solution of Navier's equations
in $(\Theta, T)$.

\medskip

\ninerm

  In fact, Cantor's diagonal method (\S 4, p. 201) allows the extraction
of a sequence of functions
$u_i^*(x,t)$ which, with their derivatives, converge for any given rational
values of  
$x_1$, $x_2$, $x_3$, $t$.
This subsequence has the properties stated in the lemma.

\rm

\medskip

  {\bf 17. } The quantities $W(t)$ and $J_1(t)$ --- which from now
on we write simply as $J(t)$ --- 
are linked by an important relation. It is obtained by replacing
$X_i$ in (2.21) by  $u_k {\partial u_i\over\partial x_k}$ and remarking
that

  $$\intpi u_i(x,t')u_k(x,t'){\partial u_i(x,t')\over\partial x_k}\dee{x}
      =
    {1\over  2}
     \intpi u_k(x,t'){\partial u_i(x,t')u_i(x,t') \over\partial x_k}\dee{x}
      = 0;$$
It is the {\it ``energy dissipation relation'' }

  $$\nu\int_{t_0}^t J^2(t')\,dt' + {1\over  2}W(t)={1\over  2}W(t_0).
  \leqno(3.4)$$
This relation and the two paragraphs above show that the functions
$W(t)$, $V(t)$, and $J(t)$ play an essential role. We will  especially
point out, of 
all the inequalities one can deduce from chapter II, some which involve
these three functions without any longer occupying ourselves with the
quantities
$C_m(t)$, $D_m(t)$, $\dots$

  Before writing the fundamental inequalities, we make the
{\it definition: }

  {\it A solution $u_i(x,t)$ of Navier's equations  will be called
regular for $\Theta\leq t < T$ } when it is regular for $\Theta < t < T$ 
and if in addition the following  conditions

\vfill
\eject

 % f221

\foddheader{221}

\noindent se pr\'esenteront: les fonctions $u_i(x,t)$ et
$ {\partial u_i(x,t)\over\partial x_j}$ sont continues par rapport aux
variables $x_1$, $x_2$, $x_3$, $t$ m\^eme pour $t=\Theta$; elles sont 
fortement continues en $t$ m\^eme pour $t=\Theta$; les fonctions $u_i(x,t)$
restent born\'ees quand $t$ tend vers $\Theta$.

  Dans ces conditions la relation (3.2) vaut pour $\Theta\leq t_0 < t < T$
(la valeur $\Theta$ \'etait jusq'\`a pr\'esent interdite \`a $t_0$); le 
chapitre II permet d'en d\'eduire {\it deux in\'egalit\'es fondamentales;}
ce sont, le symbole $\{B; C\}$ nous servant \`a repr\'esenter la plus petite
 des deux quantit\'es $B$ et $C$; $A'$, $A''$, $A'''$ \'etant des constantes
num\'eriques:

  $$V(t)<A'
        \int_{t_0}^t  \{
            {V^2(t')\over\sqrt{\nu(t-t')}};
            {W(t')\over [\nu(t-t')]^2}
           \}\,dt'
        +
        \{V(t_0); {A'''J(t_0)\over [\nu(t-t_0)]^{1\over 4}}\}
  \leqno(3.5)$$

  $$ J(t)< A''
        \int_{t_0}^t 
            {J(t')V(t')\over\sqrt{\nu(t-t')}}
        \,dt'
        +
      J(t_0) \qquad (\Theta\leq t_0 < t < T).
  \leqno(3.6)$$

\medskip

  {\bf 18. } {\it Comparaison de deux solutions r\'eguli\`eres. }

\smallskip

  Consid\'erons deux solutions des \'equations de Navier, $u_i$ et 
$u_i+v_i$, r\'eguli\`eres pour $\Theta < t < T$. Nous avons:

  $$ \nu\Delta v_i - {\partial v_i\over\partial t}
             -{1\over\rho}{\partial q\over\partial x_i}
    =
    v_k{\partial u_i\over\partial x_k}
    +
    (u_k+v_k){\partial v_i\over\partial x_k};\quad 
        {\partial v_k\over\partial x_k} = 0.$$
Posons:

  $$w(t)=\intpi v_i(x,t)v_i(x,t)\dee{x}; \qquad 
          j^2(t)=\intpi
             {\partial v_i(x,t)\over\partial x_k}
             {\partial v_i(x,t)\over\partial x_k}
                \dee{x}.$$
Appliquons la relation (2.21) qui nous a d\'ej\`a fourni la relation
fondamentale (3.4); elle donne ici:

  $$\nu j^2(t) + {1\over 2}{dw\over dt}
   =
    \intpi v_i v_k {\partial u_i\over\partial x_k}\dee{x}
    +
    \intpi v_i (u_k+v_k) {\partial v_i\over\partial x_k}\dee{x}.$$
Or nous avons:

  $$\intpi v_i (u_k+v_k){\partial v_i\over\partial x_k}\dee{x}
    =
   {1\over 2}\intpi 
      (u_k+v_k)
      {\partial (v_i v_i) \over\partial x_k}\dee{x}
    =0;$$
et

  $$\intpi v_i v_k {\partial u_i\over\partial x_k}\dee{x}
   =
    -\intpi {\partial v_i\over\partial x_k}v_k u_i \dee{x}
   <
   j(t)\sqrt{w(t)}V(t).$$

\vfill
\eject

 % e221

\oddheader{221}

\noindent are satisfied:
The functions  $u_i(x,t)$ and  $ {\partial u_i(x,t)\over\partial x_j}$ 
are continuous with respect to 
$x_1$, $x_2$, $x_3$, $t$ also for $t=\Theta$, they are strongly
continuous in $t$ also for 
$t=\Theta$, and the $u_i(x,t)$
remain bounded when $t$ approaches $\Theta$.

  In these conditions the relation (3.2) holds for $\Theta\leq t_0 < t < T$
(the value  $\Theta$ was not allowed to be $t_0$ until now). Chapter II
allows us to deduce {\it two fundamental inequalities.}
In these, the symbol 
$\{B; C\}$ is the smaller of  $B$ and $C$, and $A'$, $A''$, $A'''$ are
numerical constants. The inequalities are

  $$V(t)<A'
        \int_{t_0}^t  \{
            {V^2(t')\over\sqrt{\nu(t-t')}};
            {W(t')\over [\nu(t-t')]^2}
           \}\,dt'
        +
        \{V(t_0); {A'''J(t_0)\over [\nu(t-t_0)]^{1\over 4}}\}
  \leqno(3.5)$$

  $$ J(t)< A''
        \int_{t_0}^t 
            {J(t')V(t')\over\sqrt{\nu(t-t')}}
        \,dt'
        +
      J(t_0) \qquad (\Theta\leq t_0 < t < T).
  \leqno(3.6)$$

\medskip

  {\bf 18. } {\it Comparison of two regular solutions. }

\smallskip

 We consider two solutions of Navier's equations, $u_i$ and 
$u_i+v_i$, regular for $\Theta < t < T$. We have

  $$ \nu\Delta v_i - {\partial v_i\over\partial t}
             -{1\over\rho}{\partial q\over\partial x_i}
    =
    v_k{\partial u_i\over\partial x_k}
    +
    (u_k+v_k){\partial v_i\over\partial x_k};\quad 
        {\partial v_k\over\partial x_k} = 0.$$
Let

  $$w(t)=\intpi v_i(x,t)v_i(x,t)\dee{x}; \qquad 
          j^2(t)=\intpi
             {\partial v_i(x,t)\over\partial x_k}
             {\partial v_i(x,t)\over\partial x_k}
                \dee{x}.$$
We apply (2.21) which has already given us the fundamental relation
(3.4). Here it gives

  $$\nu j^2(t) + {1\over 2}{dw\over dt}
   =
    \intpi v_i v_k {\partial u_i\over\partial x_k}\dee{x}
    +
    \intpi v_i (u_k+v_k) {\partial v_i\over\partial x_k}\dee{x}.$$
Now we have

  $$\intpi v_i (u_k+v_k){\partial v_i\over\partial x_k}\dee{x}
    =
   {1\over 2}\intpi 
      (u_k+v_k)
      {\partial (v_i v_i) \over\partial x_k}\dee{x}
    =0$$
and

  $$\intpi v_i v_k {\partial u_i\over\partial x_k}\dee{x}
   =
    -\intpi {\partial v_i\over\partial x_k}v_k u_i \dee{x}
   <
   j(t)\sqrt{w(t)}V(t).$$

\vfill
\eject

 % f222

\evenheader{222}

\noindent Donc:

  $$ \nu j^2(t) + {1\over 2}{dw\over dt} 
    <
     j(t)\sqrt{w(t)}V(t)$$
d'o\`u:

  $$ 2\nu{dw\over dt} < w(t)V^2(t)$$
et finalement:

  $$ w(t)<w(t_0)e^{{1\over 2\nu}\int_{t_0}^t V^2(t')\,dt'}
     \qquad (\Theta < t_0 < t < T).
  \leqno(3.7)$$

  De cette relation importante r\'esulte en particulier:

  {\it Une th\'eor\`eme d'unicit\'e: } Deux solutions des \'equations
de Navier r\'eguli\`eres pour $\Theta \leq t < T$ sont n\'ecessairement
identiques pour ces valeurs dt $t$ si leurs \'etats de vitesse le sont
pour $t=\Theta$.

\medskip

  {\bf 19. } Donnons-nous {\it un \'etat initial r\'eguli\`er, } c'est-\`a-dire
un vecteur de divergence nulle, $u_i(x,t)$, continu, ainsi que les
d\'eriv\'ees premi\`eres de ses composantes, et tel que les quantit\'es
$W(0)$, $V(0)$, $J(0)$ soient finies. Le but de ce paragraphe est
d'\'etablir la proposition suivante:

\medskip

  {\it Th\'eor\`eme d'existence: A tout \'etat initial r\'eguli\`er,
$u_i(x,0)$, correspond une solution des \'equations
de Navier, $u_i(x,t)$, qui est d\'efinie pour des valeurs
$0\leq t < \tau$ de $t$ et qui se r\'eduit \`a $u_i(x,0)$ pour $t=0$.}

\medskip

  Formons les {\it approximations successives: }

  $$u_i^{(0)}(x,t)={1\over (2\sqrt{\pi})^3}
            \intpi 
           {e^{-{{\bf r}^2\over 4\nu t}}\over (\nu t)^{{3\over 2}}}
           u_i(y,0)\dee{y},$$

  $$\ldots\ldots\ldots\ldots\ldots\ldots\ldots\ldots\ldots\ldots\ldots\ldots$$

  $$u_i^{(n+1)}(x,t)={\partial\over\partial x_k}
         \int_0^t\,dt'\,
         \intpi
         T_{ij}(x-y,t-t')u_k^{(n)}(y,t')u_j^{(n)}(y,t')
         \dee{y}
       +
      u_i^{(0)}(x,t),$$

  $$\ldots\ldots\ldots\ldots\ldots\ldots\ldots\ldots\ldots\ldots\ldots\ldots$$

  Ecrivons en premier lieu les in\'egalit\'es, d\'eduites de (2.3)
et (2.13):

  $$V^0(t) \leq V(0)$$

  $$V^{(n+1)}(t)
      \leq
     A'\int_0^t {[V^{(n)}(t')]^2
                   \over
                 \sqrt{\nu(t-t')}}\,dt'
     + V(0);$$

\vfill
\eject

 % e222

\evenheader{222}

\noindent Therefore 

  $$ \nu j^2(t) + {1\over 2}{dw\over dt} 
    <
     j(t)\sqrt{w(t)}V(t)$$
from which

  $$ 2\nu{dw\over dt} < w(t)V^2(t)$$
and finally

  $$ w(t)<w(t_0)e^{{1\over 2\nu}\int_{t_0}^t V^2(t')\,dt'}
     \qquad (\Theta < t_0 < t < T).
  \leqno(3.7)$$

  From this important relation we get in particular

  {\it A uniqueness theorem: } Two regular solutions of Navier's 
equations for
$\Theta \leq t < T$ are necessarily identical for these  $t$ if 
their initial velocities are the same for $t=\Theta$.

\medskip

  {\bf 19.} Suppose we are given a {\it regular initial state,} which
is to say a continuous vector $u_i(x,t)$ with continuous first derivatives,
having zero divergence
and such that the quantities $W(0)$, $V(0)$, $J(0)$ are finite. 
The goal of this paragraph is to establish the following proposition.

\medskip

  {\it Existence theorem: To each regular initial state
$u_i(x,0)$ there corresponds a solution $u_i(x,t)$ to  Navier's equations,
defined for 
$0\leq t < \tau$ and which reduces to  $u_i(x,0)$ for $t=0$.}

\medskip

  We form {\it successive approximations}

  $$u_i^{(0)}(x,t)={1\over (2\sqrt{\pi})^3}
            \intpi 
           {e^{-{{\bf r}^2\over 4\nu t}}\over (\nu t)^{{3\over 2}}}
           u_i(y,0)\dee{y},$$

  $$\ldots\ldots\ldots\ldots\ldots\ldots\ldots\ldots\ldots\ldots\ldots\ldots$$

  $$u_i^{(n+1)}(x,t)={\partial\over\partial x_k}
         \int_0^t\,dt'\,
         \intpi
         T_{ij}(x-y,t-t')u_k^{(n)}(y,t')u_j^{(n)}(y,t')
         \dee{y}
       +
      u_i^{(0)}(x,t),$$

  $$\ldots\ldots\ldots\ldots\ldots\ldots\ldots\ldots\ldots\ldots\ldots\ldots$$

  First we write inequalities which follow from (2.3) and  (2.13):

  $$V^0(t) \leq V(0)$$

  $$V^{(n+1)}(t)
      \leq
     A'\int_0^t {[V^{(n)}(t')]^2
                   \over
                 \sqrt{\nu(t-t')}}\,dt'
     + V(0).$$

\vfill
\eject

% f223

\foddheader{223}

\noindent elles prouvent que nous avons, quel que soit $n$:

  $$ V^{(n)}(t) \leq \varphi(t) \qquad{\rm pour \ } 0\leq t\leq \tau,$$
si $\varphi(t)$ est une fonctions continue qui v\'erifie pour ces
valeurs de $t$ l'in\'egalit\'e:

  $$\varphi(t)\geq A'\int_0^t 
        {\varphi^2(t')\over\sqrt{\nu(t-t')}}
         \,dt' + V(0);$$
nous choisirons $\varphi(t)=(1+A)V(0)$; la valeur \`a donner \`a $\tau$ est:

  $$\tau = A\nu V^{-2}(0).
  \leqno (3.8)$$

  Posons alors:

  $$v^{(n)}(t)=\max\sqrt{[u_i^{(n)}(x,t)-u_i^{(n+1)}(x,t)]
                          [u_i^{(n)}(x,t)-u_i^{(n+1)}(x,t)]}
       $$
\`a l'instant $t$.

  Nous avons:
  $$v^{(1)}(t)
     < 
     A'\int_0^\tau {V^2(0)
                   \over
                 \sqrt{\nu(\tau-t')}}\,dt'
     = A V(0);$$
  $$v^{(n+1)}(t)
      <
     A\int_0^\tau {\varphi(t')v^{(n)}(t')
                   \over
                 \sqrt{\nu(\tau-t')}}\,dt'
     = A V(0)\int_0^\tau {v^{(n)}(t')
                   \over
                 \sqrt{\nu(\tau-t')}}\,dt';$$
d'o\`u r\'esulte que, pour $0\leq t\leq\tau$, les fonctions $u_i^{(n)}(x,t)$
convergent uniform\'ement vers des limites continues, $u_i(x,t)$.

  On d\'emontre sans difficult\'e qu'\`a l'int\'erieur de l'intervalle
$(0,\tau)$ chacune des d\'eriv\'ees des fonctions $u_i^{(n)}(x,t)$
converge uniform\'ement vers la d\'eriv\'ee correspondante des fonctions 
$u_i(x,t)$; les raisonnements sont trop proches de ceux du paragraphe 15
 pour que nous les reproduisions. Les fonctions $u_i(x,t)$ satisfont
donc les \'equations de Navier pour $0<t<\tau$.

\smallskip

  V\'erifions que l'int\'egrale 
$W(t)=\intpi u_i(x,t)u_i(x,t)\dee{x}$ est 
inf\'erieure \`a une fonction continue de $t$: les in\'egalit\'es 
(2.5) et (2.12) fournissent les suivantes, o\`u $A_0$ repr\'esente une
constante:

  $$\sqrt{W^{(0)}(t)}
             \leq
           \sqrt{W(0)}$$
  $$\sqrt{W^{(n+1)}(t)}
             \leq
      A_0\int_0^t {\varphi(t')\sqrt{W^{(n)}(t')}
                   \over
                 \sqrt{\nu(t-t')}}\,dt'
           + \sqrt{W(0)};$$

\vfill
\eject

 % e223

\oddheader{223}

\noindent These show that we have for all 
$n$

  $$ V^{(n)}(t) \leq \varphi(t) \qquad{\rm for\ } 0\leq t\leq \tau,$$
if $\varphi(t)$ is a continuous function satisfying

  $$\varphi(t)\geq A'\int_0^t 
        {\varphi^2(t')\over\sqrt{\nu(t-t')}}
         \,dt' + V(0).$$
We choose $\varphi(t)=(1+A)V(0)$ which gives $\tau$ the value

  $$\tau = A\nu V^{-2}(0).
  \leqno (3.8)$$

Then let

  $$v^{(n)}(t)=\max\sqrt{[u_i^{(n)}(x,t)-u_i^{(n+1)}(x,t)]
                          [u_i^{(n)}(x,t)-u_i^{(n+1)}(x,t)]}
       $$
at time $t$.

  We have
  $$v^{(1)}(t)
     < 
     A'\int_0^\tau {V^2(0)
                   \over
                 \sqrt{\nu(\tau-t')}}\,dt'
     = A V(0)$$
  $$v^{(n+1)}(t)
      <
     A\int_0^\tau {\varphi(t')v^{(n)}(t')
                   \over
                 \sqrt{\nu(\tau-t')}}\,dt'
     = A V(0)\int_0^\tau {v^{(n)}(t')
                   \over
                 \sqrt{\nu(\tau-t')}}\,dt'$$
From this we get that the functions 
$u_i^{(n)}(x,t)$
converge uniformly to continuous limits
$u_i(x,t)$
for $0\leq t\leq\tau$.

  One shows without difficulty that in the interior of the interval, each
of the derivatives of the $u_i^{(n)}(x,t)$ converges uniformly to the
 corresponding derivative of the  $u_i(x,t)$; the reasoning is too close
to that of paragraph 15 to repeat it.
 The functions $u_i(x,t)$ therefore satisfy Navier's equations for
$0<t<\tau$.

\smallskip

  We verify that the integral
$W(t)=\intpi u_i(x,t)u_i(x,t)\dee{x}$ is less than a continuous 
function of  $t$. Inequalities 
(2.5) and  (2.12) give the following, where $A_0$ is a 
constant

  $$\sqrt{W^{(0)}(t)}
             \leq
           \sqrt{W(0)}$$
  $$\sqrt{W^{(n+1)}(t)}
             \leq
      A_0\int_0^t {\varphi(t')\sqrt{W^{(n)}(t')}
                   \over
                 \sqrt{\nu(t-t')}}\,dt'
           + \sqrt{W(0)}.$$

\vfill
\eject

 % f224

\evenheader{224}

\noindent la th\'eorie des \'equations lin\'eaires nous apprend
l'existence d'une fontion positive $\theta(t)$ solution de 
l'\'equation:

  $$\theta(t)=A_0\int_0^t
                 {\varphi(t')\theta(t)\over
                  \sqrt{\nu(t-t')}}
                  \,dt'
           + \sqrt{W(0)};$$
nous avons $W^{(n)}(t)\leq \theta^2(t)$; donc $W(t)\leq \theta^2(t)$.

  Il nous reste \`a pr\'eciser comment les fonctions $u_i(x,t)$
se comportent quand $t$ tend vers z\'ero. Nous savons d\'eja qu'elles
se r\'eduisent alors aux donn\'ees
$u_i(x,0)$, en restant continues m\^eme pour $t=0$. Pour prouver
qu'elles demeurent fortement continues en $t$ quand $t$ s'annule,
il suffit d'apr\`es le lemme 1 d'\'etablir que:

  $$\limsup_{t\rightarrow 0} W(t) \leq W(0);$$
or cette in\'egalit\'e a manifestement lieu, puisque
$\theta^2(0)=W(0)$. On prouve de m\^eme que les fonctions 
${\partial u_i(x,t)\over\partial x_k}$ sont fortement continues en $t$,
m\^eme pour $t=0$.

  D\`es lors la d\'emonstration du th\'eorem d'existence \'enonc\'e
ci-dessus est achev\'ee.

  Mais la formule (3.8) nous fournit un second r\'esultat: Convenons de
dire {\it qu'une solution des \'equations de Navier, } r\'eguli\`ere dans
un intervalle $(\Theta,T)$, {\it devient irr\'eguli\`ere \`a 
l'\'epoque T} quand $T$ est fini et qu'il est impossible de d\'efinir
cette solution r\'eguli\`ere dans un intervalle $(\Theta,T')$ plus grand
que $(\Theta,T)$. La formule (3.8) r\'ev\`ele:

\medskip

  {\it Un premier caract\`ere des irr\'egularit\'es:} Si une solution
des \'equations de Navier devient irr\'eguli\`ere \`a l'\'epoque $T$,
alors $V(t)$ augmente ind\'efiniment quand $t$ tend vers $T$; et plus
pr\'ecis\'ement:

  $$V(t) > A \sqrt{\nu\over T-t}.
  \leqno (3.9)$$

\medskip

  {\bf 20. } {\it Il serait important de savoir s'il existe des solutions
des \'equations de Navier qui deviennent irr\'eguli\`eres.}
S'il ne s'en trouvait pas, la solution r\'eguli\`ere unique qui
correspond \`a un \'etat initial r\'egulier, $u_i(x,0)$, existerait
pour toutes les valeurs positives de $t$.

  Aucune solution ne pourrait devenir irr\'eguli\`ere si l'in\'egalit\'e
(3.9) \'etait incompatible avec les relations fondamentales (3.4),
(3.5) et (3.6); mais il n'en est rien, comme en le voit en choisissant:

\vfill
\eject

 % e224

\evenheader{224}

\noindent By the theory of linear equations there is a positive function
$\theta(t)$ satisfying 

  $$\theta(t)=A_0\int_0^t
                 {\varphi(t')\theta(t)\over
                  \sqrt{\nu(t-t')}}
                  \,dt'
           + \sqrt{W(0)}.$$
We have
$W^{(n)}(t)\leq \theta^2(t)$, so $W(t)\leq \theta^2(t)$.

  It rests upon us to make precise how the $u_i(x,t)$
behave when $t$ tends to zero. We already know  that
they reduce to the given 
$u_i(x,0)$, remaining continuous for $t=0$. To show that they remain
strongly continuous in 
$t$ when $t$ approaches zero, it suffices by lemma 1 to prove 

  $$\limsup_{t\rightarrow 0} W(t) \leq W(0).$$
This inequality is clear since
$\theta^2(0)=W(0)$. One shows in the same way that
the
${\partial u_i(x,t)\over\partial x_k}$ are strongly continuous
in $t$, even for
$t=0$.

  At this point the proof of the existence theorem announced above 
is complete.

  But formula (3.8) furnishes a second result: Let us say that 
{\it a solution of Navier's equations, } regular in a interval
$(\Theta,T)$, {\it  becomes irregular at time T}
when $T$ is finite and it is impossible to  extend the regular solution
to any larger interval
$(\Theta,T')$.
Formula (3.8) reveals

\medskip

  {\it A first characterization of irregularities} If a solution
of Navier's equations becomes irregular at time $T$, then 
$V(t)$ becomes arbitrarily large as $t$ tends to $T$, and  more 
precisely

  $$V(t) > A \sqrt{\nu\over T-t}.
  \leqno (3.9)$$

\medskip

  {\bf 20. } {\it It will be important to know whether there are 
solutions which become irregular.}
If these cannot be found to exist, then the regular solution corresponding
to a regular initial state 
$u_i(x,0)$ will exist for all positive values of
$t$.

  No solution can become irregular if inequality
(3.9) is incompatible with the fundamental relations
(3.4),
(3.5) and (3.6), but this is not an issue as one sees by choosing

\vfill
\eject

 % f225

\foddheader{225}

  $$V(t)=A'_0[\nu(T-t)]^{-{1\over 2}}; 
        W(t)=A''_0[\nu(T-t)]^{1\over 2};
          J(t)={\sqrt{A''_0}\over 2}[\nu(T-t)]^{-{1\over 4}}; 
  \leqno (3.10)$$
et en
v\'erifiant que pour des valeurs suffisamment fortes des contantes
$A'_0$ et $A''_0$ l'in\'egalit\'e (3.9) et la relation (3.4) sont
v\'erifi\'ees, ainsi que les deux in\'egalit\'es suivantes, qui sont 
plus strictes que (3.5) et (3.6):

  $$V(t)<A'\int_{t_0}^t
            \{
              {V^2(t')\over \sqrt{\nu(T-t')}};
             {W(t')\over [\nu(T-t')]^2}
            \}\,dt'
       +
            \{
              V(t_0);{A'''J(t_0)\over [\nu(t-t_0)]^{1\over 4}}
            \}$$

  $$J(t)< A''\int_{t_0}^t
              {J(t')V(t')\over \sqrt{\nu(T-t')}}
             \,dt'
           + J(t_0) \qquad (t_0<t<T).$$

  Les \'equations de Navier poss\`edent s\^urement une solution qui
devient irr\'eguli\`ere et pour laquelle les fonctions 
$W(t)$, $V(t)$ et $J(t)$ sont du type (3.10) si le syst\`eme:

  $$\nu\Delta U_i(x)
      -\alpha\left[U_i(x)+x_k{\partial U_i(x)\over\partial x_k}\right]
       -{1\over \rho}{\partial P(x)\over\partial x_i}
     =
       U_k(x){\partial U_i(x)\over\partial x_k};
    \leqno (3.11)$$
  $${\partial U_k(x)\over\partial x_k}=0,$$
o\`u $\alpha$ d\'esigne une constante positive, poss\`ede une solution
non nulle, les $U_i(x,t)$ \'etant born\'es et les int\'egrales
$\intpi U_i(x,t)U_i(x,t)\dee{x}$ finies; la solution des \'equations de Navier
dont il s'agit est:

  $$u_i(x,t)=
     [2\alpha(T-t)]^{-{1\over 2}} U_i[(2\alpha(T-t))^{-{1\over 2}}x]
      \qquad (t<T)
    \leqno (3.12)$$
($\lambda x$ d\'esigne le point de coordonn\'ees
$\lambda x_1$, $\lambda x_2$, $\lambda x_3$.)

  Je n'ai malheureusement pas r\'eussi \`a faire l'\'etude du syst\`eme
(3.11). Nous laisserons donc en suspens cette question de savoir si des
irr\'egularit\'es peuvent ou non se pr\'esenter.

\medskip

   {\bf 21. }{\it Cons\'equences diverses des relations fondamentales}
(3.4), (3.5) {\it et} (3.6). Soit une solution des \'equations de Navier,
r\'eguli\`ere pour 
$\Theta\leq t < T$ qui, lorsque $t$ tend vers $T$, devient irr\'eguli\`ere,
\`a moins que $T$ ne soit \'egal \`a $+\infty$. Des relations fondamentales 
(3.4) et (3.5) r\'esulte l'in\'egalit\'e:

\vfill
\eject

 % e225

\oddheader{225}

  $$V(t)=A'_0[\nu(T-t)]^{-{1\over 2}}; 
        W(t)=A''_0[\nu(T-t)]^{1\over 2};
          J(t)={\sqrt{A''_0}\over 2}[\nu(T-t)]^{-{1\over 4}} 
  \leqno (3.10)$$
and from this check that for all sufficiently large values of the constants
$A'_0$ and $A''_0$ inequality (3.9) and relation  (3.4) are satisfied,
as well as the following two inequalities which are stronger than 
(3.5) and (3.6)

  $$V(t)<A'\int_{t_0}^t
            \{
              {V^2(t')\over \sqrt{\nu(T-t')}};
             {W(t')\over [\nu(T-t')]^2}
            \}\,dt'
       +
            \{
              V(t_0);{A'''J(t_0)\over [\nu(t-t_0)]^{1\over 4}}
            \}$$

  $$J(t)< A''\int_{t_0}^t
              {J(t')V(t')\over \sqrt{\nu(T-t')}}
             \,dt'
           + J(t_0) \qquad (t_0<t<T).$$

  Navier's equations certainly have a solution which becomes irregular
and for which 
$W(t)$, $V(t)$ and $J(t)$ are of the type (3.10) if the system

  $$\nu\Delta U_i(x)
      -\alpha\left[U_i(x)+x_k{\partial U_i(x)\over\partial x_k}\right]
       -{1\over \rho}{\partial P(x)\over\partial x_i}
     =
       U_k(x){\partial U_i(x)\over\partial x_k};
    \leqno (3.11)$$
  $${\partial U_k(x)\over\partial x_k}=0,$$
where 
$\alpha$ is a positive constant, has a nonzero solution
with the $U_i(x,t)$ bounded and the integrals
$\intpi U_i(x,t)U_i(x,t)\dee{x}$ 
finite. It is 

  $$u_i(x,t)=
     [2\alpha(T-t)]^{-{1\over 2}} U_i[(2\alpha(T-t))^{-{1\over 2}}x]
      \qquad (t<T)
    \leqno (3.12)$$
($\lambda x$ is the point with coordinates 
$\lambda x_1$, $\lambda x_2$, $\lambda x_3$.)

  Unfortunately I have not made a  successful study of system
(3.11). We therefore leave  in suspense the matter of knowing
whether irregularities occur or not.

\medskip

   {\bf 21. }{\it Various consequences of the fundamental relations}
(3.4), (3.5) {\it and} (3.6). Suppose we have a solution to Navier's 
equations, regular for 
$\Theta\leq t < T$
and which becomes irregular as $t$ tends to $T$, where $T$ is not
$+\infty$. From the fundamental relations 
(3.4) and (3.5) we get the inequality

\vfill
\eject

% f226

\evenheader{226}

  $$V(t)<A'\int_{t_0}^t
            \{
              {V^2(t')\over \sqrt{\nu(t-t')}};
             {W(t_0)\over [\nu(t-t')]^2}
            \}\,dt'
       +
            \{
              V(t_0);{A'''J(t_0)\over [\nu(t-t_0)]^{1\over 4}}
            \}
  \leqno (3.13)$$
 
  $$(\Theta\leq t_0 < t < T);$$
supposons qu'une fonction $\varphi(t)$, continue pour
$0<t\leq \tau$, v\'erifie pour ces valeurs de $t$
l'in\'egalit\'e:

  $$\varphi(t)\geq A'\int_0^t
            \{
              {\varphi^2(t')\over \sqrt{\nu(t-t')}};
             {W(t_0)\over [\nu(t-t')]^2}
            \}\,dt'
       +
            \{
              V(t_0);{A'''J(t_0)\over [\nu(t-t_0)]^{1\over 4}}
            \},
  \leqno (3.14)$$
nous avons alors pour les valeurs de $t$ communes aux deux intervalles
$(t_0,T)$
et
$(t_0,t_0+\tau)$:

  $$V(t)<\varphi(t-t_0);
  \leqno (3.15)$$
le premier caract\`ere des irr\'egularit\'es permet d'en d\'eduire

  $$t_0+\tau < T.
  \leqno (3.16)$$

  Supposons en outre connue une fonction $\psi(t)$ telle que

  $$\psi(t) \geq A''\int_0^t
              {\varphi(t')\psi(t')\over \sqrt{\nu(t-t')}}
            \,dt'
        + J(t_0) \qquad (0<t\leq\tau);
  \leqno (3.17)$$
de l'in\'egalit\'e (3.6) r\'esulte alors la suivante:

  $$J(t) < \psi(t-t_0) \quad {\rm pour} \quad t_0<t\leq t_0+\tau.
  \leqno (3.18)$$

\smallskip

  {\it Le premier caract\`ere des irr\'egularit\'es } se d\'eduite 
de (3.16) en choisissant:

  $$\varphi(t) = (1+A) V(t_0) \quad {\rm et } \quad \tau=A\nu V^{-2}(t_0).$$

  Le choix $\varphi(t) = (1+A) V(t_0)$ et $\tau=+\infty$ satisfait 
l'in\'egalit\'e (3.14) quand

  $$V(t_0)) > \int_0^\infty
            \{
              {AV^2(t_0)\over \sqrt{\nu t'}};
             {AW(t_0)\over (\nu t')^2}
            \}\,dt'$$
c'est-\`a-dire quand $\nu^{-3}W(t_0)V(t_0) < A$. Donc:

\smallskip

  {\it Premier cas de r\'egularit\'e: } On est assur\'e qu'une
solution r\'eguli\`ere ne devient jamais
irr\'eguli\`ere quand la quantit\'e $\nu^{-3}W(t)V(t)$ se 
trouve \^etre inf\'erieure \`a une

\vfill
\eject

 % e226

\evenheader{226}

  $$V(t)<A'\int_{t_0}^t
            \{
              {V^2(t')\over \sqrt{\nu(t-t')}};
             {W(t_0)\over [\nu(t-t')]^2}
            \}\,dt'
       +
            \{
              V(t_0);{A'''J(t_0)\over [\nu(t-t_0)]^{1\over 4}}
            \}
  \leqno (3.13)$$
 
  $$(\Theta\leq t_0 < t < T).$$
We suppose there is a continuous function $\varphi(t)$ in 
$0<t\leq \tau$, satisfying the inequality

  $$\varphi(t)\geq A'\int_0^t
            \{
              {\varphi^2(t')\over \sqrt{\nu(t-t')}};
             {W(t_0)\over [\nu(t-t')]^2}
            \}\,dt'
       +
            \{
              V(t_0);{A'''J(t_0)\over [\nu(t-t_0)]^{1\over 4}}
            \}.
  \leqno (3.14)$$
We then have
  
  $$V(t)<\varphi(t-t_0)
  \leqno (3.15)$$
for values of $t$ common to the two intervals
$(t_0,T)$
and
$(t_0,t_0+\tau)$.
Then the first characterisation of irregularity implies

  $$t_0+\tau < T.
  \leqno (3.16)$$

  Further suppose we know a function $\psi(t)$ such that

  $$\psi(t) \geq A''\int_0^t
              {\varphi(t')\psi(t')\over \sqrt{\nu(t-t')}}
            \,dt'
        + J(t_0) \qquad (0<t\leq\tau).
  \leqno (3.17)$$
Then inequality (3.6) gives

  $$J(t) < \psi(t-t_0) \quad {\rm for} \quad t_0<t\leq t_0+\tau.
  \leqno (3.18)$$

\smallskip

  {\it The first characterisation of irregularities  } follows from
(3.16)  if we choose

  $$\varphi(t) = (1+A) V(t_0) \quad {\rm and } \quad \tau=A\nu V^{-2}(t_0).$$

  The choice $\varphi(t) = (1+A) V(t_0)$ and $\tau=+\infty$ satisfies
(3.14) if 

  $$V(t_0)) > \int_0^\infty
            \{
              {AV^2(t_0)\over \sqrt{\nu t'}};
             {AW(t_0)\over (\nu t')^2}
            \}\,dt'$$
i.e. when $\nu^{-3}W(t_0)V(t_0) < A$. So:

\smallskip

  {\it First  case of regularity: } A regular solution never becomes
irregular if the quantity
$\nu^{-3}W(t)V(t)$ is less than a 

\vfill
\eject

 % f227

\foddheader{227}

\noindent certaine constante $A$ soit \`a l'instant initial, soit \`a tout autre
instant ant\'erieurement auquel cette solution n'est pas devenue
irr\'eguli\`ere.

  On peut satisfaire (3.14) et (3.17) par un choix du type

   $$\varphi(t)=AJ(t_0)[\nu (t-t_0)]^{-{1\over 4}};\quad
        \psi(t)=(1+A)J(t_0);\quad
          \tau=A\nu^3J^{-4}(t_0).
   \leqno (3.19)$$
Cette expression fournit:

\smallskip

  {\it Un second caract\`ere des irr\'egularit\'ees:} Si une solution
des \'equations de Navier devient irr\'eguli\`ere \`a l'\'epoque $T$,
alors $J(t)$ augmente ind\'efiniment quand $t$ tend vers $T$; et 
plus pr\'ecis\'ement:

  $$J(t)> {A\nu^{3\over 4}\over(T-t)^{1\over 4}}.$$

  Les in\'egalit\'es (3.15) et (3.19) prouvent qu'une solution r\'eguli\`ere
\`a un instant $t$ reste r\'eguli\`ere jusq'\`a l'instant $t_0+\tau$
et que l'on a:

  $$V(t_0+\tau)<A\nu^{-1}J^2(t_0).$$
La relation fondamentale (3.4) donne d'autre part:

  $$W(t_0+\tau)<W(t_0).$$
Donc:

  $$\nu^{-3}W(t_0+\tau)V(t_0+\tau)<A\nu^{-4}W(t_0)J^2(t_0).$$
L'application du premier cas de r\'egularit\'e \`a l'\'epoque
$t_0+\tau$ fournit d\`es lors:

\medskip

  {\it Un second cas de r\'egularit\'e: } On est assur\'e qu'une solution
r\'eguli\`ere ne devient jamais irr\'eguli\`ere quand la quantit\'e 
$\nu^{-4}W(t)J^2(t)$ se trouve \^etre inf\'erieure \`a une certaine 
constante $A$ soit \`a l'instant initial, soit \`a tout autre 
instant ant\'erieurement auquel cette solution n'est pas devenue
irr\'eguli\`ere.

\medskip

  {\bf 22. } On \'etablit de m\^eme les r\'esultats suivants dont les 
pr\'ec\'edent peuvent d'ailleurs \^etre consid\'er\'es comme des
cas particuliers:

\ninerm 

{\it Caract\`ere des irr\'egularit\'es: } Si une solution
devient irr\'eguli\`ere \`a l'\'epoque $T$, on a:

  $$\{\intpi [u_i(x,t)u_i(x,t)]^{p\over 2}\dee{x}\}^{1\over p}
       >
        { 
           A(1-{3\over p})\nu^{{1\over 2}(1+{3\over p})}
          \over
           (T-t)^{{1\over 2}(1-{3\over p})}
        }     \qquad (p>3).$$

{\it Cas de r\'egularit\'e: } On est assur\'e qu'une solution 
r\'eguli\`ere ne devient jamais irr\'eguli\`ere quand on a \`a un instant
quelconque:

\rm
\vfill
\eject

 % e227

\oddheader{227}

\noindent certain constant $A$ either initially or at any other
instant at which the solution has not become irregular.

  One can satisfy (3.14) and (3.17) by a choice of the type

   $$\varphi(t)=AJ(t_0)[\nu (t-t_0)]^{-{1\over 4}};\quad
        \psi(t)=(1+A)J(t_0);\quad
          \tau=A\nu^3J^{-4}(t_0).
   \leqno (3.19)$$
This gives

\smallskip

  {\it A second characterisation of irregularities:}
 If a solution of Navier's equations becomes irregular at time
$T$,
then $J(t)$ grows indefinitely as $t$ tends to $T$; and
more precisely 

  $$J(t)> {A\nu^{3\over 4}\over(T-t)^{1\over 4}}.$$

  Inequalities (3.15) and (3.19) show that a solution
regular at  $t$ remains regular until $t_0+\tau$
and that one has

  $$V(t_0+\tau)<A\nu^{-1}J^2(t_0).$$
The fundamental relation (3.4) further gives

  $$W(t_0+\tau)<W(t_0).$$
Therefore

  $$\nu^{-3}W(t_0+\tau)V(t_0+\tau)<A\nu^{-4}W(t_0)J^2(t_0).$$
An application of the first case of regularity to the time
$t_0+\tau$ now gives

\medskip

  {\it A second case of regularity: } A regular solution never
becomes irregular if
$$\nu^{-4}W(t)J^2(t)$$ is less than a certain constant $A$ either
initially or at all other previous instants at which the solution
has not become irregular.

\medskip

  {\bf 22. } One similarly establishes the following results, of which
the preceeding are particular cases.

\ninerm 

{\it Characterisation of irregularities: } If a solution becomes 
irregular at time $T$, one has

  $$\{\intpi [u_i(x,t)u_i(x,t)]^{p\over 2}\dee{x}\}^{1\over p}
       >
        { 
           A(1-{3\over p})\nu^{{1\over 2}(1+{3\over p})}
          \over
           (T-t)^{{1\over 2}(1-{3\over p})}
        }     \qquad (p>3).$$

{\it Case of regularity: } A regular solution never
becomes irregular if at some time

\rm
\vfill
\eject

 % f228

\evenheader{228}

\ninerm

   $$[AW(t)]^{p-3}\intpi [u_i(x,t)u_i(x,t)]^{p\over 2}\dee{x}\}
      <
        A(1-{3\over p})^3\nu^{3(p-2)} 
     \qquad (p>3).$$

  Les cas de r\'egularit\'e que nous signalons montrent
comment une solution reste toujours r\'eguli\`ere quand son \'etat initial de
vitesse est suffisamment voisin du repos. Plus g\'en\'eralement
consid\'erons un \'etat de vitesse auquel correspnd une solution
ne devenant jamais irr\'eguli\`ere; \`a tout \'etat initial suffisamment voisin
correspond une solution qui elle aussi ne devient jamais
irr\'eguli\`ere. La d\'emonstration de ce fait utilise ceux des r\'esultats
du paragraphe 34 qui concernent l'allure d'une solution des \'equations
de Navier pour les grandes valeurs de $t$.

\rm
\bigskip

\centerline{\bf IV. Etats initiaux semi-r\'eguliers.}
\medskip

  {\bf 23. } Nous serons amen\'es, dans le courant du chapitre VI,
\`a envisager des \'etats initiaux non r\'eguliers au sens du paragraphe 17.
Commen\c cons leur \'etude en remarquant que l'in\'egalit\'e
(3.7) permet d'\'enoncer un th\'eor\`eme d'unicit\'e plus g\'en\'eral
que celui du paragraphe 18: Posons \`a cet effet une d\'efinition:

\smallskip

  {\it Nous dirons qu'une solution des \'equations de Navier est
semi-r\'eguli\`ere pour } $\Theta\leq t<T$ quand elle est r\'eguli\`ere
pour $\Theta < t<T$ et que les deux circonstances suivantes se pr\'esentent:

      \qquad\qquad\qquad\qquad L'int\'egrale $\int_\Theta^t V^2(t')\,dt'$
est finie quand $\Theta < t<T$.

  \qquad\qquad Les fonctions $u_i(x,t)$ ont de fortes limites en moyenne,
$u_i(x,\Theta)$ quand $t$ tend vers $\Theta$.

  --- Nous nommerons $<<$ \'etat initial des vitesses $>>$ ce vecteur
$u_i(x,\Theta)$, dont la quasi-divergence est nulle. ---

  Le th\'eor\`eme que fournit l'in\'egalit\'e (3.7) est le suivant:

\medskip

  {\it Th\'eor\`eme d'unicit\'e: } Deux solutions des \'equations de
Navier, semi-r\'eguli\`eres pour $\Theta\leq t<T$, sont
n\'ecessairement identiques pour toutes ces valeurs de $t$ quand
leurs \'etats de vitesse \`a l'instant $\Theta$ sont presque
partout identiques.

\smallskip

  Nous dirons qu'{\it une \'etat initial de vitesse,} $u_i(x,0)$,
est semi-r\'egulier quand it lui correspond une solution
$u_i(x,t)$ semi-r\'eguli\`ere sur un intervalle $0\leq t<\tau$.

\medskip

   {\bf 24. } Soit un vecteur $U_i(x)$ de quasi-divergence nulle,
dont les composantes sont de carr\'es sommables sur $\Pi$
et poss\`edent des quasi-d\'eriv\'ees $U_{i,j}(x)$
de carr\'es sommables sur $\Pi$. Nous allons \'etablir que le champ
de vitesses $U_i(x)$ est un \'etat initial semi-r\'egulier.

\vfill
\eject

 % e228

\evenheader{228}

\ninerm

   $$[AW(t)]^{p-3}\intpi [u_i(x,t)u_i(x,t)]^{p\over 2}\dee{x}\}
      <
        A(1-{3\over p})^3\nu^{3(p-2)} 
     \qquad (p>3).$$

  The case of regularity which we are pointing out shows how a
solution always remains regular if its initial velocity state is
sufficiently near rest. More generally, consider a velocity state
to which corresponds a solution which never becomes irregular. For
all initial states sufficiently near there corresponds a solution
which also never becomes irregular. The proof makes use of those
results of paragraph 34 which concern behavior  of
solutions to Navier's equations for large values of $t$.

\rm
\bigskip

\centerline{\bf IV. Semi-regular initial states.}
\medskip

  {\bf 23. } We will be led by the current of Chapter VI to
consider initial states which are not regular in the sense of
paragraph 17. We begin their study with the remark, that inequality
(3.7) allows a uniqueness theorem which is more general than
that of paragraph 18. To this end we make a definition.

\smallskip

  {\it We say that a solution of Navier's equations is semi-regular for }
$\Theta\leq t<T$ if it is regular for  
$\Theta < t<T$ and the two following conditions hold.

\smallskip
      \qquad\qquad\qquad\qquad The integral $\int_\Theta^t V^2(t')\,dt'$
is finite when $\Theta < t<T$.

\smallskip
  \qquad\qquad The $u_i(x,t)$ have $u_i(x,\Theta)$ as strong limit in mean
as $t$ tends to $\Theta$.

\smallskip
  --- We call ``initial velocity state''
any vector $u_i(x,\Theta)$, with quasi-divergence zero. ---

  The theorem given by inequality (3.7) is the following.

\medskip

  {\it Uniqueness theorem: } Two solutions of Navier's equations
which are semi-regular for $\Theta\leq t<T$, are necessarily
identical for all values of $t$ if their initial velocity
states at time $\Theta$ are equal almost everywhere.

\smallskip

  We say that {\it an initial velocity state } $u_i(x,0)$
is semi-regular if there corresponds a semi-regular
solution
$u_i(x,t)$ on an interval $0\leq t<\tau$.

\medskip

   {\bf 24. } Suppose a vector  $U_i(x)$ has quasi-divergence zero,
components square summable on $\Pi$, and quasi-derivatives
$U_{i,j}(x)$
square summable on $\Pi$. We are going to establish that the velocity
field $U_i(x)$ is a semi-regular initial state.

\vfill
\eject

% f229

\foddheader{229}

  Posons:

  $$W(0)=\intpi U_i(x) U_i(x)\dee{x} \quad{\rm et }\quad
    J^2(0)=\intpi U_{i,j}(x) U_{i,j}(x)\dee{x}.$$
Les fonctions $\overline{U_i(x)}$ constituent un \'etat initial r\'egulier,
comme le prouvent le lemme 6 et le paragraphe 8 (p. 209 et 206); soit
$u_i^*(x,t)$ la solution r\'eguli\`ere qui correspond \`a l'\'etat
initial $\overline{U_i(x)}$; nous avons en vertu de l'in\ 'egalit\'e
(1.21) et de la relation de dissipation de l'\'energie (3.4):

  $$W^*(t) < W(0).
  \leqno (4.1)$$
Le lemme 4 nous apprend que
${\partial\overline{U_i(x)}\over\partial x_j} = U_{i,j}(x)$; 
nous avons donc
d'apr\`es (1.21):

  $$J^*(0) < J(0);$$
les relations 
(3.15), (3.18), et (3.19) permettent d'en d\'eduire que sur un m\^eme
intervalle $(0,\tau)$ les diverses solutions $u_i^*(x,t)$ sont r\'eguli\`eres
et v\'erifient les in\'egalit\'es:

  $$V^*(t) < AJ(0)(\nu t)^{-{1\over 4}}; \qquad J^*(t) < (1+A)J(0);
  \leqno (4.2)$$
nous avons d'ailleurs:

  $$\tau = A\nu^3J^{-4}(0).
  \leqno (4.3)$$

  Les in\'egalit\'es (4.1) et (4.2) nous autorisent \`a appliquer
le lemme 9 (p. 220): dans la formule de d\'efinition (1.18) de 
$\overline{U(x)}$ figure une longueur $\epsilon$; il est 
possible de la faire tendre vers z\'ero en sorte que pour 
$0<t<\tau$ les fonctions $u_i^*(x,t)$ et chacune de leurs d\'eriv\'ees
convergnet respectivement vers certaines fonctions $u_i(x,t)$ et vers
leurs d\'eriv\'ees. Ces fonctions $u_i(x,t)$ constituent une
solution des \'equation de Navier r\'eguli\`ere pour $0<t<\tau$; d'apr\`es
(4.1) et (4.2) cette solution satisfait les trois in\'egalit\'es:

  $$W(t)\leq W(0);\qquad V(t)\leq AJ(0)(\nu t)^{-{1\over 4}}; \qquad 
     J(t)\leq (1+A)J(0).
  \leqno (4.4)$$
L'int\'egrale $\int_0^t V^2(t')\,dt'$ est donc finie pour 
$0<t<\tau$. Il nous reste \`a pr\'eciser comment les fonctions
$u_i(x,t)$ se comportent quand $t$ tend vers z\'ero.

\vfill
\eject

 % e229

\oddheader{229}

  Let

  $$W(0)=\intpi U_i(x) U_i(x)\dee{x} \quad{\rm et }\quad
    J^2(0)=\intpi U_{i,j}(x) U_{i,j}(x)\dee{x}.$$
The functions
$\overline{U_i(x)}$ constitute a regular initial state, as shown
by lemma 6 and paragraph 8 (p. 209 et 206). Let
$u_i^*(x,t)$ be the regular solution which 
corresponds to the initial state
$\overline{U_i(x)}$. We have, in virtue of inequality
(1.21) and the energy dissipation relation (3.4) that

  $$W^*(t) < W(0).
  \leqno (4.1)$$
Lemma 4 shows us that 
${\partial\overline{ U_i(x)}\over\partial x_j} = U_{i,j}(x)$. Thus we
have from (1.21)

  $$J^*(0) < J(0).$$
Relations 
(3.15), (3.18), and (3.19) allow us to deduce from this
that in some interval
 $(0,\tau)$ the various solutions $u_i^*(x,t)$ are regular and satisfy
inequalities

  $$V^*(t) < AJ(0)(\nu t)^{-{1\over 4}}; \qquad J^*(t) < (1+A)J(0).
  \leqno (4.2)$$
We have further

  $$\tau = A\nu^3J^{-4}(0).
  \leqno (4.3)$$

  Inequalities (4.1) and (4.2) let us apply lemma  9 (p. 220).
There is a length $\epsilon$ in the definition
(1.18) of 
$\overline{U(x)}$. It is possible to make this tend to zero in 
such a way that for 
$0<t<\tau$ the functions $u_i^*(x,t)$ and all their derivatives
converge respectively to certain functions
 $u_i(x,t)$ and to their derivatives. These $u_i(x,t)$ are a regular 
solution
to Navier's equations for
 $0<t<\tau$. By 
(4.1) and (4.2) this solution satisfies the three inequalities

  $$W(t)\leq W(0);\qquad V(t)\leq AJ(0)(\nu t)^{-{1\over 4}}; \qquad 
     J(t)\leq (1+A)J(0).
  \leqno (4.4)$$
The integral $\int_0^t V^2(t')\,dt'$ is therefore finite for 
$0<t<\tau$. Now we must specify how the 
$u_i(x,t)$ behave as $t$ tends to zero.

\vfill
\eject

 % f230

\evenheader{230}

   Soit $a_i(x)$ un vecteur quelconque, de divergence nulle,
dont les composantes, ainsi que toutes leurs d\'eriv\'ees, sont
de carr\'es sommables sur $\Pi$. Des \'equations de Navier
r\'esulte la ralation:

  $$\intpi u_i^*(x,t)a_i(x)\dee{x} = 
      \intpi \overline{U_i(x)}a_i(x)\dee{x} + $$

  $$\qquad\qquad 
      \nu\int_0^t\,dt'\intpi u_i^*(x,t')\Delta a_i(x)\dee{x}+
        \int_0^t\,dt'\intpi 
            u_k(x,t')u_i(x,t')
         {\partial a_i(x)\over
           \partial x_k}\dee{x};$$
d'o\`u, en passant \`a la limite:

  $$\intpi u_i(x,t)a_i(x)\dee{x} = 
      \intpi U_i(x)a_i(x)\dee{x} + $$

  $$\qquad\qquad 
      \nu\int_0^t\,dt'\intpi u_i(x,t')\Delta a_i(x)\dee{x}+
        \int_0^t\,dt'\intpi 
            u_k(x,t')u_i(x,t')
         {\partial a_i(x)\over
           \partial x_k}\dee{x}.$$
Cette derni\`ere relation prouve que

  $$\intpi u_i(x,t)a_i(x)\dee{x} \qquad {\rm tend \  vers }\qquad 
       \intpi U_i(x)a_i(x)\dee{x}$$
quand $t$ tend vers z\'ero. Dans ces conditions $u_i(x,t)$ a une
faible limite en moyenne unique, qui est $U_i(x)$ 
(cf, Corollaire du lemme 7, p. 209). Mais l'in\'egalit\'e
$W(t)\leq W(0)$ nous permet d'utiliser le crit\`ere de forte
convergence \'enonc\'e p. 200; et nous constatons ainsi que
les fonctions $u_i(x,t)$ convergent fortement en moyenne vers les
fonctions $U_i(x)$ quand $t$ tend vers z\'ero.

  $u_i(x,t)$ est donc une solution semi-r\'eguli\`ere$^1$ 
pour
$0\leq t < \tau$ et elle correspond \`a l'\'etat initial 
$U_i(x)$.

\medskip

   {\bf 25. } On peut par des raisonnements analogues
traiter les deux autres cas que signale le th\'eor\`eme ci-dessous:

\smallskip

  {\it Th\'eor\`eme d'existence: } Soit un vecteur $U_i(x)$, de
quasi-divergence nulle, dont

\footrule

  $^1$ On peut m\^eme affirmer plus: les fonctions
${\partial u_i(x,t)\over\partial x_j}$ convergent fortement en
moyenne vers les fonctions $U_{i,j}(x)$ quand $t$ tend vers z\'ero.

\rm
\vfill
\eject

 % e230

\evenheader{230}

  Let $a_i(x)$ be any vector of divergence zero, for which the components
as well as all their derivatives are square summable on 
$\Pi$. From Navier's equations we get

  $$\intpi u_i^*(x,t)a_i(x)\dee{x} = 
      \intpi \overline{U_i(x)}a_i(x)\dee{x} + $$

  $$\qquad\qquad 
      \nu\int_0^t\,dt'\intpi u_i^*(x,t')\Delta a_i(x)\dee{x}+
        \int_0^t\,dt'\intpi 
            u_k(x,t')u_i(x,t')
         {\partial a_i(x)\over
           \partial x_k}\dee{x}.$$
Then passing to the limit

  $$\intpi u_i(x,t)a_i(x)\dee{x} = 
      \intpi U_i(x)a_i(x)\dee{x} + $$

  $$\qquad\qquad 
      \nu\int_0^t\,dt'\intpi u_i(x,t')\Delta a_i(x)\dee{x}+
        \int_0^t\,dt'\intpi 
            u_k(x,t')u_i(x,t')
         {\partial a_i(x)\over
           \partial x_k}\dee{x}.$$
This last relation shows that

  $$\intpi u_i(x,t)a_i(x)\dee{x} \qquad {\rm tends \  to }\qquad 
       \intpi U_i(x)a_i(x)\dee{x}$$
when $t$ tends to zero. In these conditions
$u_i(x,t)$ has a unique weak limit in mean, which is
$U_i(x)$ 
(cf. Corollary to lemma 7, p. 209). But the inequality
$W(t)\leq W(0)$ allows us to use the criteria for strong
convergence announced on
p. 200, and we also note that the
$u_i(x,t)$ converge  strongly in mean to the
$U_i(x)$ as $t$ tends to zero.

  $u_i(x,t)$ is therefore a semi-regular solution$^1$ 
for
$0\leq t < \tau$ and it corresponds to the initial state
$U_i(x)$.

\medskip

   {\bf 25. } By analogous reasoning one can treat the two other
cases pointed out in the theorem below.

\smallskip

  {\it Existence theorem: } Let the vector $U_i(x)$ have 
quasi-divergence zero, with

\footrule

  $^1$ One can similarly check that the functions
${\partial u_i(x,t)\over\partial x_j}$ converge strongly in mean
to the 
$U_{i,j}(x)$ as $t$ tends to zero.

\rm
\vfill
\eject

 % f231

\foddheader{231}

\noindent les composantes sont de carr\'es sommables sur $\Pi$; on 
peut affirmer que l'\'etat initial de vitesses qu'il d\'efinit est
semi-r\'egulier:

  a) quand les fonctions $U_i(x)$ poss\`edent des quasi-d\'eriv\'ees
de carr\'es sojmmables sur $\Pi$;

  b) quand les fonctions $U_i(x)$ sont born\'ees;

  c) ou enfin quand l'int\'egrale 
$\intpi [U_i(x)U_i(x)]^{p\over 2}\dee{x}$ est finie pour une valeur
de $p$ sup\'erieure \`a 3.

  {\it N. B. } Ce th\'eor\`eme d'existence du paragraphe 19
n'\'epuisent \'evidemment pas l'\'etude de l'allure que pr\'esente
au voisinage de l'instant initial la solution qui correspond \`a
un \'etat initial donn\'e.

\bigskip

\centerline{\bf V. Solutions turbulentes.}

  {\bf 26. } Soit un \'etat initial r\'egulier $u_i(x,0)$. Nous
n'avons pas r\'eussi \`a prouver que la solution r\'eguli\`ere des
\'equations de Navier qui lui correspond est d\'efinie pour toutes
les valeurs de $t$ 
post\'erieures \`a l'instant initial $t=0$. Mais consid\'erons
le syst\`eme:

  $$\nu\Delta u_i(x,t)-
    {\partial u_i(x,t)\over \partial t}
    -{1\over \rho}{\partial p(x,t)\over \partial x_i} 
    =
    \overline{u_k(x,t)}{\partial u_i(x,t)\over \partial x_k}; \qquad
    {\partial u_j(x,t)\over \partial x_j} = 0.
  \leqno (5.1)$$

  {\it C'est un syst\`eme qui est tr\`es voisin des \'equations de 
Navier } quand la longeur$^1$ $\epsilon$ est tr\`es courte. Tout ce
que nous avons dit au cours du chapitre III sur les \'equations de 
Navier lui est applicable sans modification, hormis les 
consid\'erations non concluantes du paragraphe 20. Par l\`a se
trouve \'etablie toute une cat\'egorie de propri\'et\'es du syst\`eme
(5.1), dans lesquelles ne figure pas la longueur $\epsilon$. D'autre
part l'in\'egalit\'e de Schwarz (1.1) nous donne:

  $$\overline{u_k(x,t)} < 
        A_0 \epsilon^{-{3\over 2}}\sqrt{W(t)},$$
$A_0$ \'etant une constante num\'erique. Cette nouvelle in\'egalit\'e
 et la relation de dissipation de l'\'energie (3.4)
autorisent \`a \'ecrire \`a c\^ot\'e de l'in\'egalit\'e (3.5) la
suivante: si une solution du syst\`eme (5.1) est r\'eguli\`ere
pour 
$0\leq t<T$, alors:

\footrule

  $^1$ Rappelons que cette longueur a \'et\'e introduite au 
\S 8 (p, 206), quand nous avons d\'efini le symbole $\overline{U(x)}$.

\rm
\vfill
\eject

 % e231

\oddheader{231}

\noindent the components are square summable on $\Pi$. One
can verify that the initial velocity state which it defines is
semi-regular

  a) if the functions $U_i(x)$ have square summable quasi-derivative
on $\Pi$;

  b) if the functions $U_i(x)$ are bounded;

  c) or finally if the integral 
$\intpi [U_i(x)U_i(x)]^{p\over 2}\dee{x}$ is finite for some value
of $p$ larger than 3.

  {\it N. B. } This theorem and the existence theorem of paragraph 
19 evidently do not allow a study of the 
attraction which is presented in the neighborhood of 
the initial instant by the solution 
which corresponds to
a given initial state.

\bigskip

\centerline{\bf V. Turbulent solutions.}

  {\bf 26. } Let $u_i(x,0)$ be a regular initial state. We have not
succeeded in proving that the corresponding regular solution
to Navier's equations is defined for all values of $t$ after the
initial instant $t=0$. But consider the system

  $$\nu\Delta u_i(x,t)-
    {\partial u_i(x,t)\over \partial t}
    -{1\over \rho}{\partial p(x,t)\over \partial x_i} 
    =
    \overline{u_k(x,t)}{\partial u_i(x,t)\over \partial x_k}; \qquad
    {\partial u_j(x,t)\over \partial x_j} = 0.
  \leqno (5.1)$$

  {\it This system is very near Navier's equations } when the 
length$^1$ $\epsilon$ is very short.
All we have said in Chapter III on Navier's equations is
applicable without modification,
other than the inconclusive considerations of paragraph 20.
Thus we know many properties of system (5.1) which are independent
of $\epsilon$. Further, the Schwarz inequality (1.1) gives us

  $$\overline{u_k(x,t)} < 
        A_0 \epsilon^{-{3\over 2}}\sqrt{W(t)},$$
$A_0$ being a numerical constant. This new inequality
and the energy dissipation relation (3.4)
allows us to write the following beside inequality (3.5) if a 
solution to system (5.1) is regular for 
$0\leq t<T$, then

\footrule

  $^1$ Recall this length was introduced in 
\S 8 (p, 206), when we defined the symbol $\overline{U(x)}$.

\rm
\vfill
\eject

% f232

\evenheader{232}

  $$V(t)<A'A_0\epsilon^{-{3\over 2}}\sqrt{W(0)}
        \int_0^t
         {V(t')\,dt'\over
          \sqrt{\nu(t-t')}}
       + V(0)
           \qquad (0<t<T).$$
De l\`a r\'esulte que sur tout intervalle de r\'egularit\'e
$(0,T)$ $V(t)$ 
reste inf\'erieur \`a la fonction
$\varphi(t)$,
continue pour 
$0\leq t$,
qui satisfait l'\'equation int\'egrale lin\'eaire du type de Volterra:

  $$\varphi(t)=A'A_0\epsilon^{-{3\over 2}}\sqrt{W(0)}
        \int_0^t
         {\varphi(t')\,dt'\over
          \sqrt{\nu(t-t')}}
       + V(0);$$
$V(t)$ reste donc born\'e quand, $T$ \'etant fini, $t$ tend vers $T$;
ceci contredit le premier caract\`ere des irr\'egularit\'es (p. 224);
en d'autres termes {\it l'unique solution des \'equations (5.1)
qui correspond \`a un \'etat initial r\'egulier donn\'e est d\'efinie
pour toutes les valeurs du temps post\'erieures \`a l'instant initial.}

\medskip

  {\bf 27. } \'Etant donn\'e un mouvement qui satisfait les \'equations
(5.1), nous aurons besoin de {\it r\'esultats concernant la
r\'epartition de son \'energie cin\'etique:}
${1\over 2}u_i(x,t)u_i(x,t)$. Ces r\'esultats devront \^etre
ind\'ependants$^1$ de $\epsilon$.

  Soit deux longueurs constantes $R_1$ et $R_2$ $(R_1<R_2)$;
introduisons la fonction $f(x)$ suivante:

  $$f(x) = 0 \quad {\rm pour} \quad r_0\leq R_1;$$

  $$f(x) = {r_0-R_1\over
            R_2-R_1} \quad {\rm pour} \quad R_1\leq r_0\leq R_2;
     \qquad (r_0^2 = x_i x_i)$$

  $$f(x) = 1 \quad {\rm pour} \quad R_2\leq r_0.$$
Un calcul analogue \`a celui qui fournit la relation de dissipation
de l'\'energie (2.21) nous donne:

  $$\nu\int_0^t \,dt'
       \intpi
         f(x){\partial u_i(x,t')\over
               \partial x_k}
              {\partial u_i(x,t')\over
               \partial x_k}
        \dee{x}
      +
     {1\over 2}\intpi
         f(x)u_i(x,t)u_i(x,t)
        \dee{x}
      = $$

  $$ = {1\over 2}\intpi
         f(x)u_i(x,0)u_i(x,0)
        \dee{x}
      -
       \nu\int_0^t \,dt'
       \intpi
         {\partial f(x)\over
               \partial x_k}
              u_i(x,t')
              {\partial u_i(x,t)\over
               \partial x_k}
        \dee{x} + $$

\footrule

  $^1$ Ils vaudront  \'egalement pour les solutions r\'eguli\`eres
des \'equations de Navier.

\rm
\vfill
\eject

 % e232

\evenheader{232}

  $$V(t)<A'A_0\epsilon^{-{3\over 2}}\sqrt{W(0)}
        \int_0^t
         {V(t')\,dt'\over
          \sqrt{\nu(t-t')}}
       + V(0)
           \qquad (0<t<T).$$
From this we get that on all intervals of regularity
$(0,T)$, $V(t)$ 
remains less than the continuous function
$\varphi(t)$ on $0\leq t$,
which satisfies the Volterra linear integral equation

  $$\varphi(t)=A'A_0\epsilon^{-{3\over 2}}\sqrt{W(0)}
        \int_0^t
         {\varphi(t')\,dt'\over
          \sqrt{\nu(t-t')}}
       + V(0).$$
$V(t)$ therefore remains bounded when, $T$ being finite, $t$ tends to $T$.
That contradicts the first caracterization of irregularity (p. 224). 
In other words, {\it the unique solution to equations (5.1)
corresponding to a given regular initial state is defined for all time
after the initial instant.}

\medskip

  {\bf 27. } Given a motion which satisfies
equations (5.1), we will need results on its {\it repartition of
kinetic energy:} 
${1\over 2}u_i(x,t)u_i(x,t)$. These must be independent$^1$ of 
$\epsilon$.

  Consider two constant lengths $R_1$ and $R_2$ $(R_1<R_2)$ and
introduce the following function $f(x)$

  $$f(x) = 0 \quad {\rm for} \quad r_0\leq R_1;$$

  $$f(x) = {r_0-R_1\over
            R_2-R_1} \quad {\rm for } \quad R_1\leq r_0\leq R_2;
     \qquad (r_0^2 = x_i x_i)$$

  $$f(x) = 1 \quad {\rm for } \quad R_2\leq r_0.$$
A calculation analogous to that giving the energy dissipationn relation
(2.21) here gives

  $$\nu\int_0^t \,dt'
       \intpi
         f(x){\partial u_i(x,t')\over
               \partial x_k}
              {\partial u_i(x,t')\over
               \partial x_k}
        \dee{x}
      +
     {1\over 2}\intpi
         f(x)u_i(x,t)u_i(x,t)
        \dee{x}
      = $$

  $$ = {1\over 2}\intpi
         f(x)u_i(x,0)u_i(x,0)
        \dee{x}
      -
       \nu\int_0^t \,dt'
       \intpi
         {\partial f(x)\over
               \partial x_k}
              u_i(x,t')
              {\partial u_i(x,t)\over
               \partial x_k}
        \dee{x} + $$

\footrule

  $^1$ They will apply equally to regular solutions of Navier's equations.

\rm
\vfill
\eject

 % f233

\foddheader{233}

  $$+{1\over\rho}\int_0^t \,dt'
       \intpi
         {\partial f(x)\over
               \partial x_i}
              p(x,t')
              u_i(x,t')
        \dee{x} + $$

  $$+{1\over 2}\int_0^t \,dt'
       \intpi
         {\partial f(x)\over
               \partial x_k}
              \overline{u_k(x,t')}
              u_i(x,t')
              u_i(x,t')
        \dee{x}.$$
Nous en d\'eduisons l'in\'egalit\'e:

  $$ {1\over 2}\tripl{r_0>R_2}
         u_i(x,t)u_i(x,t)
        \dee{x}
         <
         {1\over 2}\tripl{r_0>R_1}
         u_i(x,0)u_i(x,0)
        \dee{x} + $$ 
  $$         +{\nu\sqrt{W(0)}\over R_2-R_1}\int_0^t J(t')\,dt'
            +{1\over\rho}{\sqrt{W(0)}\over R_2-R_1}\int_0^t \,dt'
                 \sqrt{\intpi p^2(x,t')\dee{x}}
                + 
   \leqno (5.2)$$
  $$     +{\sqrt{W(0)}\over R_2-R_1}\int_0^t \,dt'
                 \sqrt{\intpi\left[
                     {1\over 2}u_i(x,t')u_i(x,t')
                            \right]^2
                               \dee{x}}.$$
Majorons les trois derniers termes: d'apr\`es l'in\'egalit\'e
de Schwarz

  $$\int_0^t J(t')\,dt'  < \sqrt{\int_0^t J^2(t')\,dt'}\sqrt{t}
                         < \sqrt{W(0)\over 2\nu}\sqrt{t}.
  \leqno (5.3)$$
D'autre part (cf. (3.3)):

  $${1\over\rho}p(x,t')=
     {1\over 4\pi}\intpi
                   {\partial ({1\over r})\over\partial x_j}
                   {\partial u_i(y,t')\over \partial y_k}
                   \overline{u_k(y,t')}
                  \dee{y},
  \leqno (5.4)$$
d'o\`u

  $${1\over\rho^2}\intpi p^2(x,t')\dee{x} = 
     {1\over 4\pi}\intpi\intpi
                      \overline{u_k(x,t')}
                      {\partial u_i(x,t')\over \partial x_k}
                      {1\over r}
                      \overline{u_j(y,t)}
                      {\partial u_i(y,t')\over \partial y_j}
                   \dee{x}\dee{y};$$
la relation (1.14) et l'in\'egalit\'e de Schwarz (1.1) nous apprennent
que

  $$\sum_i\left[\intpi
               {1\over r}\overline{u_j(y,t)}
               {\partial u_i(y,t')\over \partial y_j}\partial y
          \right]^2
      < 4J^4(t'); $$
en outre:

\vfill
\eject

 % e233

\oddheader{233}

  $$+{1\over\rho}\int_0^t \,dt'
       \intpi
         {\partial f(x)\over
               \partial x_i}
              p(x,t')
              u_i(x,t')
        \dee{x} + $$

  $$+{1\over 2}\int_0^t \,dt'
       \intpi
         {\partial f(x)\over
               \partial x_k}
              \overline{u_k(x,t')}
              u_i(x,t')
              u_i(x,t')
        \dee{x}.$$
From this we get the inequality

  $$ {1\over 2}\tripl{r_0>R_2}
         u_i(x,t)u_i(x,t)
        \dee{x}
         <
         {1\over 2}\tripl{r_0>R_1}
         u_i(x,0)u_i(x,0)
        \dee{x} + $$ 
  $$         +{\nu\sqrt{W(0)}\over R_2-R_1}\int_0^t J(t')\,dt'
            +{1\over\rho}{\sqrt{W(0)}\over R_2-R_1}\int_0^t \,dt'
                 \sqrt{\intpi p^2(x,t')\dee{x}}
                + 
   \leqno (5.2)$$
  $$     +{\sqrt{W(0)}\over R_2-R_1}\int_0^t \,dt'
                 \sqrt{\intpi\left[
                     {1\over 2}u_i(x,t')u_i(x,t')
                            \right]^2
                               \dee{x}}.$$
We majorize the last three terms. By the Schwarz inequality

  $$\int_0^t J(t')\,dt'  < \sqrt{\int_0^t J^2(t')\,dt'}\sqrt{t}
                         < \sqrt{W(0)\over 2\nu}\sqrt{t}.
  \leqno (5.3)$$
Further (cf. (3.3)):

  $${1\over\rho}p(x,t')=
     {1\over 4\pi}\intpi
                   {\partial ({1\over r})\over\partial x_j}
                   {\partial u_i(y,t')\over \partial y_k}
                   \overline{u_k(y,t')}
                  \dee{y},
  \leqno (5.4)$$
from which

  $${1\over\rho^2}\intpi p^2(x,t')\dee{x} = 
     {1\over 4\pi}\intpi\intpi
                      \overline{u_k(x,t')}
                      {\partial u_i(x,t')\over \partial x_k}
                      {1\over r}
                      \overline{u_j(y,t)}
                      {\partial u_i(y,t')\over \partial y_j}
                   \dee{x}\dee{y}.$$
Relation (1.14) and the Schwarz inequality (1.1) give 

  $$\sum_i\left[\intpi
               {1\over r}\overline{u_j(y,t)}
               {\partial u_i(y,t')\over \partial y_j}\partial y
          \right]^2
      < 4J^4(t'). $$
Further

\vfill
\eject

 % f234

\evenheader{234}

  $$\sum_i\left[\intpi
           \overline{u_k(x,t')}
           {\partial u_i(x,t')\over \partial x_k}
                \dee{x}
          \right]^2 < W(t')J^2(t');$$
donc:

  $${1\over\rho^2}\intpi p^2(x,t')\dee{x} <
      {1\over 2\pi}\sqrt{W(t')}J^3(t');$$
par suite:$^1$

  $${1\over\rho}\int_0^t \,dt'
                   \sqrt{\intpi p^2(x,t')\dee{x}} <
      {[W(0)]^{1\over 4}\over\sqrt{2\pi}}
      \int_0^t J^{3\over 2}(t')\,dt'      <
      {W(0)\over\sqrt{2\pi}(2\nu)^{3\over 4}}t^{1\over 4}.
  \leqno (5.5)$$
De (1.13) r\'esulte:

  $${1\over 2}u_i(x,t')u_i(x,t')= -{1\over 4\pi}\intpi
                   {\partial ({1\over r})\over\partial x_k}
                   u_i(y,t')
                   {\partial u_i(y,t')\over \partial y_k}
                  \dee{y};$$
cette formule analogue \`a (5.4) conduit par des calculs analogues
aux pr\'ec\'edents \`a l'in\'egalit\'e:

  $$\int_0^t \,dt'\sqrt{\intpi\left[
                     {1\over 2}u_i(x,t')u_i(x,t')
                            \right]^2
                               \dee{x}} < 
        {W(0)\over\sqrt{2\pi}(2\nu)^{3\over 4}}t^{1\over 4}.
  \leqno (5.6)$$
Tenons compte dans (5.2) des majorantes (5.3), (5.5) et (5.6); nous
obtenons:

  $${1\over 2}\tripl{r_0>R_2} u_i(x,t)u_i(x,t)\dee{x} <
    {1\over 2}\tripl{r_0>R_1} u_i(x,0)u_i(x,0)\dee{x} + $$

  $$\qquad\qquad\qquad\qquad\qquad\qquad 
     +{W(0)\sqrt{\nu t}\over\sqrt{2}(R_2-R_1)}
     +{W^{3\over 2}(0)t^{1\over 4}\over
        2^{1\over 4}\pi^{1\over 2}\nu^{3\over 4}(R_2-R_1)}.
  \leqno (5.7)$$

\footrule

  $^1$ Nous utilisons l'in\'egalit\'e:

  $$\int_0^t J^{3\over 2}(t')\,dt' <
    \left[\int_0^t J^2(t')\,dt'\right]^{3\over 4}t^{1\over 4}$$
qui est un cas particulier de $<<$l'in\'egalit\'e de H\"older$>>$:

  $$\left|\int_0^t \varphi(t')\psi(t')\,dt'\right|
     <
      \left[\int_0^t \varphi^p(t')\,dt'\right]^{1\over p}
      \left[\int_0^t \psi^q(t')\,dt'\right]^{1\over q}
       \qquad
       ({1\over p}+{1\over q}=1; 1<p, 1<q).$$
\rm
\vfill
\eject

 % e234

\evenheader{234}

  $$\sum_i\left[\intpi
           \overline{u_k(x,t')}
           {\partial u_i(x,t')\over \partial x_k}
                \dee{x}
          \right]^2 < W(t')J^2(t');$$
therefore

  $${1\over\rho^2}\intpi p^2(x,t')\dee{x} <
      {1\over 2\pi}\sqrt{W(t')}J^3(t');$$
and it follows$^1$

  $${1\over\rho}\int_0^t \,dt'
                   \sqrt{\intpi p^2(x,t')\dee{x}} <
      {[W(0)]^{1\over 4}\over\sqrt{2\pi}}
      \int_0^t J^{3\over 2}(t')\,dt'      <
      {W(0)\over\sqrt{2\pi}(2\nu)^{3\over 4}}t^{1\over 4}.
  \leqno (5.5)$$
From (1.13) we get

  $${1\over 2}u_i(x,t')u_i(x,t')= -{1\over 4\pi}\intpi
                   {\partial ({1\over r})\over\partial x_k}
                   u_i(y,t')
                   {\partial u_i(y,t')\over \partial y_k}
                  \dee{y}.$$
This formula is analogous to (5.4). By calculations like the preceeding
it leads to the inequality

  $$\int_0^t \,dt'\sqrt{\intpi\left[
                     {1\over 2}u_i(x,t')u_i(x,t')
                            \right]^2
                               \dee{x}} < 
        {W(0)\over\sqrt{2\pi}(2\nu)^{3\over 4}}t^{1\over 4}.
  \leqno (5.6)$$
Using the majorants (5.3), (5.5), and (5.6) in (5.2) we obtain

  $${1\over 2}\tripl{r_0>R_2} u_i(x,t)u_i(x,t)\dee{x} <
    {1\over 2}\tripl{r_0>R_1} u_i(x,0)u_i(x,0)\dee{x} + $$

  $$\qquad\qquad\qquad\qquad\qquad\qquad 
     +{W(0)\sqrt{\nu t}\over\sqrt{2}(R_2-R_1)}
     +{W^{3\over 2}(0)t^{1\over 4}\over
        2^{1\over 4}\pi^{1\over 2}\nu^{3\over 4}(R_2-R_1)}.
  \leqno (5.7)$$

\footrule

  $^1$ We use the inequality

  $$\int_0^t J^{3\over 2}(t')\,dt' <
    \left[\int_0^t J^2(t')\,dt'\right]^{3\over 4}t^{1\over 4}$$
which is a particular case of ``H\"older's inequality''

  $$\left|\int_0^t \varphi(t')\psi(t')\,dt'\right|
     <
      \left[\int_0^t \varphi^p(t')\,dt'\right]^{1\over p}
      \left[\int_0^t \psi^q(t')\,dt'\right]^{1\over q}
       \qquad
       ({1\over p}+{1\over q}=1; 1<p, 1<q).$$
\rm
\vfill
\eject

% f235

\foddheader{235}

\noindent Cette in\'egalit\'e renseigne sur la fa\c con dont
l'\'energie cin\'etique reste localis\'ee \`a distance finie.

\medskip

  {\bf 28. } Donnons nous \`a l'instant initial $t=0$ un \'etat
initial constitu\'e par un vecteur quelconque, $U_i(x)$, dont
les composantes sont de carr\'es sommables sur $\Pi$ 
et dont la quasi-divergence est nulle. Le vecteur
$\overline{U_i(x)}$
constitue un \'etat initial r\'egulier (cf. lemme 6 et paragraphe
8); nommons $u_i^*(x,t)$ la solution r\'eguli\`ere des \'equations
(5.1) qui lui correspond; elle est d\'efinie pour toutes les 
valeurs de $t$. {\it Le but de ce chapitre est d'\'etudier
les limites que peut avoir cette solution r\'eguli\`ere $u_i^*(x,t)$
du syst\'`eme (5.1) quand $\epsilon$ tend vers z\'ero. }

\smallskip

  Les propri\'et\'es des fonctions $u_i^*(x,t)$ dont nous ferons
usage sont les trois suivantes:

\smallskip

  1$^o$) Soit $a_i^*(x,t)$ un vecteur quelconque de divergence nulle,
dont toutes les composantes et toutes leurs d\'eriv\'ees sont 
uniform\'ement et fortement continues en $t$; nous avons d'apr\`es
(5.1):

  $$\intpi u_i^*(x,t)a_i(x,t)\dee{x} =
       \intpi \overline{U_i(x)}a_i(x,0)\dee{x} + $$

  $$ + \int_0^t\,dt'\intpi 
                     u_i^*(x,t')
                     \left[\nu\Delta a_i(x,t')
                          +{\partial a_i(x,t')\over \partial t'}
                     \right]
                     \dee{x} +
  \leqno (5.8)$$

  $$ + \int_0^t\,dt'\intpi 
                     \overline{u_i^*(x,t)}
                     u_i^*(x,t')
                     {\partial a_i(x,t')\over \partial x_k}
                    \dee{x}.$$

  2$^o$) La relation de dissipation de l'\'energie et l'in\'egalit\'e
(1.21) nous donnent:

  $$\nu\int_{t_0}^t J^{*2}(t')\,dt' + {1\over 2}W^*(t) = 
    {1\over 2}W^*(t_0) < {1\over 2}W(0).
  \leqno (5.9)$$

 $$\left({\rm Par\ definition\ }W(0)=\intpi U_i(x)U_i(x)\dee{x}\right).
  \leqno (5.10)$$

  3$^o$) L'in\'egalit\'e (5.7) et l'in\'egalit\'e $W^*(0)<W(0)$
justifient la proposition suivante:

\vfill
\eject

 % e235

\oddheader{235}

\noindent This inequality shows how kinetic energy remains
localized at finite distance.

\medskip

  {\bf 28. } Suppose we have given at $t=0$ an arbitrary initial vector
$U_i(x)$, with components square summable on 
$\Pi$ 
and quasi-divergence zero. The vector 
$\overline{U_i(x)}$
is a regular initial state (cf. lemma 6 and paragraph
8). Write $u_i^*(x,t)$ for the corresponding regular solution to equations
(5.1). It is defined for all  $t$. {\it The object of this chapter
is to study the limits which this regular solution may have
as 
$\epsilon$ tends to zero. }

\smallskip

  We will use the  following three 
properties of the $u_i^*(x,t)$.

\smallskip

  1$^o$) Let $a_i^*(x,t)$ be an arbitrary vector of divergence zero,
of which all components and all their derivatives are uniformly and
strongly continuous in $t$. By 
(5.1):

  $$\intpi u_i^*(x,t)a_i(x,t)\dee{x} =
       \intpi \overline{U_i(x)}a_i(x,0)\dee{x} + $$

  $$ + \int_0^t\,dt'\intpi 
                     u_i^*(x,t')
                     \left[\nu\Delta a_i(x,t')
                          +{\partial a_i(x,t')\over \partial t'}
                     \right]
                     \dee{x} +
  \leqno (5.8)$$

  $$ + \int_0^t\,dt'\intpi 
                     \overline{u_i^*(x,t)}
                     u_i^*(x,t')
                     {\partial a_i(x,t')\over \partial x_k}
                    \dee{x}.$$

  2$^o$) The energy dissipation relation and
(1.21) give 

  $$\nu\int_{t_0}^t J^{*2}(t')\,dt' + {1\over 2}W^*(t) = 
    {1\over 2}W^*(t_0) < {1\over 2}W(0).
  \leqno (5.9)$$

  $$ \left({\rm By\ definition\ }W(0)=\intpi U_i(x)U_i(x)\dee{x}\right).
  \leqno (5.10)$$

  3$^o$) Inequality (5.7) and the inequality $W^*(0)<W(0)$
justify the following proposition

\vfill
\eject

 % f236

\evenheader{236}

  Soit une constante arbitrairement faible
$\eta$ $(0<\eta<W(0))$; 
nommons 
$R_1(\eta)$
la longueur que:

  $$\tripl{r_0>R_1(\eta)} U_i(x)U_i(x)\dee{x} = {\eta\over 2},$$
d\'esignons par $s(\eta,t)$ la sph\`ere, qui d\'epend contin\^ument de
$\eta$ et de $t$, dont le centre est l'origine des coordon\'ees et
dont le rayon est:

  $$R_2(\eta,t)=R_1(\eta,t)+
      {4\over \eta}\left[
                   {W(0)\sqrt{\nu t}\over\sqrt{2}}
               + {W^{3\over 2}t^{1\over 4}\over
                          2^{1\over 4}\pi^{1\over 2}\nu^{3\over 4}}
                   \right],$$
nous avons:

  $$\limsup_{\epsilon\rightarrow 0}
           \tripl{\Pi - s(\eta,t)} u_i^*(x,t)u_i^*(x,t)\dee{x}
            \leq \eta.
  \leqno (5.11)$$

\medskip

  {\bf 29.} Faisons tendre $\epsilon$ vers z\'ero par une suite d\'enombrable
de valeurs:
$\epsilon_1,\epsilon_2,\ldots$ Consid\'erons les fonctions $W^*(t)$ 
qui leur correspondent; elles sont born\'ees dans leur ensemble et chacune
d'elles est d\'ecroissante. Le Proc\'ed\'e diagonal de Cantor (\S 4)
permet d'extraire de la suite 
$\epsilon_1,\epsilon_2,\ldots$
une suite partielle 
$\epsilon_{l_1},\epsilon_{l_2},\ldots$
telle que les fonctions $W^*(t)$ conrrespondantes convergent pour
chaque valeur rationnelle de $t$. Ces fonctions $W^*(t)$ convergent alors
vers une fonction d\'ecroissante, sauf peut-\^etre en des points de
discontinuit\'e de cette derni\`ere.
Les points de discontinuit\'e d'une fonction d\'ecroissante sont 
d\'enombrable. Une seconde application du Proc\'ed\'e diagonal de
Cantor permet donc d'extraire de la suite
$\epsilon_{l_1},\epsilon_{l_2},\ldots$
une suite partielle
$\epsilon_{m_1},\epsilon_{m_2},\ldots$
telle que les fonctions $W^*(t)$ 
correspondants convergent$^1$ quel que soit $t$. Nous nommerons 
$W(t)$
la fonction d\'ecroissante qui est leur limite. (Cette d\'efinition ne
contredit pas (5.10).)

  L'in\'egalit\'e $W^*(t)<W(0)$ prouve que chacune des int\'egrals:

  $$\int_{t_1}^{t_2}\,dt'\tripl{\varpi} u_i^*(x,t')\dee{x};\quad
    \int_{t_1}^{t_2}\,dt'\tripl{\varpi} 
                               \overline{u_k^*(x,t')}u_i^*(x,t')
                           \dee{x} $$
est inf\'erieure \`a une borne ind\'ependente de $\epsilon$. Par une
troisi\`eme application du Proc\'ed\'e diagonal de Cantor nous pouvons
donc extraire de la suite
$\epsilon_{m_1},\epsilon_{m_2},\ldots$
une suite partielle
$\epsilon_{n_1},\epsilon_{n_2},\ldots$
telle que chacune de ces int\'egrals ait une limite

\footrule

  $^1$ En d'autres termes nous utilisons le th\'eor\`eme de Helly.

\rm
\vfill
\eject

 % e236

\evenheader{236}

  Le $\eta$ be an arbitrarily small constant with 
$(0<\eta<W(0))$. We let
$R_1(\eta)$
be the length for which

  $$\tripl{r_0>R_1(\eta)} U_i(x)U_i(x)\dee{x} = {\eta\over 2}$$
and write $s(\eta,t)$ for the sphere with center at the origin with
radius

  $$R_2(\eta,t)=R_1(\eta,t)+
      {4\over \eta}\left[
                   {W(0)\sqrt{\nu t}\over\sqrt{2}}
               + {W^{3\over 2}t^{1\over 4}\over
                          2^{1\over 4}\pi^{1\over 2}\nu^{3\over 4}}
                   \right].$$
We have

  $$\limsup_{\epsilon\rightarrow 0}
           \tripl{\Pi - s(\eta,t)} u_i^*(x,t)u_i^*(x,t)\dee{x}
            \leq \eta.
  \leqno (5.11)$$

\medskip

  {\bf 29.} Let $\epsilon$ tend to zero through a countable  sequence of values
$\epsilon_1,\epsilon_2,\ldots$ Consider the corresponding
functions $W^*(t)$. This is a bounded set of functions and each is 
decreasing.
Cantor's diagonal method (\S 4) allows us to extract from the sequence
$\epsilon_1,\epsilon_2,\ldots$
 a subsequence
$\epsilon_{l_1},\epsilon_{l_2},\ldots$
such that the  $W^*(t)$ converge for all rational values of $t$. The
$W^*(t)$ therefore converge to a decreasing function, except maybe
at points of discontinuity of the limit.
The points discontinuity of a decreasing function are countable.
A second application of Cantor's method allows us to extract from 
$\epsilon_{l_1},\epsilon_{l_2},\ldots$
a subsequence 
$\epsilon_{m_1},\epsilon_{m_2},\ldots$
such that the corresponding $W^*(t)$  converge$^1$ for all $t$. We write 
$W(t)$
for the decreasing function which is their limit. (This definition does
not contradict (5.10).)

  The inequality  $W^*(t)<W(0)$ shows that each of the integrals

  $$\int_{t_1}^{t_2}\,dt'\tripl{\varpi} u_i^*(x,t')\dee{x};\quad
    \int_{t_1}^{t_2}\,dt'\tripl{\varpi} 
                               \overline{u_k^*(x,t')}u_i^*(x,t')
                           \dee{x} $$
is less than a bound independent of $\epsilon$. By a third use of
Cantor's diagonal method we can therefore extract from the sequence
$\epsilon_{m_1},\epsilon_{m_2},\ldots$
a subsequence
$\epsilon_{n_1},\epsilon_{n_2},\ldots$
such that each of these integrals has a unique 

\footrule

  $^1$ In other words we use Helly's theorem.

\rm
\vfill
\eject

 % f237

\foddheader{237}

\noindent unique quand $t_1$ et $t_2$ sont rationnels et que
$\varpi$ est un cube d'ar\^etes parall\`eles aux axes et de
sommets \`a coordon\'ees rationnelles. L'in\'egalit\'e
$W^*(t)<W(0)$ et les hypoth\`eses faite sur les fonctions
$a_i(x,t)$ permettent d'affirmer que dans ces conditions chacune des
int\'egrales:

  $$\int_0^t\,dt'\intpi
          u_i^*(x,t')
          \left[
           \nu\Delta a_i(x,t')+{\partial a_i(x,t')\over\partial t'}
          \right]
               \dee{x};$$

  $$\int_0^t\,dt'\intpi
          \overline{u_k^*(x,t')}
          u_i^*(x,t')
          {\partial a_i(x,t')\over\partial x_k}
               \dee{x}$$
a une limite unique.  Ce r\'esultat,  port\'e  dans 
(5.8),   nous apprend que l'int\'egrale 

  $$\intpi u_i^*(x,t)a_i(x,t)\dee{x}$$
converge vers une limite unique, quels que soient 
$a_i(x,t)$ et $t$. Donc (cf. Corollaire du lemme 7) les fonctions
$u_i^*(x,t)$ convergent faiblement en moyenne vers une limite
$U_i(x,t)$ pour chaque valeur de $t$.

  Ainsi, \'etant donn\'ee une suite de valeurs de $\epsilon$ qui
tendent vers z\'ero, on peut en extraire une suite partielle telle que les
{\it fonctions $W^*(t)$ convergent vers une limite unique} $W(t)$ et que
{\it les fonctions $u_i^*(x,t)$ aient pour chaque valeur de $t$
une faible limite en moyenne unique: } $U_i(x,t)$. Nous supposerons 
d\'esormais que $\epsilon$ tend vers z\'ero par une suite de valeurs
$\epsilon^*$ telle que ces deux circonstances se produisent.

\medskip

  {\it Remarque I } Nous avons d'apr\`es (1.9):

  $$W(t)\geq\intpi U_i(x,t)U_i(x,t)\dee{x}.$$

\footrule

  $^1$ De ces hypoth\`eses r\'esulte en effet qu'\'etant donn\'es $t$,
un nombre $\eta$ $(>0)$ et une fonction 
$\delta(x,t)$ \'egale \`a l'une des d\'eriv\'ees des fonctions 
$a_i(x,t)$ on peut trouver un entier $N$ et deux fonctions
discontinues
$\beta(x,t)$ et
$\gamma(x,t)$
qui possedent les propri\'et\'es suivantes: ces fonctions 
$\beta(x,t)$, 
$\gamma(x,t)$
restent constantes quand 
$x_1$, $x_2$, $x_3$, $t$ varient sans atteindre aucune valeur multiple de
${1\over N}$; chacune d'elles est nulle hors d'un domaine 
$\varpi$; on a:

  $$\int_0^t\,dt'\intpi [\delta(x,t')-\beta(x,t')]^2\dee{x} < \eta;
     \quad |\delta(x,t')-\gamma(x,t')|<\eta \quad {\rm pour}\quad 0<t'<t.$$

\rm
\vfill
\eject

 % e237

\oddheader{237}

\noindent limit when $t_1$ and $t_2$ are rational and
$\varpi$ is a cube with sides parallel to the axes and with vertices
having rational coordinates.
The inequality $W^*(t)<W(0)$ and the hypotheses made on the
$a_i(x,t)$ imply that the integrals 

  $$\int_0^t\,dt'\intpi
          u_i^*(x,t')
          \left[
           \nu\Delta a_i(x,t')+{\partial a_i(x,t')\over\partial t'}
          \right]
               \dee{x};$$

  $$\int_0^t\,dt'\intpi
          \overline{u_k^*(x,t')}
          u_i^*(x,t')
          {\partial a_i(x,t')\over\partial x_k}
               \dee{x}$$
have a unique limit. This result, with (5.8) shows that

  $$\intpi u_i^*(x,t)a_i(x,t)\dee{x}$$
converges to a unique limit, for all
$a_i(x,t)$ and $t$. Therefore (cf. Corollary to lemma 7) the
$u_i^*(x,t)$ converge weakly in mean to some limit 
$U_i(x,t)$ for each value of $t$.

  Also, given a sequence of values of $\epsilon$ which tend to zero,
one can extract from them a subsequence such that 
{\it the  $W^*(t)$ converge to a unique limit} $W(t)$ and that
{\it the  $u_i^*(x,t)$ have for each value of $t$ a unique  weak limit
in mean: } $U_i(x,t)$. 
We suppose from here on that
$\epsilon$ tends to zero through a sequence of values
$\epsilon^*$ such that these two conditions hold.

\medskip

  {\it Remark I } By (1.9)

  $$W(t)\geq\intpi U_i(x,t)U_i(x,t)\dee{x}.$$

\footrule

  $^1$ In fact these hypotheses imply the following.
Given $t$,
a number $\eta$ $(>0)$ and a function 
$\delta(x,t)$ equal to one of the derivatives of the
$a_i(x,t)$, one can find an integer $N$ and two discontinuous functions 
$\beta(x,t)$ and 
$\gamma(x,t)$
with the following properties.
$\beta(x,t)$ 
 and $\gamma(x,t)$
remain constant when 
$x_1$, $x_2$, $x_3$, $t$ vary without hitting(?) any multiple of
${1\over N}$, and each of them is zero outside of a domain
$\varpi$, and(?)

  $$\int_0^t\,dt'\intpi [\delta(x,t')-\beta(x,t')]^2\dee{x} < \eta;
     \quad |\delta(x,t')-\gamma(x,t')|<\eta \quad {\rm for}\quad 0<t'<t.$$

\rm
\vfill
\eject

% f238

\evenheader{238}

  {\it Remarque II. } Le vecteur $U_i(x,t)$ poss\`ede manifestement une
quasi-divergence \'egal \`a z\'ero.

\medskip

  {\bf 30.} L'in\'egalit\'e (5.9) nous donne:

  $$\nu\int_0^\infty
       [\liminf J^*(t')]^2\,dt' < {1\over 2}W(0);$$
la limite inf\'erieure de $J^*(t)$ ne peut donc \^etre $+\infty$
que pour un ensemble de valeurs de $t$ dont la mesure est nulle.
Soit $t_1$ une valeur de l'ensemble compl\'ementaire. On peut extraire
de la suite de valeurs
$\epsilon^*$ envisag\'ee une suite partielle$^1$
$\epsilon^{**}$ telle que sur
$\Pi$ les fonctions
${\partial u_i^{**}(x,t_1)\over\partial x_j}$ correspondantes
convergent faiblement en moyenne vers une limite: 
$U_{i,j}(x,t_1)$ (cf. Th\'eor\`eme fondamental de M. F. Riesz, p. 202).

  Le lemme 2 nous permet d'en d\'eduire tout d'abord que 
{\it les fonctions $U_i(x,t_1)$ ont des quasi-d\'eriv\'ees qui sont ces 
fonctions } $U_{i,j}(x,t_1)$.
Nous poserons:

  $$J(t_1)=\intpi U_{i,j}(x,t_1)U_{i,j}(x,t_1)\dee{x};$$
nous avons (cf. (1.9)):

  $$J(t_1)\leq\liminf J^*(t_1);$$
portons cette in\'egalit\'e dans (5.9); il vient:

  $$\nu\int_{t_0}^t J^2(t')\,dt'+{1\over 2}W(t)
     \leq {1\over 2}W(t_0) \leq {1\over 2}W(0)
     \qquad (0\leq t_0\leq t).
  \leqno (5.12)$$

  Le lemme 2 nous apprend ensuite que sur tout domaine $\varpi$ les
fonctions 
$u_i^{**}(x,t_1)$
convergent fortement en moyenne vers les fonctions
$U_i(x,t)$;

  $$\lim_{\epsilon^{**}\rightarrow 0}
      \tripl{\varpi}
        u_i^{**}(x,t_1)u_i^{**}(x,t_1)
      \dee{x}
     =
     \tripl{\varpi} U_i(x,t_1)U_i(x,t_1)\dee{x}.$$
Choisissons $\varpi$ identique \`a $s(\eta,t_1)$ et tenons compte
de (5.11); il vient:

\footrule

  $^1$ Cette suite partielle que nous choisissons est fonction de
l'\'epoque $t_1$ envisag\'ee.

\rm
\vfill
\eject

 % e238

\evenheader{238}

  {\it Remark II. } The vector $U_i(x,t)$ clearly has quasi-divergence zero.

\medskip

  {\bf 30.} Inequality (5.9) gives us

  $$\nu\int_0^\infty
       [\liminf J^*(t')]^2\,dt' < {1\over 2}W(0).$$
Thus the $\liminf J^*(t)$ can only be $+\infty$
for a set of values of $t$ of measure zero. Suppose $t_1$ is in the
complement of this set. One can extract from the sequence of values
$\epsilon^*$ considered here a subsequence$^1$
$\epsilon^{**}$ such that on $\Pi$ the corresponding functions 
${\partial u_i^{**}(x,t_1)\over\partial x_j}$ 
converge weakly in mean to a limit
$U_{i,j}(x,t_1)$ (cf. Fundamental Theorem of F. Riesz, p. 202).

  Lemma 2 allows us to conclude that 
{\it the $U_i(x,t_1)$ have quasi-derivatives
which are the } $U_{i,j}(x,t_1)$.
We set

  $$J(t_1)=\intpi U_{i,j}(x,t_1)U_{i,j}(x,t_1)\dee{x}.$$
We have (cf. (1.9))

  $$J(t_1)\leq\liminf J^*(t_1).$$
Using this inequality in (5.9) we obtain

  $$\nu\int_{t_0}^t J^2(t')\,dt'+{1\over 2}W(t)
     \leq {1\over 2}W(t_0) \leq {1\over 2}W(0)
     \qquad (0\leq t_0\leq t).
  \leqno (5.12)$$

  Lemma 2 teaches us finally that on all domains $\varpi$ the
$u_i^{**}(x,t_1)$
converge strongly in mean to the 
$U_i(x,t)$;

  $$\lim_{\epsilon^{**}\rightarrow 0}
      \tripl{\varpi}
        u_i^{**}(x,t_1)u_i^{**}(x,t_1)
      \dee{x}
     =
     \tripl{\varpi} U_i(x,t_1)U_i(x,t_1)\dee{x}.$$
Choosing $\varpi$ to be $s(\eta,t_1)$ and taking account of (5.11)
we get

\footrule

  $^1$ The subsequence we choose is a function of $t_1$.

\rm
\vfill
\eject

 % f239

\foddheader{239}

  $$\limsup \intpi u_i^{**}(x,t_1) u_i^{**}(x,t_1) \dee{x}
     \leq
         \tripl{s(\eta,t_1)} U_i(x,t_1) U_i(x,t_1) \dee{x}
     +\eta.$$
D'o\`u, puisque $\eta$ est arbitrairement faible et que 
$W^*(t_1)$
a une valeur limite:

  $$\lim_{\epsilon^*\rightarrow 0}
       \intpi u_i^*(x,t_1) u_i^*(x,t_1) \dee{x}
     \leq
       \intpi U_i(x,t_1) U_i(x,t_1) \dee{x}.
  \leqno (5.13)$$
Appliquons le crit\`ere de forte convergence \'enonc\'e p. 200;
nous constatons que {\it sur $\Pi$ les fonctions $u_i^*(x,t)$
convergent fortement en moyenne vers les fonctions $U_i(x,t)$ 
pours toutes les valeurs $t_1$ de $t$ qui n'appartiennent pas \`a
l'ensemble de mesure nulle sur lequel la limite inf\'erieure de 
$J^*(t)$ est $+\infty$.}

  Pour toutes ces valeurs $t$ les deux membres de (5.13) sont
\'egaux c'est-\`a-dire:

  $$W(t_1)=\intpi  U_i(x,t_1) U_i(x,t_1) \dee{x}.
  \leqno (5.14)$$

  Les fonctions $\overline{u_i^*(x,t_1)}$ elles aussi convergent
fortement en moyenne vers $U_i(x,t_1)$ (cf. G\'en\'eralisation du
lemme 3, p. 207). L'int\'egrale qui figure dans (5.8):

  $$\intpi \overline{u_k^*(x,t')} u_i^*(x,t') 
           {\partial a_i(x,t')\over\partial x_k}
    \dee{x} $$
converge donc vers:

  $$\intpi U_k(x,t') U_i(x,t') 
           {\partial a_i(x,t')\over\partial x_k}
    \dee{x} $$
pour presque toutes les valeurs de $t'$ (cf. (1.8)); cette int\'egrale
est d'autre part inf\'erieure \`a

  $$3W(0)\max |{\partial a_i(x,t')\over\partial x_k}|;$$
le Th\'eor\`eme de M. Lebesgue qui concerne le passage \`a la limite
sous le signe $\int$ permet d'en d\'eduire que:

\vfill
\eject

 % e239

\oddheader{239}

  $$\limsup \intpi u_i^{**}(x,t_1) u_i^{**}(x,t_1) \dee{x}
     \leq
         \tripl{s(\eta,t_1)} U_i(x,t_1) U_i(x,t_1) \dee{x}
     +\eta.$$
From this, since $\eta$ is arbitrarily small and since
$W^*(t_1)$
has a limit

  $$\lim_{\epsilon^*\rightarrow 0}
       \intpi u_i^*(x,t_1) u_i^*(x,t_1) \dee{x}
     \leq
       \intpi U_i(x,t_1) U_i(x,t_1) \dee{x}.
  \leqno (5.13)$$
We apply the strong convergence criterion from p. 200.
Note that  {\it on $\Pi$ the $u_i^*(x,t)$
converge strongly in mean to the $U_i(x,t)$ 
for all values $t_1$ of $t$ not belonging to the 
set of measure zero on which 
$\liminf J^*(t) = +\infty$.}

  For all these values of  $t$ the two sides of (5.13) are equal,
i.e.

  $$W(t_1)=\intpi  U_i(x,t_1) U_i(x,t_1) \dee{x}.
  \leqno (5.14)$$

  The functions $\overline{u_i^*(x,t_1)}$ also converge strongly
in mean to
$U_i(x,t_1)$ (cf. Generalisation of lemma 3, p. 207). The integral
which figures in (5.8)

  $$\intpi \overline{u_k^*(x,t')} u_i^*(x,t') 
           {\partial a_i(x,t')\over\partial x_k}
    \dee{x} $$
therefore converges to

  $$\intpi U_k(x,t') U_i(x,t') 
           {\partial a_i(x,t')\over\partial x_k}
    \dee{x} $$
for almost all values of  $t'$ (cf. (1.8)). Further, this integral is
less than 

  $$3W(0)\max |{\partial a_i(x,t')\over\partial x_k}|$$
Lebesgue's theorem concerning passage to the limit under the $\int$
sign gives

\vfill
\eject

 % f240

\evenheader{240}

  $$\lim_{\epsilon^*\rightarrow 0}
      \int_0^t\,dt'
       \intpi \overline{u_k^*(x,t')} u_i^*(x,t') 
           {\partial a_i(x,t')\over\partial x_k}
       \dee{x} = $$
  
  $$\int_0^t\,dt'
      \intpi  
         U_k(x,t') U_i(x,t')
         {\partial a_i(x,t')\over\partial x_k}
       \dee{x},$$
le second membre de cette relation pouvant \^etre mis, d'apr\`es 
lemme 5, sous la forme:

  $$-\int_0^t\,dt'
      \intpi  
         U_k(x,t') U_{i,k}(x,t')
         a_i(x,t')
       \dee{x}.$$

  D\`es le d\'ebut de ce paragraphe nous avions le droit d'affirmer que les
autres termes qui figurent dans (5.8) convergent de m\^eme
ver des limites qui s'obtiennent en substituant 
$U_i(x,t)$ 
\`a
$u_i^*(x,t)$, $U_i(x)$ \`a $\overline{U_i(x)}$.
Par suite:

  $$\intpi U_i(x,t)a_i(x,t)\dee{x} =
       \intpi U_i(x)a_i(x,0)\dee{x}  $$

  $$ + \int_0^t\,dt'\intpi 
                     U_i(x,t)
                     \left[\nu\Delta a_i(x,t')
                          +{\partial a_i(x,t')\over \partial t'}
                     \right]
                     \dee{x} 
  \leqno (5.15)$$

  $$ - \int_0^t\,dt'\intpi 
                     U_k(x,t') U_{i,k}(x,t')
                     a_i(x,t')
                    \dee{x}.$$
\medskip

  {\bf 31.} Les r\'esultats ainsi obtenus conduisent \`a la d\'efinition
suivante: Nous dirons qu'un vecteur $U_i(x,t)$, d\'efini pour
$t\geq 0$, constitue {\it une solution turbulente des \'equations de
Navier } quand les conditions que nous allons  \'enoncer se trouveront
r\'ealis\'ees, les valeurs de $t$ que nous nommerons {\it singuli\`eres }
constituant un ensemble de mesure nulle:

\medskip
\qquad\hrule
\smallskip
\qquad\hrule
\medskip

  Pour chacque valeur positive de $t$ les fonctions $U_i(x,t)$
sont de carr\'es sommables sur $\Pi$, et le vecteur $U_i(x,t)$ a une
quasi-divergence nulle.

\smallskip

  La fonction:

\vfill
\eject

 % e240

\evenheader{240}

  $$\lim_{\epsilon^*\rightarrow 0}
      \int_0^t\,dt'
       \intpi \overline{u_k^*(x,t')} u_i^*(x,t') 
           {\partial a_i(x,t')\over\partial x_k}
       \dee{x} = $$
  
  $$\int_0^t\,dt'
      \intpi  
         U_k(x,t') U_i(x,t')
         {\partial a_i(x,t')\over\partial x_k}
       \dee{x},$$
By lemma 5 the right hand side of this can be put into the form

  $$-\int_0^t\,dt'
      \intpi  
         U_k(x,t') U_{i,k}(x,t')
         a_i(x,t')
       \dee{x}.$$

  From the beginning of this paragraph we can claim that the other
terms in (5.8) similarly converge. We obtain the limits by 
substituting 
$U_i(x,t)$ 
for $u_i^*(x,t)$ and  $U_i(x)$ for $\overline{U_i(x)}$.
This gives

  $$\intpi U_i(x,t)a_i(x,t)\dee{x} =
       \intpi U_i(x)a_i(x,0)\dee{x}  $$

  $$ + \int_0^t\,dt'\intpi 
                     U_i(x,t)
                     \left[\nu\Delta a_i(x,t')
                          +{\partial a_i(x,t')\over \partial t'}
                     \right]
                     \dee{x} 
  \leqno (5.15)$$

  $$ - \int_0^t\,dt'\intpi 
                     U_k(x,t') U_{i,k}(x,t')
                     a_i(x,t')
                    \dee{x}.$$
\medskip

  {\bf 31.} These results lead to the following definition.
We say that a vector  $U_i(x,t)$ defined for
$t\geq 0$ constitutes {\it a turbulent solution to Navier's equations }
when the following conditions are realised, where values of $t$ that
we call {\it singular } form a set of measure zero.

\medskip
\qquad\hrule
\smallskip
\qquad\hrule
\medskip

  For each positive  $t$ the functions $U_i(x,t)$ are square summable on $\Pi$
and the vector $U_i(x,t)$ has quasi-divergence zero.

\smallskip

  The function

\vfill
\eject

% f241

\foddheader{241}

  $$\int_0^t\,dt'\intpi 
                     U_i(x,t')
                     \left[\nu\Delta a_i(x,t')
                          +{\partial a_i(x,t')\over \partial t'}
                     \right]
                     \dee{x}
     -\intpi U_i(x,t) a_i(x,t) \dee{x}$$

  $$- \int_0^t\,dt'\intpi 
                     U_k(x,t') U_{i,k}(x,t')
                     a_i(x,t')
                    \dee{x}$$
est constante $(t\geq 0)$. (Autrement dit la relation (5.15) a lieu.)
Pour toutes les valeurs positives de $t$, sauf \'eventuellement pour certaines
valeurs {\it singuli\`eres}, les fonctions 
$U_i(x,t)$ 
poss\`edent des quasi-d\'eriv\'ees
$U_{i,j}(x,t)$, de carr\'es sommables sur $\Pi$.

  Nous poserons:

  $$J^2(t)=\intpi U_{i,j}(x,t) U_{i,j}(x,t) \dee{x},$$
$J(t)$ se trouvant donc d\'efini pour presque toutes les valeurs positives
de $t$.

  Il existe une fonction $W(t)$, d\'efinie pour $t\geq 0$, qui poss\`ede
les deux propri\'et\'es suivantes:

  $${\rm la \  fonction \ }\nu\int_0^t J^2(t')\,dt'
        +{1\over 2}W(t)
    {\rm \  est \  non \  croissante};$$
on a: $\intpi U_i(x,t) U_i(x,t) \dee{x} \leq W(t)$,
l'in\'egalit\'e n'ayant lieu qu'\`a certaines \'epoques 
{\it singuli\`eres}, dont l'\'epoque initiale $t=0$ ne fait pas partie.

\medskip
\qquad\hrule
\smallskip
\qquad\hrule
\medskip

  Nous dirons {\it qu'une telle solution
turbulente correspond \`a l'\'etat initial} 
$U_i(x)$ quand nous aurons: $U_i(x,0) = U_i(x)$.

  La conclusion de ce chapitre peut alors se formuler comme suit:

\smallskip

  {\it Th\'eor\`eme d'existence: Supposons donn\'e \`a l'instant initial
un \'etat initial $U_i(x)$ tel que les fonctions $U_i(x)$ soient de
carr\'es sommables sur $\Pi$ et que le vecteur de composantes $U_i(x)$
poss\`ede une quasi-divergence nulle. Il correspond \`a cet \'etat 
initial au moins une solution turbulente, que est d\'efinie pour toutes
les valeurs du temps post\`erieures \`a l'instant initial.}

\bigskip

\centerline{\bf VI. Structure d'une solution turbulente.}

\medskip

  {\bf 32.} Il nous reste \`a \'etablir quel liens existent entre les
solutions r\'eguli\`eres el les solutions turbulentes des \'equations
de Navier. Il est tout d'abord

\vfill
\eject

 % e241

\oddheader{241}

  $$\int_0^t\,dt'\intpi 
                     U_i(x,t')
                     \left[\nu\Delta a_i(x,t')
                          +{\partial a_i(x,t')\over \partial t'}
                     \right]
                     \dee{x}
     -\intpi U_i(x,t) a_i(x,t) \dee{x}$$

  $$- \int_0^t\,dt'\intpi 
                     U_k(x,t') U_{i,k}(x,t')
                     a_i(x,t')
                    \dee{x}$$
is constant $(t\geq 0)$. (Equivalently, (5.15) holds.)
For all positive values of $t$ except possibly for certain
{\it singular} values, the functions
$U_i(x,t)$ 
have quasi-derivatives
$U_{i,j}(x,t)$ which are square summable on $\Pi$.
 
 Set

  $$J^2(t)=\intpi U_{i,j}(x,t) U_{i,j}(x,t) \dee{x},$$
$J(t)$ is thus defined for almost all positive 
$t$.

  There exists a function $W(t)$ defined for $t\geq 0$ which has the
two following properties.

  $${\rm the \ function \ }\nu\int_0^t J^2(t')\,dt'
        +{1\over 2}W(t)
    {\rm \  is \  nonincreasing }$$
and $\intpi U_i(x,t) U_i(x,t) \dee{x} \leq W(t)$, the inequality holding
except for certain {\it singular} times, but $t=0$ is not a singular time.

\medskip
\qquad\hrule
\smallskip
\qquad\hrule
\medskip

  We say that {\it such a turbulent solution corresponds to initial state}
$U_i(x)$ when we have  $U_i(x,0) = U_i(x)$.

  The conclusion of this chapter can then be formulated as follows.

\smallskip

  {\it Existence theorem:  Suppose an initial state $U_i(x)$ is given
such that the functions $U_i(x)$ are square summable on $\Pi$ and that
the vector having components $U_i(x)$ has quasi-divergence zero. There
corresponds to this initial state at least one turbulent solution, which
is defined for all values of $t>0$.}

\bigskip

\centerline{\bf VI. Structure of a turbulent solution.}

\medskip

  {\bf 32.} It remains to establish what connections exist between
regular solutions and turbulent solutions to Navier's equations.
It is entirely

\vfill
\eject

 % f242

\evenheader{242}

\noindent \'evident que toute solution r\'eguli\`ere constitue a
fortiori une solution turbulente. Nous allons chercher dans quels cas une
solution turbulente se trouve constituer une solution r\'eguli\`ere. 
G\'en\'eralisons \`a cet effet les raisonnements du paragraphe 
18 (p. 221).

\medskip

  {\it Comparaison d'une solution r\'eguli\`ere et d'une solution
turbulente:} Soit une solution des \'equations de Navier,
$a_i(x,t)$,
d\'efinie et semi-r\'eguli\`ere pour
$\Theta\leq t<T$; nous supposons qu'elle devient irr\'eguli\`ere quand
$t$ tend vers $T$, \`a moins que $T$ ne soit \'egal \`a
$+\infty$. Consid\'erons une solution turbulente,
$U_i(x,t)$, d\'efinie pour
$\Theta\leq t$, l'\'epoque $\Theta$ n'\'etant pas singuli\`ere.
Les symboles 
$W(t)$ et $J(t)$ se rapporteront \`a cette solution turbulente. 
Nous poserons:

  $$w(t)=W(t)-2\intpi U_i(x,t) a_i(x,t) \dee{x}
        + \intpi a_i(x,t) a_i(x,t) \dee{x}$$

  $$j^2(t)=J^2(t)-2\tripl{}
         U_{i,j}(x,t){\partial a_i(x,t)\over \partial x_j}
                    \dee{x}
          + \intpi
              {\partial a_i(x,t)\over \partial x_j}
              {\partial a_i(x,t)\over \partial x_j}
            \dee{x}.$$
Rappelons que la fonction de $t$:

  $$\nu\int_\Theta^t\,dt'\intpi
              {\partial a_i(x,t')\over \partial x_j}
              {\partial a_i(x,t')\over \partial x_j}
            \dee{x}
          + {1\over 2}\intpi a_i(x,t)a_i(x,t) \dee{x}$$
est constante et que la fonction:

  $$\nu\int_\Theta^t J^2(t')\,dt' + {1\over 2}W(t)$$
est non croissante. Il en r\'esulte que la fonction de $t$:

  $$\nu\int_\Theta^t j^2(t')\,dt' + {1\over 2}w(t)
     + 2\nu\int_\Theta^t \,dt'
       \intpi
         U_{i,k}(x,t'){\partial a_i(x,t')\over \partial x_k}
       \dee{x}
     + 
  \leqno (6.1)$$

  $$\qquad\qquad\qquad\qquad\qquad\intpi U_i(x,t) a_i(x,t) \dee{x}$$
est non croissante. Tenons compte de la relation (5.15) et de ce que
$a_i(x,t)$

\vfill
\eject

 % e242

\evenheader{242}

\noindent clear that any regular solution is a fortiori a turbulent
solution. We are going to look for those cases in which a turbulent
solution is regular. To this end we generalize the reasoning of paragraph
18 (p. 221).

\medskip

  {\it Comparison of a regular solution and a turbulent solution:}
Let 
$a_i(x,t)$ be a solution to Navier's equations, defined and semi-regular for
$\Theta\leq t<T$. We suppose that it becomes irregular when
$t$ tends to $T$, at least in the case when $T$ is not equal to 
$+\infty$. Consider a turbulent solution
$U_i(x,t)$ defined for 
$\Theta\leq t$, where $\Theta$ is not a singular time.
The symbols 
$W(t)$ and $J(t)$ correspond to the turbulent solution. 
Set

  $$w(t)=W(t)-2\intpi U_i(x,t) a_i(x,t) \dee{x}
        + \intpi a_i(x,t) a_i(x,t) \dee{x}$$

  $$j^2(t)=J^2(t)-2\tripl{}
         U_{i,j}(x,t){\partial a_i(x,t)\over \partial x_j}
                    \dee{x}
          + \intpi
              {\partial a_i(x,t)\over \partial x_j}
              {\partial a_i(x,t)\over \partial x_j}
            \dee{x}.$$
Recall that the funtion of $t$

  $$\nu\int_\Theta^t\,dt'\intpi
              {\partial a_i(x,t')\over \partial x_j}
              {\partial a_i(x,t')\over \partial x_j}
            \dee{x}
          + {1\over 2}\intpi a_i(x,t)a_i(x,t) \dee{x}$$
is constant in $t$ and that
the function

  $$\nu\int_\Theta^t J^2(t')\,dt' + {1\over 2}W(t)$$
is nonincreasing. Consequently the function of $t$

  $$\nu\int_\Theta^t j^2(t')\,dt' + {1\over 2}w(t)
     + 2\nu\int_\Theta^t \,dt'
       \intpi
         U_{i,k}(x,t'){\partial a_i(x,t')\over \partial x_k}
       \dee{x}
     + 
  \leqno (6.1)$$

  $$\qquad\qquad\qquad\qquad\qquad\intpi U_i(x,t) a_i(x,t) \dee{x}$$
is nonincreasing. Taking account of relation (5.15) and of that for
$a_i(x,t)$

\vfill
\eject

 % f243

\foddheader{243}

\noindent est une solution semi-r\'eguli\`ere des \'equations de
Navier: nous constatons que la fonction non croissante (6.1) est
\`a une constante pr\`es \'egale \`a la suivante:

  $$\nu\int_\Theta^t j^2(t')\,dt' + {1\over 2}w(t)
     + \int_\Theta^t \,dt'
       \intpi
         [a_k(x,t')-U_k(x,t')]
         U_{i,k}(x,t') a_i(x,t')
       \dee{x}.
  \leqno (6.2)$$
Or nous avons pour chaque valeur non singuli\`ere de $t$:

  $$\intpi [a_k(x,t')-U_k(x,t')]
           {\partial a_i(x,t')\over \partial x_k}
           a_i(x,t')
     \dee{x}=$$

  $${1\over 2}
     \intpi [a_k(x,t')-U_k(x,t')]
           {\partial a_i(x,t') a_i(x,t')\over \partial x_k}
     \dee{x}=0;$$
l'int\'egrale

  $$\intpi [a_k(x,t')-U_k(x,t')]
           U_{i,k}(x,t')
           a_i(x,t')
     \dee{x}$$
peut donc s'\'ecrire

  $$\intpi [a_k(x,t')-U_k(x,t')]
           \left[
             U_{i,k}(x,t')-{\partial a_i(x,t')\over \partial x_k}
           \right]
           a_i(x,t')
     \dee{x}$$
et par suite elle est inf\'erieure en valeur absolue \`a:

  $$\sqrt{w(t')}j(t')V(t'),$$
$V(t')$ d\'esignant la plus grande longueur du vecteur 
$a_i(x,t)$ \`a l'instant $t'$.
Puisque (6.2) n'est pas croissante, il en est donc a fortiori de
m\^eme pour la fonction:

  $$\nu\int_\Theta^t j^2(t')\,dt'+ {1\over 2}w(t) 
      - \int_\Theta^t \sqrt{w(t')}j(t')V(t')\,dt'.$$
Or

  $$\nu\int_\Theta^t j^2(t')\,dt' -
       \int_\Theta^t \sqrt{w(t')}j(t')V(t')\,dt'
      +{1\over 4\nu}\int_\Theta^t w(t')V^2(t')\,dt'$$
ne peut manifestement pas d\'ecroitre. Par suite la fonction:

\vfill
\eject

 % e243

\oddheader{243}

\noindent is a semi-regular solution to Navier's equations.
Note that the nonincreasing function (6.1) is 
up to a constant nearly equal to 

  $$\nu\int_\Theta^t j^2(t')\,dt' + {1\over 2}w(t)
     + \int_\Theta^t \,dt'
       \intpi
         [a_k(x,t')-U_k(x,t')]
         U_{i,k}(x,t') a_i(x,t')
       \dee{x}.
  \leqno (6.2)$$
Now we have for each nonsingular value of $t$

  $$\intpi [a_k(x,t')-U_k(x,t')]
           {\partial a_i(x,t')\over \partial x_k}
           a_i(x,t')
     \dee{x}=$$

  $${1\over 2}
     \intpi [a_k(x,t')-U_k(x,t')]
           {\partial a_i(x,t') a_i(x,t')\over \partial x_k}
     \dee{x}=0.$$
The integral

  $$\intpi [a_k(x,t')-U_k(x,t')]
           U_{i,k}(x,t')
           a_i(x,t')
     \dee{x}$$
may therefore be written

  $$\intpi [a_k(x,t')-U_k(x,t')]
           \left[
             U_{i,k}(x,t')-{\partial a_i(x,t')\over \partial x_k}
           \right]
           a_i(x,t')
     \dee{x}$$
and so it is less in absolute value than

  $$\sqrt{w(t')}j(t')V(t'),$$
where $V(t')$ is the greatest length of the vector
$a_i(x,t)$ at time $t'$.
Since (6.2) is not increasing, it is a fortiori the same for
the function

  $$\nu\int_\Theta^t j^2(t')\,dt'+ {1\over 2}w(t) 
      - \int_\Theta^t \sqrt{w(t')}j(t')V(t')\,dt'.$$
Now

  $$\nu\int_\Theta^t j^2(t')\,dt' -
       \int_\Theta^t \sqrt{w(t')}j(t')V(t')\,dt'
      +{1\over 4\nu}\int_\Theta^t w(t')V^2(t')\,dt'$$
manifestly cannot decrease. It follows that the function

\vfill
\eject

% f244

\evenheader{244}

  $${1\over 2}w(t)-{1\over 4\nu}\int_\Theta^t w(t')V^2(t')\,dt'$$
est non croissante. De l\`a r\'esulte l'in\'egalit\'e qui
g\'en\'eralise (3.7):

  $$w(t)\leq w(\Theta)e^{{1\over 2\nu}\int_\Theta^t V^2(t')\,dt'}
       \qquad (\Theta<t<T).  
  \leqno (6.3)$$

  Supposons en particulier que les solutions $U_i(x,t)$ et
$a_i(x,t)$ correspondent \`a un m\^eme \'etat initial:
$w(\Theta)=0$; de (6.3) r\'esulte alors 
$w(t)=0$;
donc 
$U_i(x,t)=a_i(x,t)$ pour $\Theta\leq t<T$. Ce r\'esultat constitue 
{\it un th\'eor\`eme d'unicit\'e}$^1$ dont ceux des paragraphes 18
et 23 (p. 222 et 228) ne sont que des cas particuliers.

\medskip

  {\bf 33.} {\it R\'egularit\'e d'une solutions turbulente pendant
certains intervalles de temps.}
\smallskip

  Consid\'erons une solution turbulente $U_i(x,t)$, d\'efinie pour
$t\geq 0$. A toute \'epoque non singuli\`ere le vecteur $U_i(x,t)$
constitue un \'etat initial semi-r\'egulier
(cf. p.231 Th\'eor\`eme d'existence, cas a)); le th\'eor\`eme
d'unicit\'e que nous venons d'\'etablir a d\`es lors la 
cons\'e-quence suivante: Soit une \'epoque non 
singuli\`ere, c'est-\`a-dire choisie hors d'un certain
ensemble de mesure nulle: cette \'epoque est l'origine d'un 
intervalle de temps \`a l'int\'erieur duquel la solution 
turbulente envisag\'ee co\"incide avec une solution r\'eguli\`ere 
des \'equations de Navier; et cette co\"incidence ne cesse 
pas tant que cette solution r\'eguli\`ere n'est pas devenue irr\'eguli\`ere.
Ce r\'esultat, compl\'et\'e par quelques autres ais\'es \`a \'etablir,
nous fournit le th\'eor\`eme ci-dessous:

\medskip
\hrule
\medskip
\hrule
\medskip

  {\it Th\'eor\`eme de structure.}
\medskip

  Pour qu'un vecteur $U_i(x,t)$ constitue, quand $t\geq 0$, une
solution turbulente des \'equations de Navier, il fout et il suffit que
ce vecteur poss\`ede les trois propri\'et\'es suivantes:

\smallskip

  a) Nommons {\it intervalle de r\'egularit\'e} tout intervalle
$\overline{\Theta_l T_l}$ de l'axe des temps \`a l'int\'erieur
duquel le vecteur
$U_i(x,t)$ constitue une solution r\'eguli\`ere des \'equations de Navier,
sans que ceci soit vrai pour aucun intervalle contenant
$\overline{\Theta_l T_l}$. Soit $O$ l'ensemble ouvert form\'e par la r\'eunion
de ces intervalles de r\'egularit\'e

\footrule

  $^1$ Je n'ai pu \'etablir de th\'eor\`eme d'unicit\'e affirmant 
qu'\`a un \'etat initial donn\'e correspond une solution turbulent unique.

\rm
\vfill
\eject

 % e244

\evenheader{244}

  $${1\over 2}w(t)-{1\over 4\nu}\int_\Theta^t w(t')V^2(t')\,dt'$$
is nonincreasing. From this we get the inequality generalizing (3.7)

  $$w(t)\leq w(\Theta)e^{{1\over 2\nu}\int_\Theta^t V^2(t')\,dt'}
       \qquad (\Theta<t<T).  
  \leqno (6.3)$$

  Suppose in particular that the solutions $U_i(x,t)$  and
$a_i(x,t)$ correspond to the same initial state. Then 
$w(\Theta)=0$ and by (6.3)  
$w(t)=0$.
therefore
$U_i(x,t)=a_i(x,t)$ for $\Theta\leq t<T$. The uniqueness theorems
of paragraphs 18 and 23 (p. 222 and 228) are special cases of this result.

\medskip

  {\bf 33.} {\it Regularity of a turbulent solution in certain time intervals.}
\smallskip

  Consider a turbulent solution $U_i(x,t)$ defined for
$t\geq 0$. For each nonsingular time
the vector
$U_i(x,t)$
is a semi-regular initial state
(cf. p.231 existence theorem, case a)). The uniqueness theorem that we are
going to establish will have the following consequence.
Consider a nonsingular time, i.e. a time chosen outside of a certain
set of measure zero. Then this is the origin of an interval of time
in the interior of which the turbulent solution coincides with a regular
solution of Navier's equation,
 and this coincidence does not end as long as the regular solution
remains so. This result, complemented by some others
which are easy to establish, gives us the
next theorem.

\medskip
\hrule
\medskip
\hrule
\medskip

  {\it Structure theorem.}
\medskip

  For a vector $U_i(x,t)$ to be a turbulent solution to Navier's equations
for $t\geq 0$, it is necessary and sufficient that it have the following
three properties.

\smallskip

  a) By an {\it interval of regularity} we mean any interval
$\overline{\Theta_l T_l}$ of time in the interior of which
the vector $U_i(x,t)$ is a regular solution to Navier's equations, and such
that this is true for no interval containing $\overline{\Theta_l T_l}$.
Let $O$ be the open set which is the union of the intervals of regularity

\footrule

  $^1$ I have not been able to establish a uniqueness theorem
stating that to a given initial state, there corresponds a unique turbulent
solution.

\rm
\vfill
\eject

 % f245

\foddheader{245}

\noindent (qui sont deux \`a deux sans point commun). $O$ ne doit 
diff\'erer du demi-axe $t\geq 0$ que par un ensemble de mesure
nulle.

\smallskip

  b) La fonction $\intpi U_i(x,t) U_i(x,t)\dee{x}$ est
d\'ecroissante sur l'ensemble que constitue $O$ et l'instant initial
$t=0$.

\smallskip

  c) Quand $t'$ tend vers $t$ les fonctions $U_i(x,t')$ doivent
converger faiblement en moyenne vers les fonctions $U_i(x,t)$.

\smallskip
\hrule
\medskip
\hrule
\medskip

  {\it Compl\'ements:}

\smallskip

  1) Une solution turbulente correspondant \`a un \'etat initial 
semi-r\'egulier coincide avec la solution semi-r\'eguli\`ere qui
correspond \`a cet \'etat initial aussi longtemps que celle-ci existe.

\smallskip

  2) Faisons tendre en croissant $t$ vers l'extr\'emit\'e $T_l$ d'un
intervalle de r\'egularit\'e. La solution $U_i(x,t)$, qui est 
r\'eguli\`ere
pour $\Theta_l < t < T_l$ devient alors irr\'eguli\`ere.

\medskip

   Ce th\'eor\`eme de structure nous permet de {\it r\'esumer notre 
travail} en ces termes: Nous avons essay\'e d'\'etablir l'existence
d'une solution des \'equations de Navier correspondant \`a un \'etat
initial donn\'e: nous n'y avons r\'eussi qu'en renon\c cant \`a la
r\'egularit\'e de la solution en certains instants, convenablement
choisis, dont l'ensemble est de mesure nulle; en ces instants la
solution n'est assujettie qu'\`a une condition de continuit\'e tr\`es 
large (c) et \`a une condition exprimant la non-croissance de 
l'\'energie cin\'etique (b).

\bigskip

  {\it Remarque:} Si le syst\`eme (3.11) poss\`ede une solution non
nulle $U_i(x)$ cette solution permet de construire un exemple tr\`es 
simple de solution turbulente c'est le vecteur $U_i(x,t)$ \'egal \`a

  $$[2\alpha (T-t)]^{-{1\over 2}}
        U_i\left[[2\alpha (T-t)]^{-{1\over 2}}x\right]
    \ {\rm pour} \ t<T$$
et \`a 0 pour $t>T$; il existe une seule \'epoque d'irr\'egularit\'e:
$T$.

\medskip

  {\bf 34.} {\it Compl\'ements relatifs aux intervalles de
r\'egularit\'e et \`a l'allure d'une solution des \'equations de
Navier pour les grandes valeurs du temps.}

  Le chapitre IV nous fournit, outre le th\'eor\`eme d'existence
utilis\'e au paragraphe pr\'ec\'edent, l'in\'egalit\'e (4.3);
d'o\`u r\'esulte la proposition suivante: consid\'erons

\vfill
\eject

 % e245

\oddheader{245}

\noindent (no two of which have a point in common). $O$ must only
differ from
the half axis $t\geq 0$ by a set of measure zero.

\smallskip

  b) The function $\intpi U_i(x,t) U_i(x,t)\dee{x}$ is decreasing
on the set consisting of $O$  and 
$t=0$.

\smallskip

  c) As $t'$ tends to $t$ the $U_i(x,t')$ must converge weakly in mean
to the
$U_i(x,t)$.

\smallskip
\hrule
\medskip
\hrule
\medskip

  {\it Supplementary information}

\smallskip

  1) A turbulent solution corresponding to a semi-regular
initial state coincides with the semi-regular solution having
that initial state, for as long a time as the semi-regular solution exists.

\smallskip

  2) Make $t$ increase to $T_l$ in an interval of regularity. Then the
solution $U_i(x,t)$ which is regular for $\Theta_l < t < T_l$ becomes
irregular. 

\medskip
 
   This structure theorem allows us to 
{\it summarize our work} in these terms:
We have tried to establish the existence of a solution to Navier's 
equations corresponding to a given initial state. We have had to give up
regularity of the solution at a set of times 
of measure zero. At these times
the solution is only subject to a very weak continuity condition (c)
and to  condition (b) expressing the nonincrease of kinetic energy.

\bigskip

  {\it Remark:} If system (3.11) has a nonzero solution $U_i(x)$
then we can
very simply construct a turbulent solution $U_i(x,t)$ equal to

  $$[2\alpha (T-t)]^{-{1\over 2}}
        U_i\left[[2\alpha (T-t)]^{-{1\over 2}}x\right]
    \ {\rm for} \ t<T$$
and to 0 for $t>T$. This has a single irregular time 
$T$.

\medskip

  {\bf 34.} {\it Supplementary information on intervals of regularity
and behavior of solutions to Navier's equations for large time.}

  Chapter IV gives inequality (4.3) in 
addition to the existence theorem used in the previous paragraph.
This results in the following proposition. Consider

\vfill
\eject

 % f246

\evenheader{246}

\noindent une solution turbulente $U_i(x,t)$;
soit une \'epoque non singuli\`ere $t$ et une \'epoque
$T_l$ post\'er-ieure; nous avons:

  $$J(t')>A_1\nu^{3\over 4}(T_l-t')^{-{1\over 4}},$$
$A_1$ \'etant une certaine constante num\'erique. Portons cette minorante
de $J(t')$ dans l'in\'egal-it\'e:

  $$\nu\int_0^{T_l} J^2(t')\,dt'
    < {1\over 2}W(0),$$
il vient:$^*$

  $$2A_1^2\nu^{5\over 2}T_l^{1\over 2}<{1\over 2}W(0).$$

  {\it Toutes les \'epoques singuli\`eres sont donc ant\'erieures \`a
l'\'epoque:}

  $$\theta={W^2(0)\over 16 A_1^4\nu^5}.
  \leqno (6.4)$$
{\it Autrement dit, il existe un intervalle de r\'egularit\'e qui
contient cette \'epoque $\theta$
et qui s'\'etend jusqu'\`a $+\infty$.} Un mouvement r\'egulier
jusqu'\`a l'\'epoque $\theta$ ne devient jamais irr\'egulier.

  Il est ais\'e de pr\'eciser ce r\'esultat:

  Soit un intervalle de r\'egularit\'e de longueur finie
$\overline{\Theta_l T_l}$; toute \'epoque $t$ int\'erieure
\`a cet intervalle est non singuli\`ere; nous avons donc d'apr\`es (4.3):

  $$J(t')>A_1\nu^{3\over 4}(T_l-t')^{-{1\over 4}}
         \ {\rm pour} \ 
        \Theta_l < t' < T_l.$$
\centerline{(cf. Second caract\'ere des irr\'egularit\'es, p. 227.)}

\medskip

  Portons cette minorante de $J(t')$ dans l'in\'egalit\'e:

  $$\nu\sum_l\int_{\Theta_l T_l} J^2(t')\,dt'
    \leq {1\over 2}W(0);$$
il vient, le signe $\sum_l'$ portant sur tous les intervalles de
longueur finie:

  $$2A_1^2\nu^{5\over 2}{\sum_l}'\sqrt{(T_l-\Theta_l)}
     < {1\over 2}W(0).
  \leqno (6.5)$$

\footrule

$^*$ {\tt translator's note:} in the original: 
  $2A_1\nu^{5\over 2}T_l^{1\over 2}<{1\over 2}W(0)$

\vfill
\eject

 % e246

\evenheader{246}

\noindent a turbulent solution $U_i(x,t)$. Let $t'$ be a nonsingular
time$^*$ and $T_l$ a later time. We have

  $$J(t')>A_1\nu^{3\over 4}(T_l-t')^{-{1\over 4}},$$
$A_1$ being a certain numerical constant.
Using this lower bound for $J(t')$ in the inequality

  $$\nu\int_0^{T_l} J^2(t')\,dt'
    < {1\over 2}W(0)$$
we get$^{**}$

  $$2A_1^2\nu^{5\over 2}T_l^{1\over 2}<{1\over 2}W(0).$$

  {\it All singular times occur prior to the time}

  $$\theta={W^2(0)\over 16 A_1^4\nu^5}.
  \leqno (6.4)$$
{\it In other words, there is an interval
of regularity that contains 
$\theta$
and which extends to $+\infty$.}  A motion which is  regular 
up to time $\theta$ never becomes irregular.

  It is easy to make this more precise.

  Let $\overline{\Theta_l T_l}$ be an interval of regularity of
finite length. All times $t$ interior to this interval are nonsingular.
By (4.3)

  $$J(t')>A_1\nu^{3\over 4}(T_l-t')^{-{1\over 4}}
         \ {\rm for} \ 
        \Theta_l < t' < T_l.$$
\centerline{(cf. Second characterisation of irregularity, p. 227.)}

\medskip

  Using this lower bound on $J(t')$ in

  $$\nu\sum_l\int_{\Theta_l T_l} J^2(t')\,dt'
    \leq {1\over 2}W(0).$$
Summing over all the intervals of finite length we get 

  $$2A_1^2\nu^{5\over 2}{\sum_l}'\sqrt{(T_l-\Theta_l)}
     < {1\over 2}W(0).
  \leqno (6.5)$$

\footrule

$^{**}$ {\tt translator's note:} in the original: 
  $2A_1\nu^{5\over 2}T_l^{1\over 2}<{1\over 2}W(0)$

$^*$ {\tt translator's note:} was $t$ in the original

\vfill
\eject

% f247

\foddheader{247}

  Le chapitre IV nous donne, \`a c\^ot\'e de l'in\'egalit\'e
(4.3) (4.4). Il en r\'esulte la propri\'et\'e suivante:
Consid\'erons une solution turbulente, une \'epoque non
singuli\'ere $t'$, une \'epoque post\'erieure $t$. Nous avons:

  $${\rm soit} \ t-t'>A_1^4\nu^3J^{-4}(t'),
      \ {\rm soit} \
       J(t)<(1+A)J(t')$$
c'est-\'`a-dire:$^1$

  $$J(t')>
   \{
     A_1\nu^{3\over 4}(t-t')^{-{1\over 4}};
     {1\over 1+A}J(t)
   \}.$$
Portons cette minorante de $J(t')$ dans l'in\'egalit\'e:

  $$\nu\int_0^t J^2(t')\,dt' \leq {1\over 2}W(0);$$
nous obtenons:

  $$\nu\int_0^t
       \{
     A_1^2\nu^{3\over 2}(t-t')^{-{1\over 2}};
     {1\over (1+A)^2}J^2(t)
      \}
     \,dt'
       \leq {1\over 2}W(0).
  \leqno (6.6)$$
Cette in\'egalit\'e (6.6) fournit pour les valeurs de $t$ sup\'erieures \`a
$\theta$ une majorante de $J(t)$; cette majorante a d'ailleurs une
expression analytique assez compliqu\'ee.

  Nous nous contenterons de remarquer que de (6.6) r\'esulte l'in\'egalit\'e
moins pr\'ecise:

  $$\nu\int_0^t
       \{
     A_1^2\nu^{3\over 2}t^{-{1\over 2}};
     {1\over (1+A)^2}J^2(t)
      \}
     \,dt'
       \leq {1\over 2}W(0);$$
cette derni\`ere exprime tout simplement que

  $$J^2(t) < {(1+A)^2\over 2}{W(0)\over \nu t}
    \ {\rm pour} \
     t>{W^2(0)\over 4 A_1^4\nu^5}.
  \leqno (6.7)$$

  Compl\'etons ce r\'esultat, qui est relatif \`a l'allure asymptotique de
$J(t)$, par un autre relatif \`a l'allure asymptotique de $V(t)$:
Les in\'egalit\'es (4.3) et (4.4) nous apprennent que:

  $$V(t)<AJ(t')[\nu(t-t')]^{-{1\over 4}}
    \ {\rm pour} \
   t-t'<A_1^4\nu^3J^{-4}(t').$$

\footrule

  $^1$ Rappelons que le symbole $\{B; C\}$ nous sert \`a d\'esigner la 
plus petite des quantit\'es $B$ et $C$.

\rm
\vfill
\eject

 % e247

\oddheader{247}

  In Chapter IV  we found the pair of inequalities
(4.3) (4.4). These imply the following.
Consider a turbulent solution, a nonsingular time $t'$,
and a later time $t$. We have

  $${\rm either } \ t-t'>A_1^4\nu^3J^{-4}(t'),
      \ {\rm or } \
       J(t)<(1+A)J(t')$$
in other words$^1$

  $$J(t')>
   \{
     A_1\nu^{3\over 4}(t-t')^{-{1\over 4}};
     {1\over 1+A}J(t)
   \}.$$
Using this lower bound for $J(t')$ in

  $$\nu\int_0^t J^2(t')\,dt' \leq {1\over 2}W(0);$$
we get

  $$\nu\int_0^t
       \{
     A_1^2\nu^{3\over 2}(t-t')^{-{1\over 2}};
     {1\over (1+A)^2}J^2(t)
      \}
     \,dt'
       \leq {1\over 2}W(0).
  \leqno (6.6)$$
This gives an upper bound for $J(t)$ for values of $t$ larger than $\theta$.
However this bound has a rather complicated analytic expression.

  We content ourselves by remarking that (6.6) gives the less precise 

  $$\nu\int_0^t
       \{
     A_1^2\nu^{3\over 2}t^{-{1\over 2}};
     {1\over (1+A)^2}J^2(t)
      \}
     \,dt'
       \leq {1\over 2}W(0).$$
This can most simply be expressed as

  $$J^2(t) < {(1+A)^2\over 2}{W(0)\over \nu t}
    \ {\rm for} \
     t>{W^2(0)\over 4 A_1^4\nu^5}.
  \leqno (6.7)$$

  Complementing this result on asymptotic behavior of $J(t)$ there is
another on $V(t)$. Inequalities (4.3) and (4.4) give

  $$V(t)<AJ(t')[\nu(t-t')]^{-{1\over 4}}
    \ {\rm for} \
   t-t'<A_1^4\nu^3J^{-4}(t').$$

\footrule

  $^1$ Recall that the symbol  $\{B; C\}$ denotes the smaller of $B$ and $C$.

\rm
\vfill
\eject

 % f248

\evenheader{248}

\noindent D'apr\`es (6.7) cette derni\`ere in\'egalit\'e est
satisfaite pour 
$t'={1\over 2}t$ 
quand on prend 
$t>A{W^2(0)\over\nu^5}$; on a donc pour ces valeurs de $t$:

  $$V(t)<A\sqrt{W(0)}(\nu t)^{-{3\over 4}}.$$

  {\it En r\'esum\'e} il existe des constantes $A$ telles que l'on ait:

  $$J(t)<A\sqrt{W(0)}(\nu t)^{-{1\over 2}}
      \ {\rm et} \ V(t)<A\sqrt{W(0)}(\nu t)^{-{3\over 4}}
      \ {\rm pour} \ t>A{W^2(0)\over\nu^5}.$$

  {\it N. B.} J'ignore si $W(t)$ tend n\'ecessairement vers 0 quand
$t$ augmente ind\'efiniment.

\bigskip
\bigskip
\centerline{-------------------------------------}

\vfill
\eject

% fin

 % e248

\evenheader{248}

\noindent By (6.7) this last inequality is satisfied for
$t'={1\over 2}t$ 
if one takes
$t>A{W^2(0)\over\nu^5}$. One therefore has for these $t$

  $$V(t)<A\sqrt{W(0)}(\nu t)^{-{3\over 4}}.$$

  {\it In summary} there exist some constants
$A$  such that

  $$J(t)<A\sqrt{W(0)}(\nu t)^{-{1\over 2}}
      \ {\rm and} \ V(t)<A\sqrt{W(0)}(\nu t)^{-{3\over 4}}
      \ {\rm for} \ t>A{W^2(0)\over\nu^5}.$$

  {\it N. B.} I am ignoring the case in which  $W(t)$ necessarily tends 
to 0 as 
$t$ becomes indefinitely large.

\bigskip
\bigskip
\centerline{-------------------------------------}

\vfill
\eject

% fin

\end